\documentclass[12pt,twoside,a4paper]{report}
\usepackage{amssymb}
\usepackage{amsmath}
\usepackage[dvips]{graphicx}
\usepackage[all,2cell,ps]{xy}
\usepackage{amscd}

\author{Eran Nevo}
 \newtheorem{thm}{Theorem}[section]
 \newtheorem{cor}[thm]{Corollary}
 \newtheorem{lem}[thm]{Lemma}
\newtheorem{de}[thm]{Definition}
\newtheorem{ex}[thm]{Example}
 \newtheorem{prop}[thm]{Proposition}
\newtheorem{obs}[thm]{Observation}
\newtheorem{prob}[thm]{Problem}
\newtheorem{conj}[thm]{Conjecture}
\DeclareMathOperator{\sgn}{sgn} \DeclareMathOperator{\antist}{ast}
\DeclareMathOperator{\kspan}{span_k}
\DeclareMathOperator{\sspan}{span} \DeclareMathOperator{\ddim}{dim}
\DeclareMathOperator{\rank}{rank} \DeclareMathOperator{\Ker}{Ker}
\DeclareMathOperator{\im}{Im} \DeclareMathOperator{\GIN}{GIN}
\DeclareMathOperator{\ddeg}{deg} \DeclareMathOperator{\mmin}{min}
\DeclareMathOperator{\iinit}{init} \DeclareMathOperator{\st}{st}
\DeclareMathOperator{\clst}{\overline{st}}
\DeclareMathOperator{\lk}{lk} \DeclareMathOperator{\Cone}{Cone}
\DeclareMathOperator{\cchar}{char} \DeclareMathOperator{\supp}{supp}
\DeclareMathOperator{\mmod}{mod} \DeclareMathOperator{\ssum}{sum}
\DeclareMathOperator{\Rig}{Rig} \DeclareMathOperator{\soc}{soc}
\DeclareMathOperator{\dom}{dom} \DeclareMathOperator{\res}{res}
\DeclareMathOperator{\ext}{ext}
\DeclareMathOperator{\Stellar}{Stellar}
\DeclareMathOperator{\iif}{if}

\begin{document}

\begin{titlepage}
\begin{center}
\bigskip
\bigskip
\bigskip
\bigskip
\bigskip
\bigskip
\bigskip
\bigskip
\bigskip

\huge\textbf{Algebraic Shifting and $f$-Vector Theory}
\\
\bigskip
\bigskip
\bigskip
\bigskip
\bigskip
\bigskip
\bigskip
\bigskip
\bigskip
\bigskip
\bigskip
\bigskip
\bigskip
\bigskip

\large
Thesis submitted for the degree of
\\
"Doctor of Philosophy"
\\
\bigskip
by
\\
\textbf{Eran Nevo}
\\
\bigskip
\bigskip
\bigskip
\bigskip
\bigskip
\bigskip
\bigskip
\bigskip
\bigskip
\bigskip
\bigskip
\bigskip
\bigskip
\bigskip
\bigskip
\bigskip
\bigskip
\bigskip
\bigskip
\bigskip
\bigskip
Submitted to the Senate of the Hebrew University
\\
\today
\end{center}
\normalsize
\end{titlepage}

\begin{center}
This work was carried out under the supervision of
\\
Prof. Gil Kalai.
\end{center} 
\section*{Abstract}
This manuscript focusses on algebraic shifting and its applications
to $f$-vector theory of simplicial complexes and more general graded
posets. It includes attempts to use algebraic shifting for solving
the $g$-conjecture for simplicial spheres, which is considered by
many as the main open problem in $f$-vector theory. While this goal
has not been achieved, related results of independent interest were
obtained, and are presented here.

The operator algebraic shifting was introduced by Kalai over 20
years ago, with applications mainly in $f$-vector theory. Since
then, connections and applications of this operator to other areas
of mathematics, like algebraic topology and combinatorics, have been
found by different researchers. See Kalai's recent survey
\cite{skira}. We try to find (with partial success) relations
between algebraic shifting and the following  other areas:
\begin{itemize}
\item{}Topological constructions on simplicial complexes.
\item{}Embeddability of simplicial complexes: into spheres and other manifolds.
\item{}$f$-vector theory for simplicial spheres, and more general complexes.
\item{}$f$-vector theory for (non-simplicial) graded partially ordered sets.
\item{}Graph minors.
\end{itemize}
Combinatorially, a (finite) \emph{simplicial complex} is a finite
collection of finite sets which is closed under inclusion. This
basic object has been subjected to extensive research. Its elements
are called \emph{faces}. Its \emph{$f$-vector} $(f_0,f_1,f_2,...)$
counts the number of faces according to their dimension, where $f_i$
is the number of its faces of size $i+1$. $f$-vector theory tries to
characterize the possible $f$-vectors, by means of numerical
relations between the components of the vector, for  interesting
families of simplicial complexes (and more general objects); for
example  for  simplicial complexes which topologically are spheres.

Algebraic shifting associates with each simplicial complex $K$ a
\emph{shifted} simplicial complex, denoted by  $\Delta(K)$, which is
combinatorially simpler. This is an invariant which on the one hand
preserves important invariants of $K$, like its $f$-vector and Betti
numbers, while on the other hand loses other invariants, like the
topological, and even homotopical, type of $K$. A general problem is
to understand which invariants of $K$ can be read off from the faces
of $\Delta(K)$, and how. There are two different variations of this
operator: one is based on the exterior algebra, the other - on the
symmetric algebra; both were introduced by Kalai. Many statements
are true, or conjectured to be true, for both variations.
(Definitions appear in the next chapter.)

The main open problem in $f$-vector theory is to characterize the
$f$-vectors of \emph{simplicial spheres} (i.e. simplicial complexes
which are homeomorphic to spheres). The widely believed conjecture
for the last 25 years, known as the $g$-conjecture, is that the
characterization for simplicial convex polytopes, proved by Stanley
(necessity)\cite{St}, and by Billera and Lee
(sufficiency)\cite{Billera-Lee}, holds also for the wider class of
simplicial spheres, and even for all homology spheres, i.e.
Gorenstein$^*$ complexes. Its open part is to show that the
$g$-vector, which is determined by the $f$-vector, is an
$M$-sequence for these simplicial complexes.

The algebraic properties of face rings hard-Lefschetz and
weak-Lefschetz translate into certain properties of the
symmetrically shifted complex. Having any of these properties in the
face ring of simplicial spheres would imply the $g$-conjecture. A
conjecture by Kalai and by Sarkaria, stating which faces are never
in $\Delta(K)$ if $K$ can be embedded in a sphere, would also imply
the $g$-conjecture for simplicial spheres \cite{skira}. The well
known lower bound and upper bound theorems for $f$-vectors of
simplicial spheres have algebraic shifting conjectured refinements,
which immediately imply these theorems. Details appear in Chapters
\ref{chapter:HL_WL} and \ref{chapter:KalaiSarkaria}. Partial results
on these conjectures include:
\begin{itemize}
\item{}
The known lower bound inequalities for simplicial spheres are shown
to hold for the larger class of doubly Cohen-Macaulay (2-CM)
complexes, by using algebraic shifting / rigidity theory for graphs
and Fogelsanger's theory of \emph{minimal cycle complexes}
\cite{Fogelsanger}. Moreover, the initial part $(g_0,g_1,g_2)$ of
the $g$-vector of a $2$-CM complex (of dimension $\geq 3$) is shown
to be an $M$-sequence. This supports the conjecture by Bj\"{o}rner
and Swartz that the entire $g$-vector of a $2$-CM complex is an
$M$-sequence \cite{Swartz}. See Section \ref{sec:2CM}.
\item{}
The non-negativity of the $g$-vector, which translates to the
generalized lower bound inequalities on the $f$-vector, is proved
for a special class of simplicial spheres, by using special edge
contractions. These contractions are well behaved with respect to
properties of the face rings of those simplicial complexes. To
obtain this result, we first answer affirmatively a problem asked by
Dey et. al. \cite{Dey} concerning topology-preserving edge
contractions in PL-manifolds. See Sections
\ref{sec:TopPreservingEdgeContractions} and \ref{sec:WL&Stellar}.
\item{}
The hard-Lefschetz property is preserved under the combinatorial
operations of join, Stellar subdivisions and connected sum of
spheres; thus supporting the $g$-conjecture. See Sections
\ref{sec:HL&Join}, \ref{sec:WL&Stellar} and \ref{sec:WL,HL&Connected
Sum}.
\end{itemize}
The (generic) rigidity property for a graph mapped into a Euclidean
space of fixed dimension is equivalent to the existence of a certain
edge in the symmetric algebraic shifting of the graph. Similarly,
the dimension of the space of stresses in a generic embedding equals
the number of edges of a certain type in its symmetric shifting.
This follows from a work by Lee \cite{Lee}. Analogues for exterior
shifting involves Kalai's notion of hyperconnectivity \cite{56}. We
use these connections, together with graph minors, to conclude the
following:
\begin{itemize}
\item{} Shifting can tell minors: for every
$2\leq r \leq 6$ and every graph $G$, if $\{r-1,r\} \in \Delta(G)$
then $G$ has a $K_r$ minor. As a corollary, obstructions to
embeddability are obtained. See Section
\ref{sec:ShiftingTellsMinors}.

\item{} Higher dimensional analogues: we generalized the notion of
\emph{minors} in graphs to the class of arbitrary simplicial
complexes, and proved that it 'respects' the Van-Kampen obstruction
in equivariant cohomology. This suggests a new approach for proving
the Kalai-Sarkaria conjecture (and hence the $g$-conjecture).
Details appear in Chapter \ref{chapter:KalaiSarkaria}.
\end{itemize}
Algebraic shifting of more general graded partially ordered sets
than simplicial complexes may be used to prove some of their
properties by looking at the shifted object. For example, an open
problem is to show that the (toric) $g$-vector of a general polytope
is an $M$-sequence. The above approach may be useful in proving it.
Recently Karu has proved that this $g$-vector is non-negative, by
algebraic means. We obtained the following progress in this
direction:
\begin{itemize}
\item{}
We defined an algebraic shifting operator for geometric meet
semi-lattices (simplicial complexes are an important example of
these objects), by constructing face rings for these objects. This
generalizes the construction for simplicial complexes. As an
application, we reprove the fact that their $f$-vector satisfies the
Kruskal-Katona inequalities, proved by Wegner \cite{Wegner}.  See
Section \ref{sec:ShiftGeomLattice}, and the rest of Chapter
\ref{chapter:ShiftingPosets} for further results in the same spirit.
\end{itemize}
Apart from applications, algebraic shifting became an interesting
research object by its own right, as indicated by numerous recent
papers done by multiple researchers.
\begin{itemize}
\item{}
We describe the behavior of algebraic shifting with respect to some
basic constructions on simplicial complexes, like union, cone and
more generally, join. For this, a 'homological' point of view on
algebraic shifting was developed. Interestingly, a multiplicative
formula obtained for exterior shifting of joins, fails for symmetric
shifting. The main applications are as follows; see Chapter
\ref{chapter:ABC} for details.
    \begin{itemize}
    \item{} Proving Kalai's conjecture \cite{skira} that if $K$ and $L$ are disjoint simplicial complexes, then $\Delta(K\cup L)= \Delta(\Delta(K)\cup\Delta(L))$.
    \item{} Disproving Kalai's conjecture for joins \cite{skira}, by providing examples where $\Delta(K*L)$ is not equal to $\Delta(\Delta(K)*\Delta(L))$.
    \item{} A new proof for Kalai's formula for exterior shifting of a cone $\Delta^e(C(K))=C(\Delta^e(K))$.
    \end{itemize}
\end{itemize}
To summarize, the operator algebraic shifting is a powerful tool for proving claims in $f$-vector theory and has relations to the above mentioned areas in mathematics. A better understanding of this operator and these relations may be used to prove well known open problems like the ones indicated here, and is also interesting by its own right. Partial success in achieving this goal is presented in this manuscript. However, it seems that the potential of this tool has not yet been exhausted.\\
\\
Most of the results presented here can be found in papers (except for those in Chapter \ref{chapter:HL_WL}), as follows: most of Chapters \ref{chapter:BasicDef&Concepts} and \ref{chapter:ABC} in \cite{Nevo-ABC}; of Chapter \ref{chapter:rigidity} in \cite{Nevo-2CM} and \cite{Nevo-Stresses}; most of Chapter \ref{chapter:KalaiSarkaria} in \cite{Nevo-VK}; of Chapter \ref{chapter:ShiftingPosets} in \cite{Nevo-GeneralizedMacaulay}.
Each chapter ends with related open problems and conjectures.\\
\\
I hope you will enjoy the reading.

\section*{Acknowledgements}
I would like to express my profound thanks to my advisor prof. Gil
Kalai, for numerous helpful and inspiring mathematical discussions,
as well as for his encouragement, support and care during the years
of work on this thesis.

During those years I had fruitful discussions with various
mathematicians, which reflect directly in this work. Discussions
with Eric Babson about the content of Chapters 4 and 5 led to the
results about join and Stellar subdivisions in Chapter 4; with
Yhonatan Iron to Proposition \ref{prop:Yonatan}; with Carsten
Thomassen to the proofs in Subsection \ref{SubSecMinors}; as well as
helpful discussions with Ed Swartz about the content of Chapter 4.

With many other mathematicians as well I had interesting and
inspiring discussions concerning this work, including Nir Avni, Uri
Bader, Anders Bj{\"o}rner, Nati Linial, Isabella Novik, Vic Reiner,
Uli Wagner and Volkmar Welker.  I am thankful to all of them.

During 2005 I attended the program on algebraic combinatorics held
in Institut Mittag-Leffler. This period was very fruitful and
pleasant for me. I deeply thank the organizers Anders Bj\"{o}rner
and Richard Stanley for the invitation, and my host Svante Linusson
at the Swedish node of the ``Algebraic Combinatorics in Europe''
program, which kindly supported my stay there.

In the last two years I have been supported by a Charles Clore
Foundation Fellowship for Ph.D. students. This fellowship allowed me
to participate in conferences, to devote more time to research and
to study tranquilly. I thank Dame Vivien Duffield and the Clore
Foundation for this support.

I also wish to thank the people of the math department at the Hebrew
University for a pleasant and supportive atmosphere to work in.

Finally, I take this opportunity to thank my family and friends,
Anna and Marni.

\tableofcontents

\chapter{Basic Definitions and Concepts}\label{chapter:BasicDef&Concepts}
\section{Basics of simplicial complexes}
Let $[n]=\{1,2,...,n\}$, and $\binom{[n]}{k}$ denote the subsets of
$[n]$ of size $k$. A collection $K$ of subsets of $[n]$ is called a
(finite abstract) \emph{simplicial complex} if it is closed under
inclusion, i.e. $S\subseteq T\in K$ implies $S\in K$. Note that if
$K$ is not empty (which we will assume from now on) then $\emptyset
\in K$. The $i$-th skeleton of $K$ is $K_{i}=\{S\in K:
|S|=i+1\}=K\cap \binom{[n]}{i+1}$. The elements of $K$ are called
\emph{faces}; those in $K_i$ have \emph{dimension i}. The
$0$-dimensional faces are called \emph{vertices}, the
$1$-dimensional faces are called \emph{edges} and the maximal faces
with respect to inclusion are called \emph{facets}. If all the
facets have the same dimension, $K$ is \emph{pure}. The
$f$-\emph{vector} (face vector) of $K$ is
$f(K)=(f_{-1},f_0,f_1,...)$ where $f_i=|K_i|$. The \emph{dimension
of K} is $\ddim(K):=\max\{i: f_i(K)\neq 0\}$; e.g. a 1-dimensional
simplicial complex is a simple graph. The $f$-\emph{polynomial} of
$K$ is $f(K,t)=\sum_{i\geq 0}f_{i-1}t^i$.

The \emph{link} of a face $S\in K$ is $\rm{lk}(S,K)=\{T\in K: T\cap
S=\emptyset,\ T\cup S\in K\}$. Note that $\rm{lk}(S,K)$ is also a
simplicial complex, and that $\rm{lk}(\emptyset,K)=K$. The (open)
\emph{star} of $S\in K$ is $\st(S,K)=\{T\in K: S\subseteq T\}$,
which is not a simplicial complex; the \emph{closed star} of $S$ in
$K$ is $\clst(S,K)= \{T\in K: T\cup S\in K\}$, which is a simplicial
complex. The \emph{anti star} of $S\in K$ is $\antist(S,K)=\{T\in K:
S\cap T=\emptyset\}$, which is a simplicial complex. The join of two
simplicial complexes $K,L$ with disjoint sets of vertices is the
simplicial complex $K*L=\{S\cup T: S\in K, T\in L\}$. Note that
$f(K*L,t)=f(K,t)f(L,t)$.

A \emph{simplex} in $\mathbb{R}^N$ is the convex hull of some
affinely independent points in $\mathbb{R}^N$. Its intersection with
a supporting hyperplane is a \emph{face} of it, as well as the empty
face. The 0-dimensional simplices are called vertices. A (finite)
\emph{geometric simplicial complex} $L$ in $\mathbb{R}^N$ is a
finite collection of simplices in $\mathbb{R}^N$ such that:
\\(a) Every face of a simplex in $L$ is in $L$.
\\(b) The intersection of any two simplices in $L$ is a face of each of them.
\\We endow the union of simplices in $L$ with the induces topology as a subspace of the Euclidian space $\mathbb{R}^N$,
and call it the \emph{topology of $L$}. A \emph{geometric
realization} of an abstract simplicial complex $K$ is a geometric
simplicial complex $L$ in some $\mathbb{R}^N$ such that $L$ is
combinatorially isomorphic to $K$, i.e. as posets w.r.t. inclusion.
Any two geometric realizations of the same simplicial complex, $K$,
are (piecewise linearly) homeomorphic, and we denote this
topological space by $||K||$. We refer to topological properties of
$||K||$ as properties of $K$; e.g. $K_5$, the complete graph on 5
vertices, is not embeddable in $\mathbb{R}^2$. We say that a
simplicial complex $K$ is a \emph{triangulation} of a topological
space $X$ if $||K||$ is homeomorphic to $X$.

Let $\tilde{H}_i(K;k)$ denote the reduced $i$-th (simplicial)
homology group of $K$ with field $k$ coefficients.
$\beta_i=\beta_i(K;k)=\ddim_k(\tilde{H}_i(K,k))$ is the $i$-th
\emph{Betti number} of $K$ with $k$ coefficients, and
$\beta(K;k)=(\beta_0,\beta_1,...)$ is its \emph{Betti vector}.

Let $<$ denote the usual order on the natural numbers. A simplicial
complex $K$ with vertices $[n]$ is \emph{shifted} if for every
$i<j$, $j\in S\in K$, also $(S\setminus\{j\})\cup \{i\}\in K$. Let
$<_P$ be the product partial order on equal sized ordered subsets of
$\mathbb{N}$. That is, for $S=\{s_1<...<s_i\}$ and
$T=\{t_1<...<t_i\}$ $S\leq_P T$ iff $s_j\leq t_j$ for every $1\leq
j\leq i$. Then $K$ is shifted iff $S<_PT\in K$ implies $S\in K$. It
is easy to see that every simplicial complex has a shifted
simplicial complex with the same $f$-vector: for some $1\leq i<j\leq
n$ apply $S\mapsto (S\setminus\{j\})\cup \{i\}$ for all $j\in S\in
K, i\notin S$ to obtain $K'$, which is also a simplicial complex,
with the same f-vector. Repeat this process as long as possible, to
obtain a shifted simplicial complex $\Delta^c(K)$. Note that
$\Delta^c(K)$ depends on the order of choices of pairs $i<j$. The
operation $K\mapsto \Delta^c(K)$ is called \emph{combinatorial
shifting}, introduced by Erd\"{s} Ko and Rado \cite{Erdos-Ko-Rado}.

Note that a shifted complex $K$ is homotopy equivalent to a wedge of
spheres, where the number of $i$-dimensional spheres in this wedge
is $|\{S\in K_i: S\cup\{1\}\notin K \}|$. In particular, its Betti
numbers are easily read off from its combinatorics:
$\beta_i(K;k)=|\{S\in K_i: S\cup\{1\}\notin K\}|$ for every field
$k$.

For further details about simplicial complexes, and about simplicial
homology, we refer to Munkres' book \cite{Munkres}.

Another useful way to encode the information in the $f$-vector
$f(K)$ is by the \emph{h-vector}. Let $K$ be $(d-1)$-dimensional,
and define
$$\sum_{0\leq i\leq d}h_i(K)x^{d-i}= \sum_{0\leq i\leq d}f_{i-1}(K)(x-1)^{d-i}.$$
Equivalently, $h_k=\sum_{0\leq i\leq
k}(-1)^{k-i}\binom{d-i}{k-i}f_{i-1}$ and $f_{k-1}=\sum_{0\leq i\leq
k}\binom{d-i}{k-i}h_{i}$. The $h$-\emph{polynomial} of $K$ is
$h(K,t)=\sum_{i\geq 0}h_i(K)t^i$, hence
$(1+t)^dh(K,\frac{1}{t+1})=t^df(K,\frac{1}{t})$. If $||K||$ is a
sphere, then $h(K)$ is symmetric, i.e. $h_i(K)=h_{d-i}(K)$ for every
$0\leq i\leq d$. Equivalently, $h(K,t)=t^dh(K,\frac{1}{t})$. These
relations are known as Dehn-Sommerville equations. In this case
$f(K)$ can be recovered from the \emph{g-vector}:
$g_0(K):=h_0(K)=1$, $g_i(K):=h_i(K)-h_{i-1}(K)$ for $1\leq i\leq
\lfloor d/2\rfloor$. $g(K):=(g_0(K),...,g_{\lfloor d/2\rfloor}(K))$.

For more information about the importance of the $h$- and
$g$-vectors, we refer to Stanley's book \cite{St}, and to Chapter
\ref{chapter:HL_WL} below.

\section{Basics of algebraic shifting}\label{sec:BasicsShifting}
We give the definition of exterior and symmetric algebraic shifting
and state some of there basic properties. We develop a 'dual' point
of view which leads to equivalent definitions, and use these 'dual'
definitions to cut the faces of the shifted complex into 'intervals'
which play a crucial role in the proofs of the results of Chapter
\ref{chapter:ABC}.

\subsection{Exterior shifting}
\subsubsection{Via the exterior algebra}
Let $\mathbb{F}$ be a field and let $k$ be a field extension of
$\mathbb{F}$ of transcendental degree $\geq n^2$ (e.g.
$\mathbb{F}=\mathbb{Q}$ and $k=\mathbb{R}$, or
$\mathbb{F}=\mathbb{Z}_2$ and $k=\mathbb{Z}_2(x_{ij})_{1\leq i,j\leq
n}$ where $x_{ij}$ are intermediates). Let $V$ be an $n$-dimensional
vector space over $k$ with basis $\{e_{1},\dots,e_{n}\}$. Let
$\bigwedge V$ be the graded exterior algebra over $V$. Denote
$e_{S}=e_{s_{1}}\wedge\dots\wedge e_{s_{j}}$ where $S=
\{s_{1}<\dots<s_{j}\}$. Then $\{e_{S}: S\in (_{\ j}^{[n]})\}$ is a
basis for $\bigwedge^{j} V$. Note that as $K$ is a simplicial
complex, the ideal $(e_{S}:S\notin K)$ of $\bigwedge V$ and the
vector subspace $\rm{span}\{e_{S}:S\notin K\}$ of $\bigwedge V$
consist of the same set of elements in $\bigwedge V$. Define the
exterior algebra of $K$ by
$$\bigwedge K=(\bigwedge V)/(e_{S}:S\notin K).$$ Let
$\{f_{1},\dots,f_{n}\}$ be a basis of $V$, generic over $\mathbb{F}$
with respect to $\{e_{1},\dots,e_{n}\}$, which means that the
entries of the corresponding transition matrix $A$ ($e_{i}A=f_{i}$
for all $i$) are algebraically independent over $\mathbb{F}$. Let
$\tilde{f}_{S}$ be the image of $f_{S}\in \bigwedge V$ in
$\bigwedge(K)$. Let $<_L$ be the lexicographic order on equal sized
subsets of $\mathbb{N}$, i.e. $S<_LT$ iff $\mmin(S\triangle T)\in
S$. Define
$$\Delta^e(K)=\Delta^e_{A}(K)=\{S: \tilde{f}_{S}\notin \rm{span}\{\tilde{f}_{S'}:S'<_{L}S\}\}$$
to be the shifted complex, introduced by Kalai \cite{55}.
 The construction is canonical, i.e. it is independent of the choice of the generic matrix
 $A$, and for a permutation $\pi:[n]\rightarrow [n]$ the induced simplicial complex
  $\pi(K)$ satisfies $\Delta^e(\pi(K))=\Delta^e(K)$.
 It results in a shifted simplicial complex,
 having the same face vector and Betti vector as $K$'s \cite{BK}.
 Some more of its basic properties are detailed in subsection \ref{subsec:basic properties}.

\subsubsection{Via dual setting}
Fixing the basis $\{e_{1},\dots,e_{n}\}$ of $V$ induces the basis
$\{e_{S}: S\subseteq [n]\}$ of $\bigwedge V$ as a $k$-vector space,
which in turn induces the dual basis $\{e_{T}^{*}: T\subseteq [n]\}$
of $(\bigwedge V)^{*}$ by defining $e_{T}^{*}(e_{S})=\delta_{T,S}$
and extending bilinearly. $(\bigwedge V)^{*}$ stands for the space
of $k$-linear functionals on $\bigwedge V$. For $f,g\in \bigwedge V$
$<f,g>$ will denote $f^{*}(g)$. Define the so called left interior
product of $g$ on $f$ \cite{56}, where $g,f \in \wedge V$, denoted
$g\lfloor f$, by the requirement that for all $h\in\bigwedge V$
$$<h,g\lfloor f>=<h\wedge g,f>.$$
$g\lfloor\cdot$ is the adjoint operator of $\cdot\wedge g$ w.r.t.
the bilinear form $<\cdot,\cdot>$ on $\bigwedge V$. Thus,
$g\lfloor f$ is a bilinear function, satisfying
\begin{equation}\label{floorE}
e_{T}\lfloor e_{S}=\{^{(-1)^{a(T,S)} e_{S\backslash T}\  if\
T\subseteq S} _{0 \ otherwise}
\end{equation}
where $a(T,S)=|\{(s,t)\in S\times T: s\notin T, t<s\}|$. This
implies in particular that for a monomial $g$ (i.e. $g$ is a wedge
product of elements of degree 1) $g \lfloor$ is a boundary
operation on $\bigwedge V$, and in particular on $\kspan\{e_{S}:
S\in K\}$ \cite{56}.
Let $\bigwedge^{j}(K)=\kspan\{e_{S}: S\in
K_{j-1}\}$ and $\bigwedge (K)=\kspan\{e_{S}: S\in K\}$.
We denote:
$$\Ker_{j}f_{R}\lfloor(K) = \Ker_{j}f_{R}\lfloor=\Ker (f_{R}\lfloor: \bigwedge^{j+1}(K) \rightarrow \bigwedge^{j+1-|R|}(K)).$$
Note that the definition of $\bigwedge (K)$ makes sense more
generally when $K_{0}\subseteq [n]$ (and not merely when
$K_{0}=[n]$), and still $f_{R}\lfloor$ operates on the subspace
$\bigwedge (K)$ of $\bigwedge V$ for every $R\subseteq [n]$. (Recall
that $f_{i}=\alpha_{i1}e_{1}+...+\alpha_{in}e_{n}$ where
$A=(\alpha_{ij})_{1\leq i,j\leq n}$ is a generic matrix.) Define
$f^{0}_{i}=\sum_{j\in K_{0}}\alpha_{ij}e_{j}$, and
$f^{0}_{S}=f^{0}_{s_{1}}\wedge\dots\wedge f^{0}_{s_{j}}$ where $S=
\{s_{1}<\dots<s_{j}\}$. By equation (\ref{floorE}), the following
equality of operators on $\bigwedge (K)$ holds:
\begin{equation} \label{floorF}
\forall S\subseteq [n] \ \ f_{S}\lfloor=f^{0}_{S}\lfloor.
\end{equation}

Now we give some equivalent descriptions of the
exterior shifting operator, using the kernels defined above.
This approach will be used in Chapter \ref{chapter:ABC}.
The following generalizes a result for graphs \cite{56} (the
proof is similar):
\begin{prop} \label{prop.1}
Let $R\subseteq [n]$, $|R|<j+1$. Then: $$\Ker_{j}f_{R}\lfloor =
\bigcap_{i:[n]\ni i\notin R}\Ker_{j}f_{R\cup{i}}\lfloor.$$
\end{prop}
$Proof$: Recall that $h \lfloor(g\lfloor f)=(h\wedge g)\lfloor f$.
Thus, if $f_{R}\lfloor m=0$ then
$$f_{R\cup i}\lfloor m=\pm f_{i}\lfloor (f_{R}\lfloor
m)=f_{i}\lfloor 0=0.$$ Now suppose $m\in \bigwedge^{j+1}(K)\setminus \Ker_{j}f_{R}\lfloor$. The set $\{f_{Q}^{*}:
Q\subseteq[n],|Q|=j+1-|R|\}$ forms a basis of $(\bigwedge
^{j+1-|R|}V)^{*}$,
 so there
is some $f_{R'}$, $R'\subseteq[n], |R'|=j+1-|R|$, (note that
$R'\neq \emptyset$) such that
$$<f_{R'}\wedge f_{R},m>=<f_{R'},f_{R}\lfloor m>\neq0.$$
We get that for $i_{0}\in R'$: $i_{0}\notin R$ and $<f_{R'\setminus i_{0}},f_{R\cup
i_{0}}\lfloor m>\neq 0 $. Thus $m\notin \Ker_{j}f_{R\cup i_{0}}\lfloor$ which completes
the proof.$\square$

In the next proposition we determine the shifting of a simplicial
complex by looking at the intersection of kernels of boundary
operations (actually only at their dimensions): Let $S$ be a
subset of $[n]$ of size $s$. For $R\subseteq [n]$,$|R|=s$, we look
at $f_{R}\lfloor:\bigwedge^{s}(K)\rightarrow
\bigwedge^{s-|R|}(K)=k$.
\begin{prop} \label{prop.2}
  Let $K_{0}, S\subseteq[n], |K_{0}|=k, |S|=s$. The following
  quantities are equal:
\begin{eqnarray}\label{KerShift}
\ddim \bigcap_{R<_{L}S, |R|=s, R\subseteq [n]}\Ker_{s-1}f_{R}\lfloor,\label{KS1}\\
\ddim \bigcap_{R<_{L}S, |R|=s, R\subseteq [k]}\Ker_{s-1}f^{0}_{R}\lfloor,\label{KS2}\\
|\{T\in \Delta^e(K): |T|=s, S\leq_{L} T\}|.\label{KS3}
\end{eqnarray}
In particular, $S\in \Delta^e(K)$ iff $$\ddim\bigcap_{R<_{L}S, |R|=s,
R\subseteq [n]}\Ker_{s-1}f_{R}\lfloor > \ddim\bigcap_{R\leq_{L}S,
|R|=s, R\subseteq [n]}\Ker_{s-1}f_{R}\lfloor$$ (equivalently, $S\in
\Delta^e(K)$\ iff\ \  $\bigcap_{R<_{L}S, |R|=s, R\subseteq
[n]}\Ker_{s-1}f_{R}\lfloor \nsubseteq \Ker_{s-1}f_{S}\lfloor$).
\end{prop}
$Proof$:
First we show that (\ref{KS1}) equals (\ref{KS2}). For every
$T\subseteq [n]$, $T<_{L}S$, decompose $T=T_{1}\cup T_{2}$,
 where $T_{1}\subseteq [k], T_{2}\cap [k]=\emptyset$. For each
 $T_{3}$ satisfying $T_{3}\subseteq [k], T_{3}\supseteq T_{1},
 |T_{3}|=|T|$,
 we have $T_{3}\leq_{L}T$. Each $f^{0}_{t}$, where $t\in T_{2}$, is a linear
 combination of the $f^{0}_{i}$'s, $1\leq i \leq k$, so $f^{0}_{T}$ is a
 linear combination of such $f^{0}_{T_{3}}$'s. Thus, for every
 $j\geq s-1$,
 $$\bigcap_{R\leq_{L}T,R\subseteq
 [k]}\Ker_{j}f^{0}_{R}\lfloor \subseteq \Ker_{j}f^{0}_{T}\lfloor,$$ and
 hence $$\bigcap_{R<_{L}S,R\subseteq
 [k]}\Ker_{j}f^{0}_{R}\lfloor = \bigcap_{R<_{L}S,R\subseteq
 [n]}\Ker_{j}f^{0}_{R}\lfloor.$$ Combining with (\ref{floorF}) the desired equality follows.

Next we show that (\ref{KS3}) equals (\ref{KS2}). Let $m\in
\bigwedge^{s}(K)$ and $R\subseteq[n], |R|=s$. Let us express $m$
and $f_{R}$ in the basis $\{e_{S}: S\subseteq [n]\}$:
$$m=\sum_{T\in K,|T|=s} \gamma_{T}e_{T},\ f_{R}=\sum_{S'\subseteq[n],|S'|=s}A_{RS'}e_{S'}$$
where $A_{RS'}$ is the minor of $A$ (transition matrix) with
respect to the rows $R$ and columns $S'$, and where $\gamma_{T}$
is a scalar in $\mathbb{R}$.

By bilinearity we get $$f^{0}_{R}\lfloor m =f_{R}\lfloor m
=\sum_{T\in K,|T|=s}\gamma_{T}A_{RT}.$$ Thus
(\ref{KS2}) equals the dimension of the solution space of the
system $B_{S}x=0$, where $B_{S}$ is the matrix $(A_{RT})$, where
$R<_{L}S$, $R\subseteq [k]$, $|R|=s$ and $T\in K,|T|=s$. But,
since the row indices of $B_{S}$ are an initial set with respect
to the lexicographic order, the intersection of $\Delta^e(K)$ with
this set of indices determines a basis of the row space of
$B_{S}$. Thus, $\rank(B_{S})=|\{R\in \Delta^e(K): |R|=s,R<_{L}S\}|$.
But $K$ and $\Delta^e(K)$ have the same $f$-vector, so we get:
$$\ddim \bigcap_{R<_{L}S, R\subseteq [k], |R|=s}\Ker_{s-1}f^{0}_{R}\lfloor =
f_{s-1}(K) - \rank(B_{S})=$$ $$= |\{T\in \Delta^e(K):|T|=s, S\leq_{L}
T\}|$$ as desired.$\square$

\subsubsection{Dividing $\Delta^e(K)$ into intervals}
For each $j>0$ and $S\subseteq [n]$, $|S|\geq j$ we define
$\iinit_{j}(S)$ to be the set of lexicographically least $j$
  elements in $S$, and for every $i>0$ define
 $$I_{S}^{i} = I_{S}^{i}(n) = \{T: T\subseteq [n], |T|=|S|+i, \iinit_{|S|}(T)=S\}.$$ Let
  $S_{(i)}^{(m)} = S_{(i)}^{(m)}(n)=\min_{<_{L}}I_{S}^{i}(n)$ and
 $S_{(i)}^{(M)} = S_{(i)}^{(M)}(n)=\max_{<_{L}}I_{S}^{i}(n)$.
In the sequel, all the sets of numbers we consider are subsets of
$[n]$. In order to simplify notation, we will often omit noting
that. We get the following information about the partition of the
faces in the shifted complex into 'intervals':
\begin{prop} \label{prop.*}
Let $K_{0}\subseteq [n]$, $S\subseteq [n]$, $i>0$. Then:
$$|I_{S}^{i}\cap \Delta^e(K)| =  \ddim \bigcap_{R<_{L}S}\Ker_{|S|+i-1}f_{R}\lfloor(K) -
\ddim \bigcap_{R\leq_{L}S}\Ker_{|S|+i-1}f_{R}\lfloor(K).$$
\end{prop}
$Proof$:
   By Proposition \ref{prop.1},$$\ddim \bigcap_{R<_{L}S}\Ker_{|S|+i-1}f_{R}\lfloor=
   \ddim \bigcap_{R<_{L}S}\bigcap_{j\notin R, j\in [n]}\Ker_{|S|+i-1}f_{R\cup j}\lfloor =...=$$
$$\ddim \bigcap_{R<_{L}S}\ \bigcap_{T:T\cap R=\emptyset, |T|=i}\Ker_{|S|+i-1}f_{R\cup
T}\lfloor =\ddim \bigcap_{R<_{L}S_{(i)}^{(m)}}\Ker_{|S_{(i)}^{(m)}|-1}f_{R}\lfloor.$$
To see that the last equation is true, one needs to check that
$\{R\cup T: T\cap R=\emptyset, |T|=i, R<_{L}S\} = \{Q:
Q<_{L}S_{(i)}^{(m)}\}$. By Proposition \ref{prop.2},$$\ddim
\bigcap_{R<_{L}S_{(i)}^{(m)}}\Ker_{|S_{(i)}^{(m)}|-1}f_{R}\lfloor =
|\{Q\in\Delta^e(K):|Q|=|S|+i,S_{(i)}^{(m)}\leq_{L}Q\}|.$$ Similarly,
$$\ddim \bigcap_{R\leq_{L}S}\Ker_{|S|+i-1}f_{R}\lfloor(K)=|\{F\in\Delta^e(K):|F|=|S|+i,S_{(i)}^{(M)}<_{L}F\}|.$$
Here one checks
that $\{R\cup T: T\cap R=\emptyset, |T|=i, R\leq_{L}S\} = \{F:
F\leq_{L}S_{(i)}^{(M)}\}$. Thus, the proof of the proposition is
completed. $\square$

Note that on $I_{S}^{1}$ the lexicographic order and the partial
order $<_P$ coincide, since all sets in $I_{S}^{1}$ have the same
$|S|$ least elements. As $\Delta^e(K)$ is shifted,
$I_{S}^{1}\cap\Delta^e(K)$ is an initial set of $I_{S}^{1}$ with
respect to $<_{L}$. Denote for short
$$D(S)=D_{K}(S)=|I^{1}_{\iinit_{|S|-1}(S)}(n)\cap \Delta^e(K)|.$$
$D_{K}(S)$ is indeed independent of the particular $n$ we choose,
as long as $K_{0}\subseteq [n]$. We observe that:
\begin{prop} \label{prop.3}
Let $K_{0}$ and $S=\{ s_{1}<\dots<s_{j}<s_{j+1}\}$ be subsets of
$[n]$. Then: $S\in \Delta^e(K)\Leftrightarrow s_{j+1}-s_{j}\leq
D(S).$ $\square$
\end{prop}
Another easy preparatory lemma is the following:
\begin{prop} \label{prop.5}
Let $K_{0}, S\subseteq [n]$. Then: $D_{\Delta^e(K)}(S)=D_{K}(S).$
\end{prop}
$Proof$: It follows from the fact that $\Delta^e\circ \Delta^e=\Delta^e$ (Kalai
\cite{Kalai-SymmMatroids}, or later on here in Corollary \ref{prop.13
delta^2}). $\square$

\subsection{Symmetric shifting}

\subsubsection{Via the symmetric algebra}
For symmetric shifting, let us look on the face ring
(Stanley-Reisner ring) of $K$ $k[K]=k[x_{1},..,x_{n}]/I_{K}$ where
$I_{K}$ is the homogenous ideal generated by the monomials whose
support is not in $K$, $\{\prod_{i\in S}x_i:\ S\notin K\}$. $k[K]$
is graded by degree. Let $\mathbb{F}\subseteq k$ be fields as before
and let $y_{1},\dots,y_{n}$ be generic linear combinations of
$x_{1},\dots,x_{n}$ w.r.t. $\mathbb{F}$. We choose a basis for each
graded component of $k[K]$, up to degree $\ddim(K)+1$, from the
canonic projection of the monomials in the $y_{i}$'s on $k[K]$, in
the greedy way:
$$\GIN (K)=\{m: \tilde{m}\notin \kspan   \{\tilde{m'}:\ddeg (m')=\ddeg (m), m'<_{L}m\}\}$$
where $\prod y_{i}^{a_{i}}<_{L}\prod y_{i}^{b_{i}}$ iff for $j=\mmin
\{i: a_{i}\neq b_{i}\}$ $a_{j}>b_{j}$). The combinatorial
information in $\GIN (K)$ is redundant: if $m\in \GIN (K)$ is of
degree $i\leq \ddim (K)$ then $y_{1}m,..,y_{i}m$ are also in $\GIN
(K)$. Thus, $\GIN (K)$ can be reconstructed from its monomials of
the form $m=y_{i_{1}}\cdot y_{i_{2}}\cdot..\cdot y_{i_{r}}$ where
$r\leq i_{1}\leq i_{2}\leq..\leq i_{r}$, $r\leq \ddim(K)+1$. Denote
this set by $\rm{gin} (K)$, and define
$S(m)=\{i_{1}-r+1,i_{2}-r+2,..,i_{r}\}$ for such $m$. The collection
of sets
$$\Delta^{s}(K)=\cup \{S(m): m\in \rm{gin}(K)\}$$
carries the same combinatorial information as $\GIN (K)$. $\Delta^{s}(K)$ is a
simplicial complex. Again, the construction is canonic, in
the same sense as for exterior shifting.
If $k$ has characteristic zero then $\Delta^{s}(K)$ is shifted \cite{Kalai-SymmMatroids}.
Further basic properties of $\Delta^{s}(K)$ are detailed in subsection \ref{subsec:basic properties}.

\subsubsection{Via dual setting}\label{subsec:dual symm}
Denote monomials in the graded polynomial ring
$R=k[x_{1},..,x_{n}]=k\oplus R_1\oplus R_2\oplus...$ by
$x^{a}=\prod_{1\leq i\leq n}x_i^{a_i}$ where $a_i\in
\mathbb{Z}_+=\{0,1,2,...\}$, and define a bilinear form on $R$ by
$<x^a,x^b>=\delta_{a,b}$. For a subspace $A$ of $R$ denote its
orthogonal subspace by $A^{\bot}$. Every element $m\in R$ defines a
map on $R$ by multiplication $m:\ r\mapsto mr$, and denote its
adjoint map by $m^*$. This induces a bilinear map on $R$, $m\cdot
u=m^*(u)$ which satisfies $$x^a\cdot x^b=\{^{x^{b-a}\ \iif \ a\leq
b}_{0\ \rm{otherwise}}$$ where $a\leq b$ means that componentwise
$a_i\leq b_i$. Thus, for a simplicial complex $K$ with $K_0\subseteq
[n]$, the restriction  $m^*_{|I_K^{\bot}}$ is into $I_K^{\bot}$, as
a basis for this subspace is $\{x^a: \supp(a)\in K\}$ where
$\supp(a)=\{i: a_i>0\}$. Denote the degree $i$ part of $I_K^{\bot}$
by $(I_K^{\bot})_i$, and the degree of an element $m$ by $\ddeg
(m)$. For a homogenous element $m$, let $\Ker_j(m^*)=\Ker(m^*:\
(I_K^{\bot})_{j+1}\rightarrow (I_K^{\bot})_{j+1-\ddeg (m)})$.

Let $Y=\{y_1,...,y_n\}$ be a generic basis for $R_1$ w.r.t. the basis $X=\{x_1,...,x_n\}$.
\begin{prop} \label{prop.1s}
Let $m$ be a monomial in the basis $Y$ (i.e. $m=\prod_i y_i^{a_i}$), and $j\geq \ddeg (m)$. Then $\Ker_{j}m^* =
\bigcap_{i\in[n]}\Ker_{j}(my_i)^*$.
\end{prop}
$Proof$: By the associativity and commutativity of multiplication in
$R$, when passing to adjoint maps one obtains $(my_i)^*=y_i^*\circ
m^*$, and hence $\Ker_{j}m^* \subseteq
\bigcap_{i\in[n]}\Ker_{j}(my_i)^*$. Conversely, if $x\in
(I_K^{\bot})_{j+1}\setminus \Ker_{j}m^*$, there exists a monomial
$y$ in the basis $Y$ such that $<y,m^*(x)>\neq 0$ (as such monomials
$y$ span $R$ over $k$). Write $y=y_{i_0}y'$ for suitable $i_0$ (this
is possible as $j\geq \ddeg (m)$). Then $0\neq
<y,m^*(x)>=<y',y_{i_0}^*(m^*(x))>$, in particular $x\notin \Ker
(my_{i_0})^*$ $\square$
\\
\textbf{Convention}: From now on when writing the relation of
monomials $z<_Ly$ it will mean that we assume $\ddeg  (z)=\ddeg (y)$
even if we do not say so explicitly.

\begin{prop} \label{prop.2s}
Let $y,z$ be monomials in the basis $Y$. Then
$$
|\{z\in \GIN (K): \ddeg (z)=\ddeg  (y), y\leq_{L}z\}| =
\ddim \bigcap_{z<_{L}y}\Ker_{\ddeg (y)-1} z^*.
$$
In particular,
$$y\in \GIN (K)\ \iff\ \bigcap_{z<_{L}y}\Ker_{\ddeg  (y)-1}z^* \nsubseteq \Ker_{\ddeg  (y)-1}y^*.$$
\end{prop}
$Proof$: First note that the intersection on RHS does not change
when replacing the $z$'s with their orthogonal projection on
$I_K^{\bot}$ (as for $i\in [n]\setminus K_0$,
$x_i^*(I_K^{\bot})=0$), thus we may assume $[n]=K_0$.

Consider the $|\{z:z<_Ly,\ddeg  (z)=\ddeg  (y)\}|\times
\ddim(I_K^{\bot})_{\ddeg  (y)}$ matrix $A$ with $A_{z,x^a}=<z,x^a>$.
Then $RHS= \ddim_k(\Ker(A))= \ddim_k((I_K^{\bot})_{\ddeg
(y)})-\ddim_k(\im(A))=LHS$. $\square$

\subsubsection{Dividing $\GIN (K)$ into intervals}
For a monomial $y^a$ and $i\leq \ddeg  (y^a)$ let $\iinit_i (y^a)$
be the lexicographically least monomial of degree $i$ which divides
$y^a$. For every $i>0$ define the following subsets of monomials
with variables in $Y$:
$$J_{y}^{i} = J_{y}^{i}(n) = \{m\in \mathbb{Z}_{+}^{Y}: y|m, \ddeg  (m)=\ddeg  (y)+i, \iinit_{\ddeg  (y)}(m)=y\}.$$
Let
  $y_{(i)}^{(m)} = y_{(i)}^{(m)}(n)=\min_{<_{L}}J_{y}^{i}(n)$ and
 $y_{(i)}^{(M)} = y_{(i)}^{(M)}(n)=\max_{<_{L}}J_{y}^{i}(n)$.

\begin{prop} \label{prop.*s}
Let $K_{0}\subseteq [n]$, $S\subseteq [n]$, $i>0$. Then:
$$|J_{y}^{i}\cap \GIN (K)| =  \ddim \bigcap_{z<_{L}y}\Ker_{\ddeg  (y)+i-1} z^* -
\ddim \bigcap_{z\leq_{L}y}\Ker_{\ddeg  (y)+i-1} z^*.$$
\end{prop}
$Proof$:
   By Proposition \ref{prop.1s},
$$\ddim \bigcap_{z<_{L}y}\Ker_{\ddeg  (y)+i-1} z^*=
   \ddim \bigcap_{z<_{L}y}\bigcap_{j\in [n]}\Ker_{\ddeg  (y)+i-1} (zy_j)^* =...=$$
$$\ddim \bigcap_{z<_{L}y}\bigcap_{t\in (\mathbb{Z}_+^Y)_i}\Ker_{\ddeg  (y)+i-1} (zt)^* =
\ddim \bigcap_{z<_{L}y_{(i)}^{(m)}}\Ker_{\ddeg  (y_{(i)}^{(m)})-1} z^*
$$
which, by Proposition \ref{prop.2s}, equals
   $$|\{m\in\GIN (K): \ddeg  (m)=\ddeg  (y)+i,\ y_{(i)}^{(m)}\leq_{L}m\}|.$$
Similarly,
$$\ddim \bigcap_{z\leq_{L}y}\Ker_{\ddeg  (y)+i-1} z^*=
\ddim \bigcap_{z\leq_{L}y_{(i)}^{(M)}}\Ker_{\ddeg  (y_{(i)}^{(M)})-1} z^* =$$
$$
|\{m\in\GIN (K): \ddeg  (m)=\ddeg  (y)+i,\ y_{(i)}^{(M)}<_{L}m\}|.
$$
$\square$

\subsection{Basic properties}\label{subsec:basic properties}
Both exterior and symmetric shifting share the following basic
properties; see \cite{skira} for more details and references to the
original proofs. We denote both shifting operators by $K\mapsto
\Delta(K)$.
\begin{thm}(Kalai)\label{thm:basic properties}
Let $K$ and $L$ be simplicial complexes, and $\Delta$ denotes
algebraic shifting. Then:
\begin{enumerate}
\item $\Delta(K)$ is a simplicial complex.
\item $\Delta(K)=\Delta(L)$ for $L$ combinatorially isomorphic to $K$.
\item $f(K)=f(\Delta(K))$.
\item $\beta(K)=\beta(\Delta(K))$.
\item $\Delta(K)$ is shifted (assume $\rm{char}(K)=0$ for $\Delta^s$).
\item $\Delta^2=\Delta$ (assume $\rm{char}(K)=0$ for $\Delta^s$).
\item $\Delta(K)$ is the same for fields with the same characteristic.
\item If $L\subseteq K$ then $\Delta(L)\subseteq \Delta(K)$.
\end{enumerate}
\end{thm}
Later we will prove further properties of algebraic shifting, and
will exhibit properties which hold only for one of the two versions.


\chapter{Algebraic Shifting and Basic Constructions on Simplicial Complexes}\label{chapter:ABC}

\section{Shifting union of simplicial complexes}\label{sec2}
\subsection{General unions}
\begin{prob}\label{U problem}(\cite{skira}, Problem 13)
Given two simplicial complexes $K$ and $L$, find all possible
connections between $\Delta(K \cup L)$, $\Delta(K)$, $\Delta(L)$
and $\Delta(K \cap L)$.
\end{prob}

\subsubsection{Exterior shifting}
 We look on $\bigwedge(K\cup L)$, $\bigwedge(K\cap
L)$, $\bigwedge(K)$ and
  $\bigwedge(L)$ as subspaces of $\bigwedge(V)$ where
  $V=\kspan\{e_{1},...,e_{n}\}$ and $[n]=(K\cup L)_{0}$. As before,
  the $f_{i}$'s are generic linear combinations of the $e_{j}$'s
  w.r.t. $\mathbb{F}$. Let $S\subseteq [n]$, $|S|=s$ and $1\leq j$.
First we find a connection between boundary operations on the
spaces associated with $K$, $L$, $K\cap L$ and $K\cup L$ via the
following commutative diagram of exact sequences:
\begin{equation} \label{diagram}
\begin{CD}
0@> >>\bigwedge^{j+s}(K\cap L)@>i>>\bigwedge^{j+s}K \bigoplus \bigwedge^{j+s}L@>p>>\bigwedge^{j+s}(K\cup L)@> >>0\\
@VV V @VVfV @VVgV @VVhV @VV V \\
0@> >>\oplus\bigwedge^{j}(K\cap L)@>\oplus i>>\oplus\bigwedge^{j}K \bigoplus \oplus\bigwedge^{j}L@>\oplus p>>\oplus\bigwedge^{j}(K\cup L)@> >>0\\
\end{CD}
\end{equation}
where all sums $\oplus$ in the bottom sequence are taken over
$\{A: A<_{L}S, |A|=s\}$ and $i(m)=(m,-m)$, $p((a,b))=a+b$, $\oplus
i(m)=(m,-m)$, $\oplus p((a,b))=a+b$,
$f=\oplus_{A<_{L}S}f_{A}\lfloor(K\cap L)$,
$g=(\oplus_{A<_{L}S}f_{A}\lfloor(K),\oplus_{A<_{L}S}f_{A}\lfloor(L))$
and $h=\oplus_{A<_{L}S}f_{A}\lfloor(K\cup L)$.

By the snake lemma, (\ref{diagram}) gives rise to the following
exact sequence:
\begin{equation} \label{snake}
\begin{CD}
0@> >>ker (f)@> >>ker (g)@> >>ker (h)@>\delta>> \\
coker (f)@> >>coker (g)@> >>coker (h)@> >>0\\
\end{CD}
\end{equation}
where $\delta$ is the connecting homomorphism. Let
($\ref{diagram}'$) be the diagram obtained from (\ref{diagram}) by
replacing $A<_{L}S$ with $A\leq_{L}S$ everywhere, and renaming the
maps by adding a superscript to each of them. Let ($\ref{snake}'$)
be the sequence derived from ($\ref{diagram}'$) by applying to it
the snake lemma. If $\delta=0$ in (\ref{snake}), and also the
connecting homomorphism $\delta'=0$ in ($\ref{snake}'$), then by
Proposition \ref{prop.*} the following additive formula holds:
\begin{equation} \label{additive-formula}
|I_{S}^{j}\cap \Delta^e(K\cup L)|=|I_{S}^{j}\cap
\Delta^e(K)|+|I_{S}^{j}\cap \Delta^e(L)|-|I_{S}^{j}\cap \Delta^e(K\cap
L)|.
\end{equation}

\begin{thm} \label{gen}
Let $K$ and $L$ be two simplicial complexes, and let $d$ be the
dimension of $K\cap L$. For every
   subset $A$ of the set of vertices $[n]=(K\cup L)_{0}$, the
   following additive formula holds:
   \begin{eqnarray}\label{eq.genKuL}
|I_{A}^{d+2}\cap \Delta^e(K\cup L)|=|I_{A}^{d+2}\cap
\Delta^e(K)|+|I_{A}^{d+2}\cap \Delta^e(L)|.
   \end{eqnarray}
\end{thm}

$Proof$: Put $j=d+2$ in (\ref{diagram})
and in ($\ref{diagram}'$). Thus, the range and domain of $f$ in
(\ref{diagram}) and of $f'$ in ($\ref{diagram}'$) are zero, hence
$ker f=coker f=0$ and $ker f'=coker f'=0$, and the theorem follows. $\square$
\\

It would be interesting to understand what extra
information about $\Delta(K\cup L)$ we can derive by using more of
the structure of $\Delta(K\cap L)$, and not merely its dimension.
In particular, it would be interesting to find combinatorial
conditions that imply the vanishing  of $\delta$ in (\ref{snake}).
The proof of Theorem \ref{clique-sum} provides a step in this direction.
 The Mayer-Vietoris long exact sequence (e.g. \cite{Munkres} p.186) gives some information of this type,
  by the interpretation of the Betti vector using the shifted complex \cite{BK}, mentioned in the Introduction.

\subsubsection{Symmetric shifting}
Let $S$ be a monomial of degree $s$ in the generic basis
$Y=\{y_i\}_{i}$ w.r.t. the basis $X=\{x_i\}_{i}$ of $R=k[x_i: i\in
(K\cup L)_0]$, and let $j>0$ be an integer. Analogously to
(\ref{diagram}), the following  is a commutative diagram of exact
sequences of subspaces of $R$:
\begin{equation} \label{diagram s}
\begin{CD}
0@> >>(I_{K\cap L}^{\bot})_{j+s}@>i>>(I_{K}^{\bot})_{j+s} \bigoplus (I_{L}^{\bot})_{j+s}@>p>>(I_{K\cup L}^{\bot})_{j+s}@> >>0\\
@VV V @VVf^*V @VVg^*V @VVh^*V @VV V \\
0@> >>\oplus(I_{K\cap L}^{\bot})_{j}@>\oplus i>>\oplus(I_{K}^{\bot})_{j} \bigoplus \oplus(I_{L}^{\bot})_{j}@>\oplus p>>\oplus(I_{K\cup L}^{\bot})_{j}@> >>0\\
\end{CD}
\end{equation}
where all sums $\oplus$ in the bottom sequence are taken over
$\{m: m<_{L}S, \ddeg(m)=s\}$ and $i(m)=(m,-m)$, $p((a,b))=a+b$, $\oplus
i(m)=(m,-m)$, $\oplus p((a,b))=a+b$, and
$f^*,g^*,h^*$ are $\oplus_{m<_{L}S}m^*$ restricted to the relevant subspaces of $R$.

As in the exterior case, by the snake lemma we obtain the following
exact sequence:
\begin{equation} \label{snake s}
\begin{CD}
0@> >>ker (f^*)@> >>ker (g^*)@> >>ker (h^*)@>\delta>> \\
coker (f^*)@> >>coker (g^*)@> >>coker (h^*)@> >>0\\
\end{CD}
\end{equation}
where $\delta$ is the connecting homomorphism, and similarly we get
($\ref{snake s}'$) with $\delta'$ as the connecting homomorphism
when applying the snake lemma to the diagram obtained from
(\ref{diagram s}) by replacing $m<_{L}S$ with $m\leq_{L}S$
everywhere.

Under which conditions do we get $\delta=0=\delta'$, which provides
an additive formula? Such conditions are given in the next two
subsections.

We do not know whether the symmetric analogue of Theorem \ref{gen} holds or not.


\subsection{How to shift a disjoint union?}\label{disj}
As a corollary to Theorem \ref{gen} we get the following
combinatorial formula for shifting the disjoint union of
simplicial complexes:
\begin{cor} \label{prop.6}
  Let $(K\dot{\cup}L)_{0}=[n], [n]\supseteq S=\{ s_{1}<\dots<s_{j}<s_{j+1}\}$. Then: $$S\in\Delta^e(K\dot{\cup}L)
  \Leftrightarrow s_{j+1}-s_{j}\leq |I^{1}_{\iinit_{|S|-1}(S)}\cap\Delta^e(K)|+|I^{1}_{\iinit_{|S|-1}(S)}\cap\Delta^e(L)|.$$
\end{cor}
 $Proof$: Put $d=-1$ and $A=\iinit_{|S|-1}(S)$ in Theorem \ref{gen},
 and by Proposition \ref{prop.3} we are done.
 $\square$
\\
Similarly, in the symmetric case note that $I_{\emptyset}^{\bot}=k$,
hence for disjoint union and $j=1,s>0$ $(I_{K\cap
L}^{\bot})_{j+s}=0$ in ($\ref{diagram s}$), $ker f^*=coker f^*=0$ in
($\ref{snake s}$) and $ker (f^*)'=coker (f^*)'=0$ in ($\ref{snake
s}'$). By Proposition \ref{prop.*s}, for every monomial $y\neq 1$ in
the basis $Y$, $$|J_{y}^{1}\cap \GIN(K\dot{\cup}L)|=|J_{y}^{1}\cap
\GIN(K)|+|J_{y}^{1}\cap \GIN(L)|.$$ Translating this into terms of
symmetric shifting, we obtain
\begin{cor} \label{prop.6s}
  Let $(K\dot{\cup}L)_{0}=[n], [n]\supseteq S=\{ s_{1}<\dots<s_{j}<s_{j+1}\}$. Then: $$S\in\Delta^s(K\dot{\cup}L)
  \Leftrightarrow s_{j+1}-s_{j}\leq |I^{1}_{\iinit_{|S|-1}(S)}\cap\Delta^s(K)|+|I^{1}_{\iinit_{|S|-1}(S)}\cap\Delta^s(L)|.$$
$\square$
\end{cor}

  As a corollary, we get the following nice equation, proposed by
Kalai \cite{skira}, for both versions of algebraic shifting (in the symmetric case assume $\cchar (k)=0$):

\begin{cor} \label{prop7}
  $\Delta(K\dot{\cup}L)=\Delta(\Delta(K)\dot{\cup}\Delta(L)).$
\end{cor}
$Proof$: $S\in\Delta(K\cup L)$ iff (by Corollaries \ref{prop.6},
\ref{prop.6s}) $s_{j+1}-s_{j}\leq D_{K}(S)+D_{L}(S)$ iff (by
Proposition \ref{prop.5} and its symmetric analogue)
$s_{j+1}-s_{j}\leq D_{\Delta(K)}(S)+D_{\Delta(L)}(S)$ iff (by
Corollaries \ref{prop.6}, \ref{prop.6s}) $S\in
\Delta(\Delta(K)\dot{\cup}\Delta(L))$. $\square$
\\
\\
\textbf{Remarks}: (1) Above a high enough dimension (to be
specified) all faces of the shifting of a union are determined by
the shifting of its components: Let $st(K\cap L)=\{\sigma\in K\cup
L: \sigma\cap(K\cap L)_{0}\neq \emptyset \}$. Then $\Delta(K)$ and
$\Delta(L)$ determine all faces of $\Delta(K\cup L)$ of dimension
$>\ddim(\st(K\cap L))$, by applying Corollary \ref{prop7} to the
subcomplex of $K\cup L$ spanned by the vertices $(K\cup
L)_{0}-(K\cap L)_{0}$, and using the basic properties Theorem \ref{thm:basic properties}(3,8).

    (2) Let $X$ be a $(t+l)\times (t+l)$
generic block matrix, with an upper block of size $t\times t$ and
a lower block of size $l\times l$. Although we defined the
shifting operator $\Delta=\Delta_{A}$ with respect to a generic
matrix $A$, the definition makes sense for any nonsingular matrix
(but in that case the resulting complex may not be shifted). Let
$K_{0}=[t]$ and $L_{0}=[t+1,t+l]$. Corollary \ref{prop7} can be
formulated as
$$\Delta\circ\Delta_{X}(K\dot{\cup}L)=\Delta(K\dot{\cup}L)$$
because $\Delta_{X}(K\dot{\cup}L)=\Delta(K)\dot{\cup}\Delta(L)$
(on the right hand side of the equation the vertices of the two
shifted complexes are considered as two disjoint sets). However,
there are simplicial complexes $C$ on $t+l$ vertices , for which
$\Delta\circ\Delta_{X}(C)\neq\Delta(C)$. For example, let $t=l=3$
and take the graph $G$ of the octahedron
$\{\{1\},\{4\}\}*\{\{2\},\{5\}\}*\{\{3\},\{6\}\}$. Then
$\Delta\circ\Delta_{X}(G)\ni \{4,5\} \notin \Delta (G)$ over $k=\mathbb{R}$, for both versions of shifting.

    (3) By induction, we get from Corollary \ref{prop7} that:
$$\Delta(\dot{\cup}_{1\leq i\leq l}K^{i})=\Delta(\dot{\cup}_{1\leq i\leq l}\Delta(K^{i}))$$
for any positive integer $l$ and disjoint simplicial complexes
$K^{i}$.

\subsection{How to shift a union over a simplex?}\label{KsimplexL}
In the case where $K\cap L=<\sigma>$ is a simplex and all of its
subsets, we also get a formula for $\Delta(K \cup L)$ in terms of
$\Delta(K)$, $\Delta(L)$ and $\Delta(K\cap L)$. This case
corresponds to the topological operation called connected sum; see
Example \ref{ex:stacked}.

\begin{thm}\label{clique-sum}
Let $K$ and $L$ be two simplicial complexes, where $K\cap
L=<\sigma>$ is the simplicial complex consisting of the set
$\sigma$ and all of its subsets. For every $i>0$ and every
   subset $S$ of the set of vertices $[n]=(K\cup L)_{0}$, the
   following additive formula holds:
$$|I_{S}^{i}\cap \Delta(K\cup_{\sigma}L)|=|I_{S}^{i}\cap\Delta(K)|+|I_{S}^{i}\cap\Delta(L)|-|I_{S}^{i}\cap 2^{[|\sigma|]}|.$$
\end{thm}

This theorem gives an explicit combinatorial description of
$\Delta(K\cup L)$ in terms of $\Delta(K)$, $\Delta(L)$ and
$\ddim(\sigma)$. In particular, any gluing of $K$ and $L$ along a
$d$-simplex results in the same shifted complex $\Delta(K\cup L)$,
depending only on $\Delta(K)$, $\Delta(L)$ and $d$.

\begin{cor} \label{corKcupL}
 Let $K$ and $L$ be simplicial complexes where $K\cap
L=<\sigma>$ is a complete simplicial complex.
  Let $(K\cup L)_{0}=[n], [n]\supseteq T=\{ t_{1}<\dots<t_{j}<t_{j+1}\}$. Then: $$T\in\Delta(K\cup L)
  \Leftrightarrow t_{j+1}-t_{j}\leq D_{K}(T)+D_{L}(T)-D_{<\sigma>}(T).$$
\end{cor}
$Proof$: Put $i=1$ and $S=\iinit_{|T|-1}(T)$ in Theorem
\ref{clique-sum}, and by Proposition \ref{prop.3} (valid for the symmetric case as well) we are done.
 $\square$

\begin{ex}\label{ex:stacked}
Let $S(d,n)$ be a $(d-1)$-dimensional stacked sphere on $n$
vertices. Then $\Delta(S(d,n))=\rm{span}\{ \{1,3,4,..,d,n\}, \{2,3,..,d+1\}
\}$
where $\rm{span}$ means taking the closure under the product
partial order $<_P$ and under inclusion.
\end{ex}
$proof$: Let $\sigma$ be the $d$-simplex and $\partial\sigma$ its
boundary. Clearly $\Delta(\partial\sigma)=\partial\sigma$, hence the
case $n=d+1$ follows. Proceed by induction on $n$: use Corollary
\ref{corKcupL} to determine the shifting of the union
$S(d,n)\cup\partial\sigma$ over a common facet. To obtain
$\Delta(S(d,n+1))$ one needs to delete from it one facet. This must
be $\{2,3,...,d,d+2\}$, which represent the extra top homology.
$\square$

\textbf{Remark}:
A more complicated proof of Example \ref{ex:stacked} was given very recently by Murai \cite{MuraiCyclic}.
\\
\\
$Proof\ of\ Theorem\ \ref{clique-sum}$:  For a simplicial
complex $H$, let $\bar{H}$ denote the complete simplicial complex
$2^{H_{0}}$.

\textbf{Exterior case}: The inclusions $H\hookrightarrow \bar{H}$
for $H=K,L,<\sigma>$ induce a morphism from the commutative diagram
(\ref{diagram}) of $K$ and $L$ into the analogous commutative
diagram $(\overline{\ref{diagram}})$ of $\bar{K}$ and $\bar{L}$. By
functoriality of the sequence of the snake lemma, we obtain the
following commutative diagram:
\begin{equation} \label{snake*2}
\begin{CD}
0@> >>ker (f)@> >>ker (g)@> >>ker (h)@>\delta>>coker (f)@> >>... \\
@VV V @VVidV @VV V @VV V @VVidV \\
0@> >>ker (\bar{f})@> >>ker (\bar{g})@> >>ker (\bar{h})@>\bar{\delta}>>coker (\bar{f})@> >>... \\
\end{CD}
\end{equation}
where the bars indicate that ($\overline{\ref{diagram}}$) is obtained
from (\ref{diagram}) by putting bars over all the complexes and
renaming the maps by adding a bar over each map. Thus, if
$\bar{\delta}=0$ then also $\delta=0$, which, as we have seen,
implies (\ref{additive-formula}). The fact that
$\Delta(<\sigma>)=2^{[|\sigma|]}$ completes the proof.

We show now that $\bar{\delta}=0$. To simplify notation, assume
that $K$ and $L$ are complete complexes whose intersection is
$\sigma$ (which is a complete complex). Consider (\ref{diagram})
with $j=1$. (It is enough to prove Theorem \ref{clique-sum} for
$i=1$ as for every $i>1$, $S\subseteq [n]$ and $H$ a simplicial
complex on $[n]$, $I_{S}^{i}\cap H = \biguplus_{T\in
I_{S}^{i-1}}(I_{T}^{1}\cap H)$.) Let $m=m_{K}+m_{L}\in ker (h)$
where $\supp(m_{K})\subseteq K\setminus <\sigma>$. By commutativity
of the middle right square of (\ref{diagram}),
$\oplus_{A<_{L}S}f_{A}\lfloor(m_{K}),\oplus_{A<_{L}S}f_{A}\lfloor(m_{L})\in
\oplus_{A<_{L}S}\bigwedge^{1}(\sigma)$. If we show that
\begin{equation} \label{capIm}
\oplus_{A<_{L}S}f_{A}\lfloor(\bigwedge^{1+|S|}\sigma)=\oplus_{A<_{L}S}f_{A}\lfloor(\bigwedge^{1+|S|}K)
\bigcap \oplus_{A<_{L}S}\bigwedge^{1}(\sigma),
\end{equation}
then there exists $m'\in \bigwedge^{1+|S|}(\sigma)$ such that
$\oplus_{A<_{L}S}f_{A}\lfloor(m_{K}) =
\oplus_{A<_{L}S}f_{A}\lfloor(m') = f(m')$, hence
$\delta(m)=[f(m')]=0$ (where $[c]$ denotes the image of $c$ under
the projection onto $coker (f)$) i.e. $\delta=0$. (\ref{capIm})
follows from the intrinsic characterization of the image of the
maps it involves, given in Proposition \ref{intrinsicIm}. By
Proposition \ref{intrinsicIm}, the right hand side of
(\ref{capIm}) consists of all $x\in
\oplus_{R<_{L}S}\bigwedge^{1}K$ that satisfy $(a)$ and $(b)$ of
Proposition \ref{intrinsicIm} which are actually in
$\oplus_{A<_{L}S}\bigwedge^{1}(\sigma)$. By Proposition
\ref{intrinsicIm}, this is exactly the left hand side of
(\ref{capIm}). $\square$

\textbf{Symmetric case}: Repeating the arguments for the exterior
case, we need to show the following analogue of (\ref{capIm}) for
every monomial $S$ in the basis $Y$ of degree $s>0$:
\begin{equation} \label{capIm s}
\oplus_{m<_{L}S}m^*((I_{\sigma}^{\bot})_{s+1})=\oplus_{m<_{L}S}m^*((I_K^{\bot})_{s+1})
\bigcap \oplus_{m<_{L}S}(I_{\sigma}^{\bot})_{1}.
\end{equation}
This will follow from the intrinsic characterization of the image of the
maps it involves, given in Proposition \ref{intrinsicIm s}.
$\square$

The following generalizes a result of Kalai for graphs (\cite{56},
Lemma 3.7).
\begin{prop}\label{intrinsicIm}
Let $H$ be a complete simplicial complex with $H_{0}\subseteq
[n]$,  and let $S \subseteq [n]$, $|S|=s$. Then
$\oplus_{R<_{L}S}f_{R}\lfloor(\bigwedge^{1+s}H)$ is the set of all
$x=(x_{R}: R<_{L}S)\in \oplus_{R<_{L}S}\bigwedge^{1}H$ satisfying
the following:

(a) For all pairs $(i,R)$ such that $i\in R<_{L}S$:
 $<f_{i},x_{R}>=0$.

(b) For all pairs $(A,B)$ such that $A<_{L}S ,B<_{L}S$ and
$|A\vartriangle B|=2$: Denote $\{a\}=A\setminus B$ and
$\{b\}=B\setminus A$. Then $$-<f_{b},x_{A}>=(-1)^{\sgn_{A\cup
B}(a,b)}<f_{a},x_{B}>$$ where $\sgn_{A\cup B}(a,b)$ is the number
modulo $2$ of elements between $a$ and $b$ in the ordered set
$A\cup B$.
\end{prop}

$Proof$:
Let us verify first that every element in
$\im=\oplus_{R<_{L}S}f_{R}\lfloor(\bigwedge^{1+s}H)$ satisfies
$(a)$ and $(b)$. Let $y\in \bigwedge^{1+s}H$. If $i\in R<_{L}S$
then $<f_i,f_{R}\lfloor y>=<f_i\wedge f_{R},y>=<0,y>=0$, hence
$(a)$ holds. For $i\in T\subseteq [n]$ for some $n$, let
$\sgn(i,T)=|\{t\in T: t<i\}|(\mmod\ 2)$. If $A,B<_{L}S$,
$\{a\}=A\setminus B$ and $\{b\}=B\setminus A$ then
$-<f_b,f_{A}\lfloor y>=-<f_b\wedge f_{A},y>=-(-1)^{\sgn(b,A\cup
B)}<f_{A\cup B},y>=-(-1)^{\sgn(b,A\cup B)}(-1)^{\sgn(a,A\cup
B)}<f_a\wedge f_{B},y>=(-1)^{\sgn_{A\cup B}(a,b)}<f_a,f_{B}\lfloor
y>$, hence $(b)$ holds.

We showed that every element of $\im$ satisfies $(a)$ and $(b)$.
Denote by $X$ the space of all $x\in
\oplus_{R<_{L}S}\bigwedge^{1}H$ satisfying $(a)$ and $(b)$. It
remains to show that $\ddim(X)=\ddim(\im)$.

Following the proof of Proposition \ref{prop.*},
$\ddim(\im)=\ddim(\bigwedge^{1+s}H)-\ddim(\bigcap_{R<_{L}S}Ker_{s}f_{R}\lfloor(\bigwedge^{1+s}H)=
|\{T\in \Delta(H): |T|=s+1\}|-|\{T\in \Delta(H): |T|=s+1,
S_{(1)}^{(m)}\leq_{L}T\}|= |\{T\in \Delta(H): |T|=s+1,
\iinit_s(T)<_{L}S \}|$. Let $h=|H_0|$ and $\ssum(T)=|\{ t\in T:
T\setminus \{t\}<_{L}S \}|$. Note that $\Delta(H)=2^{[h]}$.
Counting according to the initial $s$-sets, we conclude that in
case $s<h$,
\begin{eqnarray} \label{dimIm}
\ddim(\im)=|\{R: R<_{L}S, R\subseteq [h]\}|(h-s)- \nonumber\\
\sum\{\ssum(T)-1 : T\subseteq [h], |T|=s+1, \iinit_{s}(T)<_{L}S,
\ssum(T)>1\}.
\end{eqnarray}
In case $s\geq h$, $\ddim(\im)=0$.

Now we calculate $\ddim(X)$. Let us observe that every $x\in X$ is
uniquely determined by its coordinates $x_R$ such that $R\subseteq
[h]$: Let $i\in R\setminus [h]$, $R<_{L}S$. Every $j\in
[h]\setminus R$ gives rise to an equation $(b)$ for the pair
$(R\cup j\setminus i,R)$ and every $j\in R\cap [h]$ gives rise to
an equation $(a)$ for the pair $(j,R)$. Recall that $x_R$ is a
linear combination of the form $x_R= \sum_{l\in
[h]}\gamma_{l,R}e_l$ with scalars $\gamma_{l,R}$. Thus, we have a
system of $h$ equations on the $h$ variables $(\gamma_{l,R})_{l\in
[h]}$ of $x_{R}$, with coefficients depending only on $x_{F}$'s
with $F<_{L}R$ (actually also $|F\cap [h]|=1+|R\cap [h]|$) and on
the generic $f_k$'s, $k\in [n]$. This system has a unique solution
as the $f_{k}$'s are generic. By repeating this argument we
conclude that $x_{R}$ is determined by the coordinates $x_F$ such
that $F\subseteq [h]$.

Let $x(h)$ be the restriction of $x\in X$ to its $\{x_R:
R\subseteq [h], R<_{L}S \}$ coordinates, and let $X(h)=\{x(h):
x\in X\}$. Then $\ddim(X(h))=\ddim(X)$.

Let $[a]$, $[b]$ be the matrices corresponding to the equation
systems $(a)$, $(b)$ with variables $(\gamma_{l,T})_{l\in
[h],T<_{L}S}$ restricted to the cases $T\subseteq [h]$ and
$A,B\subseteq [h]$, respectively. $[a]$ is an $s\cdot
|\{R\subseteq [h]: R<_{L}S\}|\times h\cdot |\{R\subseteq [h]:
R<_{L}S\}|$ matrix and $[b]$ is a $|\{A,B\subseteq [h]: A,B<_{L}S,
|A\vartriangle B|=2\}|\times h\cdot |\{R\subseteq [h]: R<_{L}S\}|$
matrix.

We observe that the row spaces of $[a]$ and $[b]$ have a zero
intersection. Indeed, for a fixed $R\subseteq [h]$, the row space
of the restriction of $[a]$ to the $h$ columns of $R$ is in
$\kspan\{f_{i}^{0}: i\in R\}$ (recall that $f_{i}^{0}$ is the
obvious projection of $f_i$ on the coordinates $\{e_j:j\in
H_0\}$), and the row space of the restriction of $[b]$ to the $h$
columns of $R$ is in $\kspan\{f_{j}^{0}: j\in [h]\setminus R\}$. But as
the $f_{k}^{0}$'s, $k\in [h]$, are generic, $\kspan\{f_{k}^{0}: k\in
[h]\}=\bigwedge^{1}H$. Hence $\kspan\{f_{i}^{0}: i\in R\}\cap
\kspan\{f_{j}^{0}: j\in [h]\setminus R\}=\{0\}$. We conclude that
the row spaces of $[a]$ and $[b]$ have a zero intersection.

$[a]$ is a diagonal block matrix whose blocks are generic of size
$s\times h$, hence
\begin{eqnarray} \label{dim[a]}
\rank([a])=s\cdot |\{R: R<_{L}S, R\subseteq [h]\}|
\end{eqnarray}
in case $s<h$.

Now we compute $\rank([b])$. For $T\subseteq [h]$, $|T|=s+1$, let
us consider the pairs in $(b)$ whose union is $T$. If $(A,B)$ and
$(C,D)$ are such pairs, and $A\neq C$, then $(A,C)$ is also such a
pair. In addition, if $A,B,C$ are different (the union of each two
of them is $T$) then the three rows in $[b]$ indexed by $(A,B)$,
$(A,C)$ and $(B,C)$ are dependent; the difference between the
first two equals the third. Thus, the row space of all pairs
$(A,B)$ with $A\cup B=T$ is spanned by the rows indexed
$(\iinit_{s}(T),B)$ where $\iinit_{s}(T)\cup B=T$.

We verify now that the rows $\bigcup\{(\iinit_{s}(T),B):
\iinit_{s}(T)\cup B=T\subseteq [h], |T|=s+1, |B|=s\}$ of $[b]$ are
independent. Suppose that we have a nontrivial linear dependence
among these rows. Let $B'$ be the lexicographically maximal
element in the set of all $B$'s appearing in the rows $(A,B)$ with
nonzero coefficient in that dependence. There are at most $h-s$
rows with nonzero coefficient whose restriction to their $h$
columns of $B'$ is nonzero (they correspond to $A$'s with
$A=\iinit_s(B\cup\{i\}$ where $i\in [h]\setminus B$). Again, as the
$f_{i}$'s are generic, this means that the restriction of the
linear dependence to the $h$ columns of $B'$ is nonzero, a
contradiction. Thus,
\begin{equation} \label{dim[b]}
\rank([b])=|\bigcup_{B<_{L}S}\{(\iinit_s(T),B): \iinit_{s}(T)\cup
B=T\subseteq [h], |T|=s+1, |B|=s\}|\nonumber
\end{equation}
\begin{equation}
= \sum\{\ssum(T)-1 : T\subseteq [h], |T|=s+1, \iinit_{s}(T)<_{L}S,
\ssum(T)>1\}.
\end{equation}
(Note that indeed $B<_{L}S$ implies $\iinit_{s}(T)<_{L}S$ as
$B\subseteq T$.)

For $s<h$, $\ddim(X(h))=h\cdot |\{R: R<_{L}S, R\subseteq
[h]\}|-\rank([a])-\rank([b])$, which by (\ref{dimIm}),
(\ref{dim[a]}) and (\ref{dim[b]}) equals $\ddim(\im)$. For $s\geq h$,
$\ddim(X(h))=0=\ddim(\im)$. This completes the proof. $\square$

For symmetric shifting, the following analogous assertion holds:
\begin{prop}\label{intrinsicIm s}
Let $H$ be a complete simplicial complex with $H_{0}\subseteq
[n]$,  and let $S$ be a monomial of degree $s$ in the generic basis $Y$ of $R_1=k[x_,...,x_n]_1$. Then
$\oplus_{m<_{L}S}m^*((I_H^{\bot})_{s+1})$ is the set of all
$x=(x_m: m<_{L}S)\in \oplus_{m<_{L}S}((I_H^{\bot})_{1})$ satisfying
the following:

(*) For all pairs of monomials in $Y$ $(A,B)$ such that $A<_{L}S
,B<_{L}S$ and $g=\gcd(A,B)$ has degree $s-1$, denote
$y_A=\frac{A}{g}$ and $y_B=\frac{B}{g}$. Then $<y_B,x_A>=<y_A,x_B>$.
\end{prop}

$Proof$: For every monomial $x^I\in I_H^{\bot}$ and $A,B$ as in (*),
indeed $<y_B,A^*(x^I)>=<y_BA,x^I>=<y_Ay_Bg,x^I>=<y_A,B^*(x^I)>$
hence all elements in $\oplus_{m<_{L}S}m^*((I_H^{\bot})_{s+1})$
satisfy (*).

For the converse implication, let us compute the dimensions of both
$k$-vector spaces. Denote for a monomial $T$ of degree $s+1$ in the
basis $Y$, $\ssum(T):=|\{i: y_i\mid T, T/y_i <S\}|$. Then:
\begin{equation}\label{dimLHSsym}
\ddim_k \oplus_{m<_{L}S}m^*((I_H^{\bot})_{s+1})= |\{R\in GIN(H): \rm{deg}(R)=s+1, \rm{\iinit}_s(R)<S\}|=\nonumber
\end{equation}
\begin{equation}
|H_0||\{R: R<S, \rm{supp}(R)\subseteq H_0\}| - \sum\{\rm{sum}T-1:\
\rm{supp}(T)\subseteq H_0, \iinit_s(T)<S, \rm{sum}T>1\}.
\end{equation}
Let $X(*)$ denote the space of elements $x=(x_m: m<_{L}S)\in
\oplus_{m<_{L}S}(I_H^{\bot})_{1}$ satisfying (*). Note that if $i\in
\rm{supp}(R)\nsubseteq H_0$ then the coordinate $x_R$ is determined
by the coordinates $\{x_{Ry_j/y_i}:j\in H_0\}$ via the equations (*)
for the pairs $(R,Ry_j/y_i)$. Iterating this argument shows that the
elements in $X(*)$ are determined by their coordinates which are
supported on $H_0$. Some of the equations in (*) are dependent, let
us find them a basis: consider the equations indexed by pairs
$(\iinit_s(T),B)$ where $\rm{supp}(T)\subseteq H_0$,
$\rm{deg}(T)=s+1$, and $\rm{lcm}(\iinit_s(T),B)=T$, as the rows of a
matrix. We now show that they are independent: assume by
contradiction the existence of a dependency, with all coefficients
being nonzero, and let $B'$ be the lexicographically maximal $B$ in
the pairs $(\iinit_s(T),B)$ which index it. Restrict the dependency
to the $H_0$ coordinates, and to the $|H_0|$ columns indexed by
$B'$. In the matrix that these columns form there are at most
$|H_0|$ nonzero rows, corresponding to pairs with $(\iinit_s(T),B')$
where $T=y_iB'$, $i\in H_0$, and they create a generic block, hence
the coefficients of these rows equal zero, a contradiction.
We conclude
\begin{equation}\label{dimRHSsym}
\ddim_k X(*)\leq
|H_0||\{R: R<S, \rm{supp}(R)\subseteq H_0\}| - \nonumber
\end{equation}
\begin{equation}
\sum\{\ssum(T)-1:\ \rm{supp}(T)\subseteq H_0, \iinit_s(T)<S, \ssum(T)>1\}.
\end{equation}
As $\oplus_{m<_{L}S}m^*((I_H^{\bot})_{s+1})\subseteq X(*)$, combining with Equation (\ref{dimLHSsym}) we conclude
$\oplus_{m<_{L}S}m^*((I_H^{\bot})_{s+1})=X(*)$. $\square$

\section{Shifting near cones}\label{sec4}
  A simplicial complex $K$ is called
a \emph{near cone} with respect to a vertex $v$ if for every $j\in S\in
K$ also $v\cup S\setminus j\in K$. We prove a
decomposition theorem for the shifted complex of a near cone, from
which the formula for shifting a cone Corollary \ref{prop.11cone} follows. As a preparatory step we introduce the
Sarkaria map, modified for homology.

\subsection{The Sarkaria map}\label{sec33}
    Let $K$ be a near cone with respect to a vertex $v=1$. Let
  $e=\sum_{i \in K_{0}}e_{i}$ and let $f=\sum_{i \in K_{0}}\alpha_{i}e_{i}$ be a linear combination
  of the $e_{i}$'s such that $\alpha_{i}\neq 0$ for every $i \in K_{0}$.
  Imitating the Sarkaria maps for cohomology \cite{SarkariaPanjab}, we get for homology
  the following linear maps:
  $$(\bigwedge K, e_{v}\lfloor) \overset{U}{\longrightarrow}(\bigwedge K, e\lfloor)
   \overset{D}{\longrightarrow}(\bigwedge K, f\lfloor)$$
 defined as follows: for $S\in K$
 $$U(e_{S})= \left\{ \begin {array}{ll}
   e_{S}-\sum_{i \in S}(-1)^{\sgn(i,S)} e_{v\cup S\setminus i} & \textrm{if
   $v\notin S$}\\
   e_{S} & \textrm{if $v\in S$} \end {array} \right. $$

$$D^{-1}(e_{S})=(\prod_{i\in S}\alpha_{i})e_{S}.$$
It is justified to write $D^{-1}$ as all the $\alpha_{i}$'s are
non zero.
\begin{prop} \label{prop.9Sarc}
The maps $U$ and $D$ are isomorphisms of chain complexes. In
addition they satisfy the following 'grading preserving' property:
if $S\cup T \in K$, $S\cap T=\emptyset$, then $$U(e_{S}\wedge
e_{T})=U(e_{S})\wedge U(e_{T})\  and\ D(e_{S}\wedge
e_{T})=D(e_{S})\wedge D(e_{T}).$$

\end{prop}
$Proof$: The check is straight forward.  First we check that $U$ and
$D$ are chain maps. Denote $\alpha_{S}=\prod_{i\in S}\alpha_{i}$.
For every $e_{S}$ where $S \in K$, $D$ satisfies
$$D\circ e\lfloor (e_{S})= D(\sum_{j\in S}(-1)^{\sgn(j,S)}e_{S\setminus j})
= \sum_{j\in
S}(-1)^{\sgn(j,S)}\frac{\alpha_{j}}{\alpha_{S}}e_{S\setminus j}$$
and
$$f\lfloor\circ D(e_{S})=f\lfloor(\frac{1}{\alpha_{S}}e_{S}) = \sum_{j\in
S}(-1)^{\sgn(j,S)}\frac{\alpha_{j}}{\alpha_{S}}e_{S\setminus j}.$$
For $U$: if $v\in S$ we have
$$U\circ e_{v}\lfloor (e_{S})=U(e_{S\setminus v}) =e_{S\setminus
v}-\sum_{i\in S\setminus v}(-1)^{\sgn(i,S\setminus v)}e_{S\setminus
i} = \sum_{j\in S}(-1)^{\sgn(j,S)}e_{S\setminus j}.$$ The last
equation holds because $v=1$. Further,
$$e\lfloor\circ U(e_{S})= e\lfloor(e_{S})=\sum_{j\in S}(-1)^{\sgn(j,S)}e_{S\setminus
j}.$$ If $v\notin S$ we have
$$U\circ e_{v}\lfloor (e_{S})=U(0)=0$$
and
$$e\lfloor\circ U(e_{S})= e\lfloor(e_{S})-e\lfloor(\sum_{j\in S}(-1)^{\sgn(j,S)}e_{S\cup v\setminus
j})=$$ $$\sum_{j\in S}(-1)^{\sgn(j,S)}e_{S\setminus j}- \sum_{i\in
S}(-1)^{\sgn(i,S)}\sum_{t\in S\cup v\setminus i}(-1)^{\sgn(t,S\cup
v\setminus i)}e_{S\cup v\setminus\{i,t\}}=$$
$$\sum_{j\in S}(-1)^{\sgn(j,S)}e_{S\setminus j}(1-(-1)^{\sgn(v,S\cup v\setminus
j)}) - \sum_{j,i\in S,i\neq j}(-1)^{\sgn(i,S)}(-1)^{\sgn(j,S\cup
v\setminus i)}e_{S\cup v\setminus \{i,j\}}.$$ In the last line,
the left sum is zero as $v=1$, and for the same reason the right
sum can be written as:
$$=\sum_{j,i\in S, i<j}((-1)^{\sgn(i,S)+\sgn(j,S\setminus i)}+(-1)^{\sgn(j,S)+\sgn(i,S\setminus
j)})e_{S\cup v\setminus \{i,j\}}.$$ As $i<j$, the $\{i,j\}$
coefficient equals
$$(-1)^{\sgn(i,S)+\sgn(j,S)+1}+(-1)^{\sgn(j,S)+\sgn(i,S)}=0,$$
hence $e\lfloor\circ U(e_{S})= U\circ e_{v}\lfloor (e_{S})$ for every
$S \in K$. By linearity of $U$ and $D$ (and of the boundary maps),
  we have that $U,D$ are chain maps. To show that $U,D$ are onto,
  it is enough to show that each $e_{S}$, where $S \in K$, is in
  their image. This is obvious for $D$. For $U$: if $v\in S$ then
  $U(e_{S})=e_{S}$, otherwise $e_{S}=U(e_{S})+\sum_{i \in S}(-1)^{\sgn(i,S)} e_{v\cup S\setminus i}$
  , which is a linear combination of elements in $\im(U)$, so $e_{S}\in
  \im(U)$ as well. Comparing dimensions, $U$ and $D$ are
  also 1-1.

We now show that $U$ 'preserves grading' in the described above
sense (for $D$ it is clear). For disjoint subsets of $[n]$ define
$\sgn(S,T)=|\{(s,t)\in S\times T: t<s\}|(\mmod\ 2)$. Let $S,T$ be
disjoint sets such that $S\cup T\in K$.  By $S\cup T$ we mean the
ordered union of $S$ and $T$ (and similarly for other set unions).

case 1: $v\notin S\cup T$.
 $$U(e_{S})\wedge U(e_{T})=$$ $$
e_{S}\wedge e_{T}+\sum_{i\in S}(-1)^{\sgn(i,S)}e_{S\cup v\setminus
i}\wedge e_{T}+\sum_{j\in T}(-1)^{\sgn(j,T)}e_{S}\wedge e_{T\cup
v\setminus j}=$$
$$(-1)^{\sgn(S,T)}(e_{S\cup T}
+\sum_{l\in S\cup T}(-1)^{\sgn(l,S\cup T)}e_{S\cup T\cup v\setminus
l})=$$
$$U(e_{S}\wedge e_{T}),$$
where the middle equation uses the fact that $v=1$, which leads to
the following sign calculation:
$$(-1)^{\sgn(i,S)}(-1)^{\sgn(S\cup v\setminus i,T)}= (-1)^{\sgn(i,S)+\sgn(S\setminus
i,T)}=$$
$$(-1)^{\sgn(i,S)+\sgn(S,T)+\sgn(i,T)}=(-1)^{\sgn(S,T)}(-1)^{\sgn(i,S\cup
T)}.$$

case 2: $v\in S\setminus T$.$$U(e_{S})\wedge U(e_{T})=e_{S}\wedge
(e_{T}-\sum_{t \in T}(-1)^{\sgn(t,T)} e_{T\cup v\setminus t}
)=e_{S}\wedge e_{T}=U(e_{S}\wedge e_{T}).$$

case 3: $v\in T\setminus S$. A similar calculation to the one for
case 2 holds.
 $\square$
\\
\textbf{Remark}: The 'grading preserving' property of $U$ and $D$
extends to the case where $S\cap T\neq \emptyset$ ($S,T\in K$),
but we won't use it here. One has to check that in this case
(where clearly $e_{S}\wedge e_{T}=0$):
$$U(e_{S})\wedge U(e_{T})=D(e_{S})\wedge D(e_{T})=0.$$

\subsection{Shifting a near cone: exterior case.}
\begin{thm} \label{thm.10'near-cone}
 Let $K$ be a near cone on a vertex set $[n]$ with respect to a vertex $v=1$.
Let $X=\{f_{i}: 1\leq i\leq n\}$ be some basis of $\bigwedge^{1} K$
such that $f_{1}$ has no zero coefficients as a linear combination
of some given basis elements $e_{i}$'s of $\bigwedge^{1} K$, and
such that for $g_{i-1}=f_{i}-<f_{i},e_{1}>e_{1}$, $Y=\{g_{i}: 1\leq
i\leq n-1\}$ is a linearly independent set.
Then:
$$\Delta_{X}(K)=(1*(\Delta_{Y}(\lk(v,K))+1))\cup B$$
where $B$ is the set $\{S\in \Delta_{X}(K): 1 \notin S\}$, $j*L:=\{j\cup T: T\in L\}$, $L+j:=\{T+j: T\in L\}$ and $T+j:=\{t+j: t\in T\}$.
\end{thm}
$Proof$: Clearly for every $l\geq 0$:
 $Ker_{l}e_{v}\lfloor =
\bigwedge^{l+1}\antist(v,K)$ and $\im_{l}e_{v}\lfloor =
\bigwedge^{l}\lk(v,K)$. Using the Sarkaria map $D\circ U$, we get that $\im f_{1}\lfloor$ is
isomorphic to $\bigwedge \lk(v,K)$ and is contained, because of
'grading preserving', in a sub-exterior-algebra generated by the
elements
$b_{i}=DU(e_{i})=\frac{1}{\alpha_{i}}e_{i}-\frac{1}{\alpha_{v}}e_{v}$,
$i\in K_{0}\setminus v$ (see
Proposition \ref{prop.9Sarc}).
 Let $S\subseteq [n], |S|=l, 1\notin
S$. Recall that $(g\wedge f)\lfloor h = g\lfloor(f\lfloor h)$. Now
we are prepared to shift.
$$\bigcap_{R<_{L}1\cup S}\Ker_{l} f_{R}\lfloor \cong \Ker_{l} f_{1}\lfloor \oplus
\bigcap_{1\notin R<_{L}S}\Ker f_{R}\lfloor (\im_{l} f_{1}\lfloor
\rightarrow k),$$ which by the Sarkaria map is isomorphic
to
\begin{equation}\label{eq3.th10}
\wedge^{l+1}\antist(v,K)\oplus \bigcap_{1\notin T<_{L}S}\Ker
(f_{T}\lfloor (DU(\wedge^{l}\lk(v,K))\rightarrow k).
\end{equation}
Denote by $\pi_{t}$ the natural projection $\pi_{t}: \kspan\{e_{R}:
|R|=t\}\rightarrow \kspan\{e_{R}: |R|=t, v\notin R\}$, and by $M$ the matrix
$(<\pi_{l}\circ(DU)^{*}f_{T} , e_{R}>)$ where $1\notin T
<_{L}S, R\in \lk(v,K)_{l-1}$.
Then
$$\bigcap_{1\notin T<_{L}S}\Ker (f_{T})\lfloor (DU(\wedge^{l}\lk(v,K))\rightarrow
k) \cong \Ker(M)$$
Let $G$ be the matrix $(<g_{T-1} , e_{R}>)$,
where $1\notin T <_{L}S, R\in \lk(v,K)_{l-1}$. Then $M$ is obtained
from $G$ by performing the following operations: multiplying rows by
nonzero scalars and adding to a row multiples of lexicographically
smaller rows. Thus, restricting to the first $m$ rows of each of
these two matrices we get matrices of equal rank, for every $m$. In
particular, $\Ker(M) \cong \Ker(G)$. Hence, using the proof of
Proposition \ref{prop.2} (note that the proof of Proposition
\ref{prop.2} can be applied to non-generic shifting as well), by
putting $Q=T-1$ in $G$ we get:
$$\ddim \bigcap_{R<_{L}1\cup S}Ker_{l} f_{R}\lfloor(K) =
\ddim \wedge^{l+1}\antist(v,K)+ \ddim \bigcap_{Q<_{L}S-1}Ker_{l-1}
g_{Q}\lfloor (\lk(v,K)).$$
 As the left summand in the right hand side is a constant
independent of $S$, it is canceled when applying the last part of
Proposition \ref{prop.2}, and we get:
$$1\dot{\!\cup}S \in \Delta_{X}(K)\  \Leftrightarrow $$
$$\ddim \bigcap_{R<_{L}1\cup S}\Ker_{l} f_{R}\lfloor(K)
> \ddim \bigcap_{R\leq_{L}1\cup S}\Ker_{l} f_{R}\lfloor(K)\
\Leftrightarrow$$ $$\ddim \bigcap_{T<_{L}S-1}\Ker_{l-1}\ g_{T}\lfloor
(\lk(v,K)) > \ddim \bigcap_{T\leq_{L}S-1}\Ker_{l-1}\
g_{T}\lfloor(\lk(v,K))\ \Leftrightarrow$$
$$S-1\in\Delta_{Y}(\lk(v,K)).$$
Thus we get the
claimed decomposition of $\Delta_{X}(K)$. $\square$
\\

As a corollary we get the following decomposition theorem for the
generic shifted complex of a near cone.
\begin{thm} \label{thm.10near-cone}
 Let $K$ be a near cone with respect to a vertex $v$. Then
$$\Delta^e(K)=(1*\Delta^e(\lk(v,K)))\cup B,$$
where $B$ is the set $\{S\in \Delta^e(K): 1 \notin S\}$.
\end{thm}
$Proof$: Apply Theorem \ref{thm.10'near-cone} for the case where $X$ is generic. In this
case, $Y$ is also generic, and the theorem follows .$\square$

As a corollary we get the following property \cite{skira}:
\begin{cor} \label{prop.11cone}
$\Delta^e\circ \Cone = \Cone \circ\Delta^e.$
\end{cor}
$Proof$: Consider a cone over $v$: $\{v\}*K$. By Theorem
\ref{thm.10near-cone}, $\{1\}*\Delta^e(K)\subseteq \Delta^e(\{v\}*K)$,
but those two simplicial complexes have equal $f$-vectors, and
hence, $\{1\}*\Delta^e(K)=\Delta^e(\{v\}*K)$.$\square$\\
\textbf{Remarks}: (1) Note that by associativity of the join
operation, we get by Corollary \ref{prop.11cone}:
$\Delta^e(K[m]*K)=K[m]*\Delta^e(K)$ for every $m$, where $K[m]$ is the
complete simplicial complex on $m$ vertices.
\\
(2) Using the notation in Theorem \ref{thm.10'near-cone} we get:
$\Delta_{X}\circ \Cone = \Cone \circ\Delta_{Y}.$
\\
(3) Recently, it was shown in \cite{Babson-Novik-Thomas-Cone} that
$\Delta^s\circ \Cone = \Cone \circ\Delta^s$ where the field is of
characteristic zero, as was claimed by Kalai in \cite{skira}.

\begin{de}\label{i near cone}
 $K$ is an $i-near \ cone$ if there exists a sequence of simplicial
 complexes
$K=K(0)\supset K(1)\supset\dots \supset K(i)$ such that for every
$1\leq j\leq i$ there is a vertex $v_{j}\in K(j-1)$ such that
$K(j)=\antist(v_{j},K(j-1))$ and $K(j-1)$ is a near cone w.r.t.
$v_{j}$.
\end{de}
\textbf{Remark}: An equivalent formulation is that there exists a
permutation $\pi: K_{0}=[n]\rightarrow [n]$ such that

$$\pi(i)\in S\in K, 1\leq
l<i \Rightarrow (S\cup \pi(l)\setminus \pi(i))\in K,$$ which is
more compact but less convenient for the proof of the
 following generalization of Theorem
\ref{thm.10near-cone}:
\begin{cor} \label{prop.12near-near-cone}
Let $K$ be an i-near cone. Then
$$\Delta^e(K)=B\cup\biguplus_{1\leq j\leq i}j*(\Delta^e(\lk(v_{j},K(j-1)))+j),$$
where $B=\{S\in \Delta^e(K): S\cap[i]=\emptyset\}$.
\end{cor}
$Proof$:
The case $i=1$ is Theorem \ref{thm.10near-cone}.
By induction hypothesis,
$\Delta(K)=\tilde{B}\cup\biguplus_{1\leq j\leq
i-1}j*(\Delta^e(K(j-1))+j)$ where $\tilde{B}=\{S\in \Delta^e(K):
S\cap[i-1]=\emptyset\}$. We have to show that
\begin{equation}\label{eq.th12}
\{S\in\Delta^e(K): \mmin\{j\in
S\}=i\}=i*(\Delta^e(\lk(v_{i},K(i-1)))+i)).
\end{equation}
For $|S|=l$ with $\mmin\{j\in S\}=i$, we have
\begin{equation}\label{eq.th12+}
\bigcap_{R<_{L}S}\Ker_{l-1} f_{R}\lfloor(K)\ =
(\bigcap_{j<i}\ \bigcap_{R:|R|=l, j\in R} \Ker_{l-1}
f_{R}\lfloor)\cap (\bigcap_{R<_{L}S: \mmin(R)=i}\Ker_{l-1}
f_{R}\lfloor).\
\end{equation}
By repeated application of Proposition \ref{prop.1}, for each
$j<i$,  $$\bigcap_{R:|R|=l, j\in R} \Ker_{l-1}
f_{R}\lfloor=\Ker_{l-1}f_{j}.$$ Hence, (\ref{eq.th12+}) equals
$$(\bigcap_{j<i}\Ker_{l-1} f_{j}\lfloor)\ \cap\ (\bigcap_{R<_{L}S: \mmin(R)=i}\Ker_{l-1}
f_{R}\lfloor)\ =$$
$$\bigcap_{R<_{L}S: \mmin(R)=i}\Ker_{l-1} f_{R}\lfloor(A_l),$$ where
$A_l=\bigcap_{j<i}\Ker f_{j}\lfloor(\bigwedge^l K)$. Let $A=\oplus_lA_l$.

By repeated application of the Sarkaria map, we get that
$A\cong\bigwedge K(i-1)$ as graded chain complexes. Now we will show
that
\begin{equation}\label{eq.th12++}
\ddim \bigcap_{R<_{L}S: \mmin(R)=i}\Ker_{l-1} f_{R}\lfloor(A)= \ddim
\bigcap_{R<_{L}S-(i-1)}\Ker_{l-1} f_{R}\lfloor(\bigwedge K(i-1)).
\end{equation}
Let $\varphi: \bigwedge K(i-1)\rightarrow A$ be the Sarkaria
isomorphism, and let $f$ be generic w.r.t. the basis
$\{e_{i},..,e_{n}\}$ of $\bigwedge^{1}K(i-1)$. Then $\varphi(f)$ is
generic w.r.t. the basis $\{\varphi(e_{i}),..,\varphi(e_{n})\}$ of
$A$. We can choose a generic $\bar{f}$ w.r.t. $\{e_{1},..,e_{n}\}$
such that $<\bar{f},\varphi(e_{j})>=<\varphi(f),\varphi(e_{j})>$ for
every $i\leq j\leq n$. Actually, we can do so for $n-i$ generic
$f_{j}$'s simultaneously (as multiplying a nonsingular matrix over a
field by a generic matrix over the same field results in a generic
matrix over that field). We get that
$$\bigcap_{R<_{L}S: min(R)=i}Ker_{l-1} \bar{f}_{R}\lfloor(A)= \bigcap_{R<_{L}S-(i-1)}Ker_{l-1} \varphi(f_{R})\lfloor(A)$$
$$\cong \bigcap_{R<_{L}S-(i-1)}Ker_{l-1} f_{R}\lfloor(\bigwedge K(i-1)).$$
As both the $f_{i}$'s and the $\bar{f}_{i}$'s are generic,
$\bigcap_{R<_{L}S: \mmin(R)=i}\Ker_{l-1} f_{R}\lfloor(A) \cong
\bigcap_{R<_{L}S: \mmin(R)=i}\Ker_{l-1} \bar{f}_{R}\lfloor(A)$ and
(\ref{eq.th12++}) follows. By applying Theorem
\ref{thm.10near-cone} to the near cone $K(i-1)$, we see that
(\ref{eq.th12}) is true, which completes the proof. $\square$
\\

From our last corollary we obtain a new proof of a well known
property of algebraic shifting, proved by Kalai \cite{Kalai-SymmMatroids}:
\begin{cor} \label{prop.13 delta^2}
$\Delta^e\circ\Delta^e=\Delta^e.$
\end{cor}
$Proof$: For every simplicial complex $K$ with $n$ vertices, $\Delta^e(K)$ is shifted,
(hence an $n$-near cone), and so are all the $\lk(i,(\Delta^e K)(i-1))$'s associated to
it. By induction on the number of vertices, $\Delta^e(\lk(i,(\Delta^e K)(i-1)))=\lk(i,(\Delta^e
K)(i-1))-i$ for all $1\leq i\leq n$. Thus, applying Corollary
\ref{prop.12near-near-cone} to the $n$-near cone $\Delta^e(K)$, we get
$\Delta^e(\Delta^e(K))= \Delta^e(K)$.$\square$


\section{Shifting join of simplicial complexes}\label{sec k*l}
Let $K,L$ be two disjoint simplicial complexes (they include the
empty set). Recall that their join is the simplicial complex $K*L=\{S\cup T:
S\in K, T\in L\}$. Using
the K\"{u}nneth theorem with field coefficients (see \cite{Munkres},
Theorem 58.8 and ex.3 on p.373) we can describe its homology in
terms of the homologies of $K$ and $L$:
$$H_{i}(K\times L)\cong \bigoplus_{k+l=i}H_{k}(K)\otimes
H_{l}(L)$$
and the exact sequence $$0\rightarrow \tilde{H}_{p+1}(K*L)\rightarrow
\tilde{H}_{p}(K\times L)\rightarrow
 \tilde{H}_{p}(K)\oplus \tilde{H}_{p}(L)\rightarrow 0.$$
Recalling that $\beta_{i}(K)=|\{S\in \Delta(K): |S|=i+1, S\cup
1\notin \Delta(K)\}|$ (\cite{BK} for exterior case, \cite{Herzog} for symmetric case), we get a
description of the number of faces in $\Delta(K*L)_{i}$ which
after union with $\{1\}$ are not in $\Delta(K*L)$, in terms of
numbers of faces of that type in $\Delta(K)$ and $\Delta(L)$. In
particular, if the dimensions of $K$ and $L$ are strictly greater
than $0$, the K\"{u}nneth theorem implies:
$$\beta_{\ddim(K*L)}(K*L)=\beta_{\ddim(K)}(K)\beta_{\ddim(L)}(L)$$
and hence
$$|\{S\in \Delta(K*L): 1\notin S, |S|=\ddim(K*L)+1\}|=$$
$$|\{S\in \Delta(K): 1\notin S, |S|=\ddim(K)+1\}|\times|\{S\in \Delta(L):
1\notin S, |S|=\ddim(L)+1\}|.$$ We now show that more can be said
about the faces of maximal size in $K*L$ that represents homology
of $K*L$:
\begin{thm}\label{K*L,max}
Let $|(K*L)_{0}|=n$. For every $i\in [n]$ $$|\{S\in \Delta^e(K*L):
[i]\cap S=\emptyset, |S|=\ddim(K*L)+1\}|=$$ $$|\{S\in \Delta^e(K):
[i]\cap S=\emptyset, |S|=\ddim(K)+1\}|\times |\{S\in \Delta^e(L):
[i]\cap S=\emptyset, |S|=\ddim(L)+1\}|.$$
\end{thm}
$Proof$: For a generic
$f=\sum_{v\in K_{0}\cup L_{0}}\alpha_{v}e_{v}$ decompose
$f=f(K)+f(L)$ with supports in $K_{0}$ and $L_{0}$ respectively.
Denote $\ddim(K)=k, \ddim(L)=l$, so $\ddim(K*L)=k+l+1$. Observe that
$(f(K)\lfloor(K*L)_{k+l+1})\cap (f(L)\lfloor(K*L)_{k+l+1})=\{0\}$.
Denote by $f\lfloor(K)$ the corresponding generic boundary
operation on $\kspan\{e_{S}:S\in K\}$, and similarly for $L$.
Looking at $\bigwedge(K*L)$ as a tensor product $(\bigwedge
K)\otimes (\bigwedge L)$ we see that $\Ker_{k+l+1}f(K)\lfloor$
equals $\Ker_{k}f(K)\lfloor\arrowvert_{\bigwedge K}\otimes
\bigwedge^{1+l} L$, and also $\Ker_{k}f(K)\lfloor|_{\bigwedge
K}\cong \Ker_{k}f\lfloor(K)$, and similarly when changing the roles
of $K$ and $L$. Hence, we get
$$\Ker_{k+l+1}f\lfloor =
\Ker_{k+l+1}f(K)\lfloor\ \cap
 \Ker_{k+l+1}f(L)\lfloor \cong
 \Ker_{k}f\lfloor(K)\otimes\Ker_{l}f\lfloor(L).$$
For the first $i$ generic $f_{j}$'s, by the
same argument, we have:
$$\bigcap_{j\in [i]}\Ker_{k+l+1}f_{j}\lfloor =
\bigcap_{j\in [i]}\Ker_{k+l+1}f_{j}(K)\lfloor\ \cap
\bigcap_{j\in [i]}\Ker_{k+l+1}f_{j}(L)\lfloor \cong$$
$$\bigcap_{j\in [i]}\Ker_{k}f_{j}\lfloor(K)\otimes \bigcap_{j\in [i]}\Ker_{l}f_{j}\lfloor(L).$$
By Propositions \ref{prop.1} and \ref{prop.2} we get the claimed assertion.
$\square$

For symmetric shifting, the analogous assertion
to Theorem \ref{K*L,max} is false:
\begin{ex}\label{ex:JoinSymmFails}
Let each of $K$
and $L$ consist of three points. Thus, $K*L=K_{3,3}$ is the
complete bipartite graph with $3$ vertices on each side. By
Theorem \ref{K*L,max}, $\{3,4\}\in \Delta^e(K_{3,3})$, but
$\{3,4\}\notin \Delta^{s}(K_{3,3})$.
\end{ex}
We now deal with the conjecture (\cite{skira}, Problem $12$)
\begin{equation}\label{eq**}
\Delta(K*L)=\Delta(\Delta(K)*\Delta(L)).
\end{equation}
We give a counterexample showing that it is false even if we
assume that one of the complexes $K$ or $L$ is shifted. Denote by
$\Sigma K$ the suspension of $K$, i.e. the join of $K$ with the
(shifted) simplicial complex consisting of two points.
\begin{ex}\label{JoinConjFails}
Let $B$ be the graph consisting of two disjoint
edges. In this case $\Delta (\Sigma(B)) \setminus \Delta (\Sigma
(\Delta((B))) = \{\{1,2,6\}\}$ and $\Delta (\Sigma (\Delta(B)))
\setminus \Delta (\Sigma (B)) = \{\{1,3,4\}\},$ for both versions of shifting. Surprisingly,
we even get that
\begin{equation}\label{neq}
\Delta (\Sigma (B))<_{L}\Delta (\Sigma (\Delta(B))),
\end{equation}
where the lexicographic partial order on simplicial complexes is
defined (as in \cite{skira}) by: $K\leq_{L} L$ iff for all $r>0$
the lexicographically first $r$-face in $K\triangle L$ (if exists)
belongs to $K$.
\end{ex}
\begin{conj}\label{suspension}
For any two simplicial complexes $K$ and $L$: $$\Delta (K*L)\leq_{L}\Delta
(\Delta K* \Delta L).$$
\end{conj}
Very recently Satoshi Murai announced a proof of  Conjecture \ref{suspension} for the exterior case \cite{Murai-Join}.

\begin{conj}\label{suspension-topo.}(Topological invariance.)
Let $K_1$ and $K_2$ be triangulations of the same topological
space. Then $\Delta (\Sigma (K_1))<_{L}\Delta (\Sigma
(\Delta(K_1)))$ iff $\Delta (\Sigma (K_2))<_{L}\Delta (\Sigma
(\Delta(K_2)))$.
\end{conj}
It would be interesting to find out when equation (\ref{eq**})
holds.
If both $K$ and $L$ are shifted, it trivially holds as
$\Delta^{2}=\Delta$. By the remark to Corollary \ref{prop.11cone}
it also holds if $K$, say, is a complete simplicial complex.

\section{Open problems}
\begin{enumerate}
\item
(Special case of Problem \ref{U problem}.) Let $L$ be a subcomplex
of $K$. What are the relations between $\Delta(K)$, $\Delta(L)$ and
$\Delta(K\cup_L \rm{Cone}(L))$?

\item \label{U3 problem}
Characterize the (face,Betti)-vectors of quadruples $(K,L,K\cup
L,K\cap L)$ (by using shifting).

Some necessary conditions on such vectors are given by the
(f,$\beta$)-vector characterization for chains of complexes w.r.t.
inclusion (Bj\"{o}rner and Kalai \cite{BK}, Duval \cite{Duval30}),
others (linear inequalities) are given by the Mayer-Vietoris exact
sequence.

\item Prove equation (\ref{additive-formula}) is the symmetric case: $$|I_{S}^{j}\cap \Delta^s(K\cup L)|=|I_{S}^{j}\cap
\Delta^s(K)|+|I_{S}^{j}\cap \Delta^s(L)|-|I_{S}^{j}\cap \Delta^s(K\cap
L)|$$
for every simplicial complexes $K,L$ and all sets $S$ and positive integers $j$.

\item Prove K\"{u}nneth theorem with field coefficients using shifting.

\item Can one recover (part of) the cohomology ring of $K$ by knowing the
shifting of suitable complexes related to $K$?

\item (\cite{skira}, Problem 15) Is algebraic shifting a functor? (onto a category with a useful
set of maps).
\end{enumerate}


\chapter{Algebraic Shifting and Rigidity of Graphs}\label{chapter:rigidity}
\section{Basics of rigidity theory of graphs}
Asimov and Roth introduced the concept of generic rigidity of graphs
\cite{Asi-Roth1,Asi-Roth2}. The presentation of rigidity here is
based mainly on the one in Kalai \cite{Kalai-LBT}.

Let $G=(V,E)$ be a simple graph. A map $f:V\rightarrow \mathbb{R}^{d}$ is
called a $d$-\emph{embedding}. It is \emph{rigid} if any small enough
perturbation of it which preserves the lengths of the edges is
induced by an isometry of $\mathbb{R}^{d}$. Formally, $f$ is
called rigid if there exists an $\varepsilon>0$ such that if
$g:V\rightarrow \mathbb{R}^{d}$ satisfies
$d(f(v),g(v))<\varepsilon$ for every $v\in V$ and
$d(g(u),g(w))=d(f(u),f(w))$ for every $\{u,w\}\in E$, then
$d(g(u),g(w))=d(f(u),f(w))$ for every $u,w\in V$ (where $d(a,b)$
denotes the Euclidean distance between the points $a$ and $b$).

$G$ is called \emph{generically d-rigid} if the set of its rigid
$d$-embeddings is open and dense in the topological vector space
of all of its $d$-embeddings.

Given a $d$-embedding $f:V\rightarrow \mathbb{R}^{d}$, a \emph{stress}
w.r.t. $f$ is a function $w:E\rightarrow \mathbb{R}$ such that for
every vertex $v\in V$
$$\sum_{u:\{v,u\}\in E}w(\{v,u\})(f(v)-f(u)) =0.$$
$G$ is called \emph{generically d-stress free} if the set of its
$d$-embeddings which have a unique stress ($w=0$) is open and dense
in the space of all of its $d$-embeddings.

Rigidity and stress freeness can be related as follows: Let $V=[n]$,
and let $\Rig(G,f)$ be the $dn\times |E|$ matrix associated with a
$d$-embedding $f$ of $V(G)$ defined as follows: for its column
corresponding to $\{v<u\}\in E$ put the vector $f(v)-f(u)$ (resp.
$f(u)-f(v)$) at the entries of the $d$ rows corresponding to $v$
(resp. $u$) and zero otherwise. $G$ is generically $d$-stress free
if $\Ker(\Rig(G,f))=0$ for a generic $f$ (i.e. for an open and dense
set of embeddings). $G$ is generically $d$-rigid if
$\im(\Rig(G,f))=\im(\Rig(K_{V},f)$ for a generic $f$, where $K_{V}$
is the complete graph on $V=V(G)$. The dimensions of the kernel and
image of $\Rig(G,f)$ are independent of the generic $f$ we choose;
we call $R(G)=\Rig(G,f)$ the \emph{rigidity matrix} of $G$.

$\im(\Rig(K_{V},f))$ can be described by the following linear
equations:
$(v_{1},..,v_{d})\in \bigoplus_{i=1}^{d}\mathbb{R}^{n}$ belongs to
$\im(\Rig(K_{V},f))$ iff
\begin{equation} \label{ImRig1}
\forall 1\leq i \neq j \leq d\ \   <f_{i},v_{j}>=<f_{j},v_{i}>
\end{equation}
\begin{equation}\label{ImRig2}
\forall 1\leq i\leq d\ \  <e,v_{i}>=0
\end{equation}
where $e$ is the all ones vector and $f_{i}$ is the vector of the
$i$-th coordinate of the $f(v)$'s, $v\in V$. From this description
it is clear that $\rank(\Rig(K_{V},f))=dn-$${d+1}\choose{2}$ (see
Asimov and Roth \cite{Asi-Roth1} for more details).

Gluck \cite{Gluck} has proven that
\begin{thm}(Gluck)\label{Gluck}
The graph of a triangulated $2$-sphere is
generically $3$-rigid. Equivalently, planar graphs are generically
$3$-stress free.
\end{thm}
The equivalence follows from the facts that every triangulated
$2$-sphere with $n$ vertices has exactly $3n-6$ edges (hence it is
generically $3$-rigid iff it is generically $3$-stress free), and
that every planar graph is a subgraph of a triangulated $2$-sphere.
Gluck's proof is based on two classical theorems: one is Cauchy's
rigidity theorem (e.g. \cite{Cromwell-CauchyThm}), which states that
any combinatorial isomorphism between two convex $3$-polytopes which
induces an isometry on their boundaries is actually induced by an
isometry of $\mathbb{R}^3$; the other is Steinitz's theorem
\cite{Steinitz}, which asserts that any polyhedral $2$-sphere is
combinatorially isomorphic to the boundary complex of some convex
$3$-polytope. Whiteley \cite{Wh} has found a proof of Gluck's
theorem which avoids convexity, based on vertex splitting. We
summarize it below.
\begin{lem}\label{Whiteley}(Whiteley)
Let $G'$ be obtained from a graph $G$ by contracting an edge
$\{u,v\}$.

 (a)If $u,v$ have at least $d-1$ common neighbors and
$G'$ is generically $d$-rigid, then $G$ is generically $d$-rigid.

(b)If $u,v$ have at most $d-1$ common neighbors and $G'$ is
generically $d$-stress free, then $G$ is generically $d$-stress
free.
\end{lem}
Lemma \ref{Whiteley} gives an alternative proof of Gluck's theorem:
starting with a triangulated $2$-sphere, repeatedly contract edges
with exactly $2$ common neighbors until the $1$-skeleton of a
tetrahedron is reached (it is not difficult to show that this is
always possible). By Theorem \ref{Whiteley}(a) it is enough to show
that the $1$-skeleton of a tetrahedron is generically $3$-rigid, as
is well known. (By definition, the graph of a simplex is generically
$d$-rigid for every $d$, and this is also true if defining rigidity
via isometries of $\mathbb{R}^d$ \cite{Asi-Roth1}).

We will need the following gluing lemma, due of Asimov and Roth \cite{Asi-Roth2}.
\begin{lem}\label{weak-clique-sum,Asimov}(Asimov and Roth)
(1) Let $G_1$ and $G_2$ be generically $d$-rigid graphs. If $G_1\cap
G_2$ contains at least $d$ vertices, then $G_1\cup G_2$ is
generically $d$-rigid.

(2) Let $G_{i}=(V_{i},E_{i})$ be generically $k$-stress free graphs, $i=1,2$ such that
$G_{1}\cap G_{2}$ is generically $k$-rigid. Then $G_{1}\cup G_{2}$ is
generically $k$-stress free.
\end{lem}
\section{Rigidity and symmetric shifting}\label{sec:Lee}
Let $G$ be the $1$-skeleton of a $(d-1)$-dimensional simplicial
complex $K$ with vertex set $[n]$. We define $d$ generic degree-one
elements in the polynomial ring $A=\mathbb{R}[x_1,..,x_n]$ as
follows: $\theta_i=\sum_{v\in [n]}f(v)_i x_v$ where $f(v)_i$ is the
projection of $f(v)$ on the $i$-th coordinate, $1\leq i\leq d$. Then
the sequence $\Theta=(\theta_1,..,\theta_d)$ is an l.s.o.p. for the
face ring $\mathbb{R}[K]=A/I_{K}$ ($I_K$ is the ideal in $A$
generated by the monomials whose support is not an element of $K$).
Let $H(K)=\mathbb{R}[K]/(\Theta)=H(K)_0\oplus H(K)_1\oplus...$ where
$(\Theta)$ is the ideal in $A$ generated by the elements of $\Theta$
and the grading is induced by the degree grading in $A$. Consider
the multiplication map $\omega: H(K)_1\longrightarrow H(K)_2$,
$m\rightarrow \omega m$ where $\omega=\sum_{v\in [n]}x_v$. Lee
\cite{Lee} proved that
\begin{equation}\label{eqLee}
\ddim_{\mathbb{R}} \Ker(\Rig(G,f)) = \ddim_{\mathbb{R}}
H(K)_2-\ddim_{\mathbb{R}} \omega(H(K)_1).
\end{equation}
Assume that $G$ is generically $d$-rigid. Then $\ddim_{\mathbb{R}}
\Ker(\Rig(G,f)) = f_1(K)-\rank(\Rig(K_{V},f)) = g_2(K) =
\ddim_{\mathbb{R}} H(K)_2-\ddim_{\mathbb{R}} H(K)_1$. Combining with
(\ref{eqLee}), the map $\omega$ is injective, and hence
$\ddim_{\mathbb{R}} (H(K)/(\omega))_i = g_i(K)$ for $i=2$; clearly
this holds for $i=0,1$ as well. Hence $(g_o(K),g_1(K),g_2(K))$ is an
$M$-sequence, i.e. the Hilbert function of a standard ring - the
sequence counting the dimensions of the graded pieces of the ring by
their degree. To summarize:
\begin{thm}\label{thm:Lee}(Lee \cite{Lee})
If a simplicial complex $K$ has a generically $(\ddim K+1)$-rigid
$1$-skeleton, then multiplication by a generic degree $1$ element
$\omega: H_1(K)\rightarrow  H_2(K)$ is injective. In particular,
$(g_o(K),g_1(K),g_2(K))$ is an $M$-sequence.
\end{thm}
Note that if the multiplication $\omega: H(K)_1\rightarrow H(K)_2$
is injective, then so is the multiplication by a generic monomial of
degree $1$, $\theta_{d+1}$, and vice versa. In terms of $GIN$, this
means that $\theta_{d+1}\theta_n\in GIN(K)$, equivalently that
$\{d,n\}\in \Delta^s(K)$. To summarize, $K$ is generically $d$-rigid
iff  $\{d,n\}\in \Delta^s(K)$. Similarly, $K$ is generically
$d$-stress free iff  $\{d+1,d+2\}\notin \Delta^s(K)$ (iff $\omega:
H(K)_1\rightarrow H(K)_2$ is onto).

\section{Hyperconnectivity and exterior shifting}
We will describe now an exterior
analogue of rigidity, namely Kalai's notion of hyperconnectivity
\cite{56}.
We keep the notation from the previous section and from Chapter \ref{chapter:BasicDef&Concepts}, and follow the
presentation in \cite{56}.

Consider the map
$$f(d,j,K): \bigwedge^{j+1}(K)\rightarrow \bigoplus_{1}^{d}\bigwedge^{j}(K) \ \ x\mapsto (f_{1}\lfloor
x,...,f_{d}\lfloor x).$$
The dimension of its kernel equals
$|\{S\in \Delta^{e}K: |S|=j+1, S\cap [d]=\emptyset\}|$;
it follows from Propositions \ref{prop.1} (with $R$ a singleton in $[d]$) and \ref{prop.*}.
 Kalai \cite{56} called a graph $G$
$d$-\emph{hyperconnected} if $\im(f(d,1,G))=\im(f(d,1,K_{V(G)}))$,
and $d$-\emph{acyclic} if $\Ker(f(d,1,G))=0$. With this terminology,
$G$ is $d$-acyclic iff $\{d+1,d+2\}\notin \Delta^{e}(G)$, and is
$d$-hyperconnected iff $\{d,n\}\in \Delta^{e}(G)$, where $n=|V(G)|$.

We shall prove now an exterior analogue of Lemma \ref{Whiteley}:
\begin{lem}\label{extWHiteley}
If $G'$ is obtained from $G$ by contracting an edge which belongs
to at most $d-1$ triangles, and $G'$ is $d$-acyclic, then so is
$G$.
\end{lem}
$Proof$: Let $\{v,u\}$ be the edge we contract, $u\mapsto v$. Consider the
$dn\times |E|$ matrix $A$ of the map $f(d,1,G)$ w.r.t. the
standard basis, where $f_{i}=\sum_{j=1}^{n}\alpha_{ij}e_{j}$,
$n=|V|$: for its column corresponding to $\{v<u\}\in E$ put the
vector $(\alpha_{1u},..,\alpha_{du})^{T}$ (resp. $-(\alpha_{1v},..,\alpha_{dv})^{T}$)
at the entries of the rows corresponding to $v$ (resp. $u$) and
zero otherwise.

Now replace in $A$ each $\alpha_{iv}$ with $\alpha_{iu}$ to obtain a
new matrix $\hat{A}$. It is enough to show that the columns of
$\hat{A}$ are independent: As the set of $dn\times |E|$ matrices
with independent columns is open (in the Euclidian topology), by
perturbing the $\alpha_{iu}$'s in the places where $\hat{A}$ differs
from $A$, we may obtain new generic $\alpha_{iv}$'s forming a matrix
with independent columns. As for every generic choice of $f_{i}$'s,
the map $f(d,1,G)$ has the same rank, we would conclude that the
columns of $A$ are independent as well.

Suppose that some linear combination of the columns of $\hat{A}$
equals zero. Let $\bar{A}$ be obtained from $\hat{A}$ by adding the
rows of $v$ to the corresponding rows of $u$, and deleting the rows
of $v$. Thus, a linear combination of the columns of $\bar{A}$ with
the same coefficients  also equals zero. $\bar{A}$ is obtained from
the matrix of $f(d,1,G')$ by adding a zero column (for the edge
$\{v,u\}$) and doubling the columns $\{v,w\}$ which correspond to
common neighbors $w$ of $v$ and $u$ in $G$. As $\Ker(f(d,1,G'))=0$,
apart from the above mentioned columns the rest have coefficient
zero, and pairs of columns we doubled have opposite sign. Let us
look at the submatrix of $\hat{A}$ consisting of the 'doubled'
columns with vertex $v$ and the column of $\{v,u\}$, restricted to
the rows of $v$: it has generic coefficients, $d$ rows and at most
$d$ columns, hence its columns are independent. Thus, all
coefficients in the above linear combination are zero. $\square$

Similarly,
\begin{lem}\label{extWHiteleyHyper}
If $G'$ is obtained from $G$ by contracting an edge which belongs to
at least $d-1$ triangles, and $G'$ is $d$-hyperconnected, then so is
$G$. $\square$
\end{lem}

We need the following exterior analogue of Lemma
\ref{weak-clique-sum,Asimov}:
\begin{lem}\label{weak-clique-sum,Kalai}(Kalai \cite{56}, Theorem 4.4)
Let $G_{i}=(V_{i},E_{i})$ be $k$-acyclic graphs, $i=1,2$ such that
$G_{1}\cap G_{2}$ is $k$-hyperconnected. Then $G_{1}\cup G_{2}$ is
$k$-acyclic.

Similarly, if $G_{i}=(V_{i},E_{i})$ are $k$-hyperconnected graphs, $i=1,2$ such that
$|G_{1}\cap G_{2}|\geq k$, then $G_{1}\cup G_{2}$ is
$k$-hyperconnected.
\end{lem}
We will also need the easy fact that the graph of a $k$-simplex is
$d$-hyperconnected for every $k\geq d$.

\section{Minimal cycle complexes}\label{sec:Fogelsanger}
We shall need the concept
of minimal cycle complexes, introduced by Fogelsanger
\cite{Fogelsanger}. We summarize his theory below.

Fix a field $k$ (or more generally, any abelian group) and consider
the formal chain complex on a ground set $[n]$, $C=(\oplus\{kT
:T\subseteq [n]\}, \partial)$, where $\partial(1T)=\sum_{t\in T}
\sgn(t,T)T\setminus \{t\}$ and $\sgn(t,T)=(-1)^{|\{s\in T: s<t\}|}$.
Define \emph{subchain}, \emph{minimal d-cycle} and \emph{minimal
d-cycle complex} as follows: $c'=\sum\{b_TT: T\subseteq [n],
|T|=d+1\}$ is a \emph{subchain} of a $d$-chain $c=\sum\{a_TT:
T\subseteq [n], |T|=d+1\}$ iff for every such $T$, $b_T=a_T$ or
$b_T=0$. A $d$-chain $c$ is a $d-cycle$ if $\partial(c)=0$, and is a
\emph{minimal d-cycle} if its only subchains which are cycles are
$c$ and $0$. A simplicial complex $K$ which is spanned by the
support of a minimal $d$-cycle is called a \emph{minimal d-cycle
complex} (over $k$), i.e. $K=\{S: \exists T\ S\subseteq T, a_T\neq
0\}$ for some minimal $d$-cycle $c$ as above. For example,
triangulations of connected manifolds without boundary are minimal
cycle complexes - fix $k=\mathbb{Z}_2$ and let the cycle be the sum
of all facets.

The following is the main result in Fogelsanger's thesis \cite{Fogelsanger}.
\begin{thm}\label{thmFog}(Fogelsanger)
For $d\geq 3$, every minimal $(d-1)$-cycle complex has a generically $d$-rigid $1$-skeleton.
\end{thm}

His proof relies on the following three properties of rigidity
solely: Lemmata \ref{Whiteley} and \ref{weak-clique-sum,Asimov} and
the fact that the graph of a $d$-simplex is generically $d$-rigid.
As these three properties hold for hyperconnectivity as well (see
Lemmata \ref{extWHiteleyHyper} and \ref{weak-clique-sum,Kalai} and
the fact that the graph of a $d$-simplex is $d$-hyperconnected),
Theorem \ref{thmFog} holds for hyperconnectivity as well. In terms
of algebraic shifting this means

\begin{thm}\label{thmFogShifting}
For $n>d\geq 3$ and every minimal $(d-1)$-cycle complex $K$ on $n$
vertices, over the field $\mathbb{R}$, $\{d,n\}\in\Delta(K)$ holds
for both versions of algebraic shifting. $\square$
\end{thm}

\section{Rigidity and doubly Cohen-Macaulay complexes}\label{sec:2CM}
\begin{de}
A simplicial complex $K$ is \emph{doubly Cohen-Macaulay} ($2-CM$ in short) over a fixed field $k$, if it
is Cohen-Macaulay and for every vertex $v\in K$, $K\setminus v$ is
Cohen-Macaulay of the same dimension as $K$.
\end{de}
Here $K\setminus v$ is the simplicial complex $\{T\in K: v\notin T\}$. By a
theorem of Reisner \cite{Reisner}, a simplicial complex $L$ is
Cohen-Macaulay over $k$ iff it is pure and for every face $T\in L$
(including the empty set) and every $i<\ddim(\lk(T,L)$,
$\tilde{H}_i((\lk(T,L);k)=0$.

For example, triangulated spheres are $2$-CM, triangulated balls are
not. A \emph{homology sphere over $k$} is a simplicial complex $K$
such that for every $F\in K$ and every $i$
$\tilde{H}_i((\lk(F,K);k)\cong \tilde{H}_i((S^{\ddim(\lk(F,K)};k)$
where $S^d$ is the $d$-dimensional sphere. Based on the fact that
homology spheres are $2$-CM and that the $g$-vector of some other
classes of $2$-CM complexes is known to be an $M$-sequence (e.g.
\cite{Swartz}), Bj\"{o}rner and Swartz \cite{Swartz} recently
suspected that
\begin{conj}\label{g-conj 2CM}(\cite{Swartz}, a weakening of Problem 4.2.)
The $g$-vector of any $2$-CM complex is an $M$-sequence.
\end{conj}
We prove a first step in this direction, namely:
\begin{thm}\label{g_2-2CM}
Let $K$ be a $(d-1)$-dimensional $2$-CM simplicial complex (over
some field) where $d\geq 4$. Then $(g_0(K),g_1(K),g_2(K))$ is an
$M$-sequence.
\end{thm}
This theorem follows from the following theorem, combined with Theorem \ref{thm:Lee}.
\begin{thm}\label{thmRigid-2CM}
Let $K$ be a $(d-1)$-dimensional $2$-CM simplicial complex (over
some field) where $d\geq 3$. Then $K$ has a generically $d$-rigid
$1$-skeleton.
\end{thm}
Kalai \cite{Kalai-LBT} showed that if a simplicial complex $K$ of
dimension $\geq 2$ satisfies the following conditions then it
satisfies Barnette's lower bound inequalities:

(a) $K$ has a generically $(\ddim(K)+1)$-rigid $1$-skeleton.

(b) For each face $F$ of $K$ of codimension $>2$, its link
$\lk(F,K)$ has a generically $(\ddim(\lk(F,K))+1)$-rigid $1$-skeleton.

(c) For each face $F$ of $K$ of codimension $2$, its link
$\lk(F,K)$ (which is a graph) has at least as many edges as
vertices.

Kalai used this observation to prove that Barnette's inequalities
hold for a large class of simplicial complexes.

Observe that the link of a vertex in a $2$-CM simplicial complex
is $2$-CM, and that a $2$-CM graph is $2$-connected. Combining it
with Theorem \ref{thmRigid-2CM} and the above result of Kalai we
conclude:
\begin{cor}\label{LBT-2CM}
Let $K$ be a $(d-1)$-dimensional $2$-CM simplicial complex where
$d\geq 3$. For all $0\leq i\leq d-1$ $f_i(K)\geq f_i(n,d)$ where
$f_i(n,d)$ is the number of $i$-faces in a (equivalently every)
stacked $d$-polytope on $n$ vertices. (Explicitly,
$f_{d-1}(n,d)=(d-1)n-(d+1)(d-2)$ and
$f_{i}(n,d)=\binom{d}{i}n-\binom{d+1}{i+1}i$ for $1\leq i\leq
d-2$.) $\square$
\end{cor}
Theorem \ref{thmRigid-2CM} is proved by decomposing $K$ into a union
of minimal $(d-1)$-cycle complexes (defined in Section
\ref{sec:Fogelsanger}). Each of these pieces has a generically
$d$-rigid $1$-skeleton by Theorem \ref{thmFog}, and the
decomposition is such that gluing the pieces together results in a
complex with a generically $d$-rigid $1$-skeleton. The decomposition
is detailed in Theorem \ref{structureLemma} below. Its proof is by
induction on $\ddim(K)$. Let us first consider the case where $K$ is
$1$-dimensional.

A (simple finite) graph is $2$-\emph{connected} if after a
deletion of any vertex from it, the remaining graph is connected
and non trivial (i.e. is not a single vertex nor empty). Note that
a graph is $2$-CM iff it is $2$-connected.
\begin{lem}\label{graphLemma}
A graph $G$ is $2$-connected iff there exists a decomposition
$G=\cup^m_{i=1} C_i$ such that each $C_i$ is a simple cycle and
for every $1<i\leq m$, $C_i\cap(\cup_{j<i}C_j)$ contains an edge.

Moreover, for each $i_0\in [m]$ the $C_i$'s can be reordered by a
permutation $\sigma:[m]\rightarrow[m]$ such that
$\sigma^{-1}(1)=i_0$ and for every $i>1$,
$C_{\sigma^{-1}(i)}\cap(\cup_{j<i}C_{\sigma^{-1}(j)})$ contains an
edge.
\end{lem}
$Proof$: Whitney \cite{Whitney21} showed that a graph $G$ is
$2$-connected iff it has an open ear decomposition, i.e. there
exists a decomposition $G=\cup^m_{i=0} P_i$ such that each $P_i$
is a simple open path, $P_0$ is an edge, $P_0\cup P_1$ is a simple
cycle and for every $1<i\leq m$ $P_i\cap(\cup_{j<i}P_j)$ equals
the $2$ end vertices of $P_i$.

Assume that $G$ is $2$-connected and consider an open ear
decomposition as above. Let $C_1=P_0\cup P_1$. For $i>1$ choose a
simple path $\tilde{P}_i$ in $\cup_{j<i}P_j$ that connects the $2$
end vertices of $P_i$, and let $C_i=P_i\cup \tilde{P_i}$.
$(C_1,...,C_m)$ is the desired decomposition sequence of $G$.

Let $C$ be the graph whose vertices are the $C_i$'s and two of
them are neighbors iff they have an edge in common. Thus, $C$ is
connected, and hence the 'Moreover' part of the Lemma is proved.

The other implication, that such a decomposition implies
$2$-connectivity, will not be used in the sequel, and its proof is
omitted. $\square$
\\

For the induction step we need the following cone lemma. For $v$ a
vertex not in the support of a $(d-1)$-chain $c$, let $v*c$ denote
the following $d$-chain: if $c=\sum\{a_T T: v\notin T\subseteq [n],
|T|=d\}$ where $a_T\in k$ for all $T$, then $v*c=\sum\{\sgn(v,T)a_T
T\cup\{v\}: v\notin T\subseteq [n], |T|=d\}$.

\begin{lem}\label{sphereLemma}
Let $s$ be a minimal $(d-1)$-cycle and let $c$ be a minimal
$d$-chain such that $\partial(c)=s$, i.e. $c$ has no proper
subchain $c'$ such that $\partial(c')=s$. For $v$ a vertex not in
any face in $\supp(c)$,the support of $c$, define
$\tilde{s}=c-v*s$. Then $\tilde{s}$ is a minimal $d$-cycle.
\end{lem}
$Proof$:
$\partial(\tilde{s})=\partial(c)-\partial(v*s)=s-(s-v*\partial(s))=0$
hence $\tilde{s}$ is a $d$-cycle. To show that it is minimal, let
$\hat{s}$ be a subchain of $\tilde{s}$ such that
$\partial(\hat{s})=0$. Note that $\supp(c)\cap
\supp(v*s)=\emptyset$.
\newline Case 1: $v$ is contained in a face in $\supp(\hat{s})$. By the minimality
of $s$, $\supp(v*s)\subseteq \supp(\hat{s})$. Thus, by the
minimality of $c$ also $\supp(c)\subseteq \supp(\hat{s})$ and hence
$\hat{s}=\tilde{s}$.
\newline Case 2: $v$ is not contained in any face in $\supp(\hat{s})$. Thus,
$\supp(\hat{s})\subseteq \supp(c)$. As $\partial(\hat{s})=0$ then
$\partial(c-\hat{s})=s$. The minimality of $c$ implies
$\hat{s}=0$. $\square$

\begin{thm}\label{structureLemma}
Let $K$ be a $d$-dimensional $2$-CM simplicial complex over a
field $k$ ($d\geq 1$). Then there exists a decomposition
$K=\cup^m_{i=1} S_i$ such that each $S_i$ is a minimal $d$-cycle
complex over $k$ and for every $i>1$, $S_i\cap(\cup_{j<i}S_j)$
contains a $d$-face.

Moreover, for each $i_0\in [m]$ the $S_i$'s can be reordered by a
permutation $\sigma:[m]\rightarrow[m]$ such that
$\sigma^{-1}(1)=i_0$ and for every $i>1$,
$S_{\sigma^{-1}(i)}\cap(\cup_{j<i}S_{\sigma^{-1}(j)})$ contains a
$d$-face.
\end{thm}
$proof$: The proof is by induction on $d$. For $d=1$, by Lemma
\ref{graphLemma} $K=\cup^{m(K)}_{i=1} C_i$ such that each $C_i$ is
a simple cycle and for every $i>1$ $C_i\cap(\cup_{j<i}C_j)$
contains an edge. Define $s_i=\sum\{\sgn_e(i)e: e\in (C_i)_1\}$,
then $s_i$ is a minimal $1$-cycle (orient the edges properly:
$\sgn_e(i)$ equals $1$ or $-1$ accordingly) whose support spans
the simplicial complex $C_i$. Moreover, by Lemma \ref{graphLemma}
each $C_{i_0}$, $i_0\in [m(K)]$, can be chosen to be the first in
such a decomposition sequence.

For $d>1$, note that the link of every vertex in a $2$-CM
simplicial complex is $2$-CM. For a vertex $v\in K$, as $\lk(v,K)$
is $2$-CM then by the induction hypothesis
$\lk(v,K)=\cup^{m(v)}_{i=1} C_i$ such that each $C_i$ is a minimal
$(d-1)$-cycle complex and for every $i>1$ $C_i\cap(\cup_{j<i}C_j)$
contains a $(d-1)$-face. Let $s_i$ be a minimal $(d-1)$-cycle
whose support spans $C_i$. As $K\setminus v$ is CM of dimension $d$,
$\tilde{H}_{d-1}(K\setminus v;k)=0$. Hence there exists a $d$-chain $c$
such that $\partial(c)=s_i$ and $\supp(c)\subseteq K\setminus v$.

Take $c_i$ to be such a chain with a support of minimal
cardinality. By Lemma \ref{sphereLemma}, $\tilde{s_i}=c_i-v*s_i$
is a minimal $d$-cycle. Let $S_i(v)$ by the simplicial complex
spanned by $\supp(\tilde{s_i})$; it is a minimal $d$-cycle complex.
By the induction hypothesis, for every $i>1$
$S_i(v)\cap(\cup_{j<i}S_j(v))$ contains a $d$-face (containing
$v$). Thus, $K(v):=\cup_{j=1}^{m(v)}S_j(v)$ has the desired
decomposition for every $v\in K$. $K=\cup_{v\in K_0}K(v)$ as
$\st(v,K)\subseteq K(v)$ for every $v$.

Let $v$ be any vertex of $K$. Since the $1$-skeleton of $K$ is
connected, we can order the vertices of $K$ such that $v_1=v$ and
for every $i>1$ $v_i$ is a neighbor of some $v_j$ where $1\leq
j<i$. Let $v_{l(i)}$ be such a neighbor of $v_i$. By the induction
hypothesis we can order the $S_j(v_i)$'s such that $S_1(v_i)$ will
contain $v_{l(i)}$, and hence, as $K$ is pure, will contain a
$d$-face which appears in $K(v_{l(i)})$ (this face contains the
edge $\{v_i,v_{l(i)}\}$). The resulting decomposition sequence
$(S_1(v_1),..,S_{m(v_1)}(v_1),S_1(v_2),..,S_{m(v_n)}(v_n))$ is as
desired.

Moreover, every $S_j(v_{i_0})$ where $i_0\in [n]$ and $j\in
[m(v_{i_0})]$ can be chosen to be the first in such a
decomposition sequence. Indeed, by the induction hypothesis
$S_j(v_{i_0})$ can be the first in the decomposition sequence of
$K(v_{i_0})$, and as mentioned before, the connectivity of the
$1$-skeleton of $K$ guarantees that each such prefix
$(S_1(v_{i_0}),..,S_{m(v_{i_0})}(v_{i_0}))$ can be completed to a
decomposition sequence of $K$ on the same $S_j(v_i)$'s. $\square$
\\
\\
$proof\ of\ Theorem\ \ref{thmRigid-2CM}$: Consider a decomposition
sequence of $K$ as guaranteed by Theorem \ref{structureLemma},
$K=\cup^m_{i=1} S_i$. By Theorem \ref{thmFog} each $S_i$ has a
generically $d$-rigid $1$-skeleton. By Lemma \ref{weak-clique-sum,Asimov}
for all $2\leq i\leq m$ $\cup^i_{j=1} S_j$ has a generically
$d$-rigid $1$-skeleton, in particular $K$ has a generically
$d$-rigid $1$-skeleton ($i=m$). $\square$
\newline

Theorem \ref{thmRigid-2CM} follows also from the following
corollary combined with Theorem \ref{thmFog}.
\begin{cor}\label{cor2CM->minCycle}
Let $K$ be a $d$-dimensional $2$-CM simplicial complex over a
field $k$ ($d\geq 1$). Then $K$ is a minimal cycle complex over
the Abelian group $\tilde{k}= k(x_1,x_2,...)$ whose elements are
finite linear combinations of the (variables) $x_i$'s with
coefficients in $k$.
\end{cor}
$Proof$: Consider a decomposition $K=\cup^m_{i=1} S_i$ as
guaranteed by Theorem \ref{structureLemma}, where
$S_i=\overline{\supp(c_i)}$, the closure w.r.t. inclusion of
$\supp(c_i)$, for some minimal $d$-cycle $c_i$ over $k$. Define
$\tilde{c_i}=x_i c_i$, thus $\tilde{c_i}$ is a minimal cycle over
$\tilde{k}$. Define $\tilde{c}=\sum_{i=1}^m \tilde{c_i}$. Clearly
$\tilde{c}$ is a cycle over $\tilde{k}$ whose support spans $K$.
It remains to show that $\tilde{c}$ is minimal. Let $\tilde{c}'$
be a subchain of $\tilde{c}$ which is a cycle, $\tilde{c}'\neq
\tilde{c}$. We need to show that $\tilde{c}'=0$. Denote by
$\tilde{\alpha_T}$ ($\tilde{\alpha_T}'$) the coefficient of the
set $T$ in $\tilde{c}$ ($\tilde{c}'$) and by $\tilde{\alpha_T}(i)$
the coefficient of the set $T$ in $\tilde{c_i}$. If
$\tilde{\alpha_T}'=0$ then for every $i$ such that
$\tilde{\alpha_T}(i)\neq 0$, the minimality of $\tilde{c_i}$
implies that $\tilde{\alpha_F}'=0$ whenever
$\tilde{\alpha_F}(i)\neq 0$. By assumption, there exists a set
$T_0$ such that $\tilde{\alpha_{T_0}}'=0\neq
\tilde{\alpha_{T_0}}$. In particular, there exists an index $i_0$
such that $\tilde{\alpha_{T_0}}(i_0)\neq 0$, hence
$\tilde{\alpha_F}'=0$ whenever $\tilde{\alpha_F}(i_0)\neq 0$. As
$S_{i_0}\cap(\cup_{j<i_0}S_j)$ contains a $d$-face in case
$i_0>1$, repeated application of the above argument implies
$\tilde{\alpha_F}'=0$ whenever $\tilde{\alpha_F}(1)\neq 0$.
Repeated application of the fact that $S_{i}\cap(\cup_{j<i}S_j)$
contains a $d$-face for $i=2,3,..$ and of the above argument shows
that $\tilde{\alpha_F}'=0$ whenever $\tilde{\alpha_F}(i)\neq 0$
for some $1\leq i\leq m$, i.e. $\tilde{c}'=0$. $\square$
\newline

A pure simplicial complex has a \emph{nowhere zero flow} if there
is an assignment of integer non-zero wights to all of its facets
which forms a $\mathbb{Z}$-cycle. This generalizes the definition
of a nowhere zero flow for graphs (e.g. \cite{SeymourNowhereZeroFlow} for a
survey).
\begin{cor}\label{corZ2CM->nowhereZeroFlow}
Let $K$ be a $d$-dimensional $2$-CM simplicial complex over
$\mathbb{Q}$ ($d\geq 1$). Then $K$ has a nowhere zero flow.
\end{cor}
$Proof$: Consider a decomposition $K=\cup^m_{i=1} S_i$ as
guaranteed by Theorem \ref{structureLemma}. Multiplying by a
common denominator, we may assume that each
$S_i=\overline{supp(c_i)}$ for some minimal $d$-cycle $c_i$ over
$\mathbb{Z}$ (instead of just over $\mathbb{Q}$). Let $N$ be the
maximal $|\alpha|$ over all nonzero coefficients $\alpha$ of the
$c_i$'s, $1\leq i\leq m$. Let $\tilde{c}=\sum_{i=1}^m
(N^m)^{i}c_i$. $\tilde{c}$ is a nowhere zero flow for $K$; we omit
the details. $\square$

\section{Shifting and minors of graphs}\label{sec:ShiftingTellsMinors}
\subsection{Shifting can tell minors}
Inspired by Lemma \ref{Whiteley}, we will show now how shifting can tell graph minors.
\begin{thm}\label{mainThm}
The following holds for symmetric and exterior shifting: for every
$2\leq r \leq 6$ and every graph $G$, if $\{r-1,r\} \in \Delta(G)$
then $G$ has a $K_r$ minor.
\end{thm}
Note that the case $r=5$ strengthens Gluck's Theorem \ref{Gluck},
via the interpretation of rigidity in terms of symmetric shifting
(see Section \ref{sec:Lee}).

The proof is by induction on the number of vertices, based on
contracting edges satisfying the conditions of Lemma \ref{Whiteley}.
We make an essential use of Mader's theorem \cite{Mader} which gives
an upper bound  $(r-2)n-$ $r-1 \choose 2$ on the number of edges in
a $K_{r}$-minor free graph with $n$ vertices, for $r\leq 7$. Indeed,
Theorem \ref{mainThm} can be regarded as a strengthening of Mader's
theorem, as $\{l+1,l+2\} \notin \Delta(G)$ implies having at most
$ln-$ $l+1 \choose 2$ edges, as is clear from the facts
$f(G)=f(\Delta(G))$ and $\Delta(G)\subseteq \rm{span}(\{d,n\})$ (as
it is shifted). This also shows that Theorem \ref{mainThm} fails for
$r\geq 8$, as is demonstrated for $r=8$ by $K_{2,2,2,2,2}$, and for
$r>8$ by repeatedly coning over the resulted graph for a smaller $r$
(e.g. \cite{Song}). It would be interesting to find a proof of
Theorem \ref{mainThm} that avoids using Mader's theorem, and derive
Mader's theorem as a corollary.

A graph is \emph{linklessly embeddable} if there exists an embedding
of it in $\mathbb{R}^3$ (where vertices and edges have disjoint
images) such that every two disjoint cycles of it are unlinked
closed curves in $\mathbb{R}^3$. As such graph is $K_6$-minor free
(e.g. \cite{SeymourLinkless}, \cite{Lovasz-Schrijver}), combining
with Theorem \ref{mainThm} we conclude:
\begin{cor}\label{Link&Stress}
Linklessly embeddable graphs are generically $4$-stress free.
\end{cor}
Let $\mu(G)$ denote the Colin de Verdi\`{e}re's parameter of a graph
$G$ \cite{Colin1}. Colin de Verdi\`{e}re \cite{Colin1} proved that a
graph $G$ is planar iff $\mu(G)\leq 3$; Lov\'{a}sz and Schrijver
\cite{Lovasz-Schrijver} proved that $G$ is linklessly embeddable iff
$\mu(G)\leq 4$. While we have seen that Theorem \ref{mainThm} fails
for $r\geq 8$, we conjecture that Theorem \ref{Gluck} and Corollary
\ref{Link&Stress} extend to:
\begin{conj}\label{Colin-conj}
Let $G$ be a graph and let $k$ be a positive integer. If
$\mu(G)\leq k$ then $G$ is generically $k$-stress free.
\end{conj}
For $k=1,2,3,4$ this is true: Colin de Verdi\`{e}re
\cite{Colin1} showed that the family $\{G: \mu(G)\leq k\}$ is
closed under taking minors for every $k$. Note that
$\mu(K_{r})=r-1$. By Theorem \ref{mainThm}, Conjecture \ref{Colin-conj} holds
for $k\leq 4$.
Conjecture \ref{Colin-conj} implies $$\mu(G)\leq k
\Rightarrow e\leq kv-(^{k+1}_{\ \ 2})$$ (where $e$ and $v$ are the
numbers of edges and vertices in $G$, respectively) which is not
known either.

Now we give a proof of Theorem \ref{mainThm} which relies on results
about graph minors which are developed in the next subsection,
\ref{SubSecMinors}.
\\
\\
\emph{Proof of Theorem \ref{mainThm}}:
For $r=2$ the
assertion of the theorem is trivial. Suppose $K_r\nprec G$, and
contract edges belonging to at most $r-3$ triangles as long as it
is possible. Denote the resulted graph by $G'$. Repeated
application of Lemmata \ref{Whiteley} and \ref{extWHiteley} asserts that if $G'$ is
generically $(r-2)$-stress free / $(r-2)$-hyperconnected, then so is $G$. In case $G'$ has
no edges, it is trivially $(r-2)$-stress free / hyperconnected. Otherwise, $G'$ has
an edge, and each edge belongs to at least $r-2$ triangles. For
$2<r<6$, by Proposition \ref{anyEdge5} $G'$ has a $K_r$ minor,
hence so has $G$, a contradiction. For $r=6$, by Proposition
\ref{anyEdge6} $G'$ either has a $K_6$ minor which leads to a
contradiction, or $G'$ is a clique sum over $K_r$ for some $r\leq
4$. In the later case, denote $G'=G_{1}\cup G_{2}$, $G_{1}\cap
G_{2}=K_r$. As the graph of a simplex is $k$-rigid/hyperconnected for any $k$, by
Lemmata \ref{weak-clique-sum,Asimov} and \ref{weak-clique-sum,Kalai} it is enough to show that
each $G_{i}$ is generically $(r-2)$-stress free / $(r-2)$-acyclic, which follows
from induction hypothesis on the number of vertices. $\square$
\\
\textbf{Remark}: We can prove the case
$r=5$ avoiding Mader's theorem, by using Wagner's structure theorem for $K_5$-minor free graphs
(\cite{Diestel}, Theorem 8.3.4) and Lemmata \ref{weak-clique-sum,Asimov} and \ref{weak-clique-sum,Kalai}. Using
Wagner's structure theorem for $K_{3,3}$-minor free graphs
(\cite{Diestel}, ex.20 on p.185) and Lemmata
\ref{weak-clique-sum,Asimov} and \ref{weak-clique-sum,Kalai}, we conclude that $K_{3,3}$-minor
free graphs are generically $4$-stress free / $4$-acyclic.

\subsection{Minors}\label{SubSecMinors}
All graphs we consider are simple, i.e. with no loops and no
multiple edges. Let $e=\{v,u\}$ be an edge in a graph $G$. By
\emph{contracting} $e$ we mean identifying the vertices $v$ and $u$ and
deleting the loop and one copy of each double edge created by this
identification, to obtain a new (simple) graph. A graph $H$ is
called a \emph{minor} of a graph $G$, denoted $H\prec G$, if by
repeated contraction of edges we can obtain $H$ from a subgraph of
$G$. In the sequel we shall make an essential use of the following
Theorem of Mader \cite{Mader}:
\begin{thm}\label{Minor&Edges}(Mader)
For $3 \leq r \leq 7$, if a graph $G$ on $n$ vertices has no $K_r$
minor then it has at most $(r-2)n-$ $r-1 \choose 2$ edges.
\end{thm}

\begin{prop}\label{anyEdge5}
For $3\leq r \leq 5$:
 If $G$ has an edge and each edge belongs to at least $r-2$ triangles,
 then $G$ has a $K_r$ minor.
\end{prop}
$Proof$: For $r=3$ $G$ actually contains $K_3$ as a subgraph. Let
$G$ have $n$ vertices and $e$ edges. Assume (by contradiction)
that $K_r\nprec G$. W.l.o.g. $G$ is connected.

For $r=4$, by Theorem \ref{Minor&Edges} $e\leq 2n-3$ hence there
is a vertex $u \in G$ with degree $d(u)\leq 3$. Denote by $N(u)$
the induced subgraph on the neighbors of $u$. For every $v\in
N(u)$, the edge $uv$ belongs to at least two triangles, hence
$N(u)$ is a triangle, and together with $u$ we obtain a $K_4$ as a
subgraph of $G$, a contradiction.

For $r=5$, by Theorem \ref{Minor&Edges} $e\leq 3n-6$ hence there
is a vertex $u \in G$ with degree $d(u)\leq 5$. Also $d(u)\geq 4$
(as $u$ is not an isolated vertex). If $d(u)=4$
then the induced subgraph on $\{u\}\cup N(u)$ is $K_{5}$, a
contradiction. Otherwise, $d(u)=5$. Every $v\in N(u)$ has degree
at least $3$ in $N(u)$, hence $e(N(u))\geq \lceil 3\cdot
5/2\rceil=8$. But $K_4\nprec N(u)$, hence $e(N(u))\leq 2\cdot
5-3=7$, a contradiction.$\square$

\begin{prop}\label{anyEdge6}
If $G$ has an edge and each edge belongs to at least $4$
triangles, then either $G$ has a $K_6$ minor, or $G$ is a clique
sum over $K_r$ for some $r\leq 4$ (i.e. $G=G_{1}\cup G_{2},
G_{1}\cap G_{2}=K_r$, $G_{i}\neq K_r$, $i=1,2$).
\end{prop}
$Proof$: We proceed as in the proof of Proposition \ref{anyEdge5}:
Assume that $K_6\nprec G$. W.l.o.g. $G$ is connected. By Theorem
\ref{Minor&Edges} $e\leq 4n-10$ hence there is a vertex $u \in G$
with degree $d(u)\leq 7$, also $d(u)\geq 5$. If $d(u)=5$ then
$N(u)=K_{5}$, a contradiction. Actually, $N(u)$ is planar: since
$N(u)$ has at most $7$ vertices, each of degree at least $4$, if
$N(u)$ were not $4$-connected, it must have exactly $7$ vertices and
two disjoint edges such that each of their $4$ vertices is adjacent
to the remaining $3$ vertices of $N(u)$ (whose removal disconnect
$N(u)$); but such graph has a $K_5$ minor. As $K_5\nprec N(u)$,
$N(u)$ is $4$-connected. Now Wagner's structure theorem for
$K_5$-minor free graphs (\cite{Diestel}, Theorem 8.3.4) asserts that
$N(u)$ is planar.

If $d(u)=6$, then $12=3\cdot 6-6 \geq e(N(u))\geq 4\cdot 6/2=12$
hence $N(u)$ is a triangulation of the $2$-sphere $S^{2}$. If
$d(u)=7$, then $15=3\cdot 7-6 \geq e(N(u))\geq 4\cdot 7/2=14$. We
will show now that $N(u)$ cannot have $14$ edges, hence it is a
triangulation of $S^{2}$: Assume that $N(u)$ has $14$ edges, so each
of its vertices has degree $4$, and $N(u)$ is a triangulation of
$S^{2}$ minus an edge. Let us look to the unique square (in a planar
embedding) and denote its vertices by $A$. The number of edges
between $A$ and $N(u)\setminus A$ is $8$. Together with the $4$
edges in the subgraph induced by $A$, leaves two edges for the
subgraph induced by $N(u)\setminus A=\{a,b,c\}$; let $a$ be their
common vertex. We now look at the neighborhood of $a$ in a planar
embedding (it is a $4$-cycle): $b,c$ must be opposite in this square
as $\{b,c\}$ is missing. Hence for $v\in A\cap N(a)$ we get that $v$
has degree $5$, a contradiction.

Now we are left to deal with the case where $N(u)$ is a
triangulation of $S^{2}$, and hence a maximal $K_5$-minor free
graph. If $G$ is the cone over $N(u)$ with apex $u$, then every
edge in $N(u)$ belongs to at least $3$ triangles in $N(u)$. By
Proposition \ref{anyEdge5}, $N(u)$ has a $K_5$ minor, a
contradiction. Hence there exists a vertex $w\neq u$, $w\in
G\setminus N(u)$.
Denote by $[w]$ the set of all vertices in $G$ connected to $w$ by
a path disjoint from $N(u)$. Denote by $N'(w)$ the induced graph
on the vertices in $N(u)$ that are neighbors of some vertex in
$[w]$. If $N'(w)$ is not a clique, there are two non-neighbors
$x,y\in N'(w)$, and a path through vertices of $[w]$ connecting
them. This path together with the cone over $N(u)$ with apex $u$
form a subgraph of $G$ with a $K_6$ minor, a contradiction.

Suppose $N'(w)$ is a clique (it has at most $4$ vertices, as
$N(u)$ is planar). Then $G$ is a clique sum of two graphs that
strictly contain $N'(w)$: Let $G_{1}$ be the induced graph on
$[w]\cup N'(w)$ and let $G_{2}$ be the induced graph on
$G\setminus [w]$. Then $G=G_{1}\cup G_{2}$ and $G_{1}\cap
G_{2}=N'(w)$. $\square$
\\
\textbf{Remark}:
In view of Theorem \ref{Minor&Edges} for the case $r=7$, we may
expect the following to be true:
\begin{prob}\label{r=7,8}
If $G$ has an edge and each edge belongs to at least $5$
triangles, then either $G$ has a $K_{7}$ minor, or $G$ is a clique
sum over $K_{l}$ for some $l\leq 6$.
\end{prob}
If true, it extends the assertion of Theorem \ref{mainThm} to the
case $r=7$. We could show only the weaker assertion
$$G\ \rm{has\ a\ generic\ 5-stress} \Rightarrow K_7^-\prec G,$$ where $K_7^-$
is $K_7$ minus an edge, by using similar
arguments to those used in this section.

\subsection{Shifting and embedding into $2$-manifolds}
Theorem \ref{Gluck} may be extended to other $2$-manifolds as
follows:
\begin{thm}\label{surface-thm}
Let $M\neq S^2$ be a compact connected $2$-manifold without
boundary, and let $G$ be a graph. Suppose that $\{r-1,r\}\in
\Delta(G)$ and $K_{r}$ can not be embedded in $M$. Then $G$ can
not be embedded in $M$.
\end{thm}
$Proof$: Let $g$ be the genus of $M$, then $g>0$ (e.g. the torus has
genus 1, the projective plane has genus 1/2). Assume by
contradiction that $G$ embeds in $M$. By looking at the rigidity
matrix we note that deleting from $G$ a vertex of degree at most
$r-2$ preserves the existence of $\{r-1,r\}$ in the shifted graph.
Deletion preserves embeddability in $M$ as well. Thus we may assume
that $G$ has minimal degree $\delta(G)\geq r-1$. By Euler formula
$e\leq 3v-6+6g$ (where $e$ and $v$ are the numbers of edges and
vertices in $G$ respectively). Also $e\geq (r-1)v/2$, hence $v\leq
\frac{12g-12}{(r-1)-6}$. Thus $(r-1)^2-5(r-1)+(6-12g)\leq 0$ which
implies $r \leq (7+\sqrt{1+48g})/2$. As $K_{r}$ can not be embedded
in $M$, by Ringel and Youngs \cite{Ringel} proof of Heawood's
map-coloring conjecture $r > (7+\sqrt{1+48g})/2$, a
contradiction.$\square$
\\
\textbf{Remark}: For any compact connected $2$-manifold without
boundary of positive genus, $M$, if $M$ is embedded in $\mathbb{R}^3$ two
linked simple closed curves on it exist. One may ask whether the
graph of any triangulated such $M$ is always not linkless.

For the
projective plane this is true. It follows from the fact that the
two minimal triangulations of the projective plane (w.r.t. edge
contraction), determined by Barnette \cite{BarnetteRP2triang}, have a minor
from the Petersen family, and hence are not linkless, by the
result of Robertson, Seymour and Thomas \cite{SeymourLinkless}. Moreover,
the graph of any polyhedral map of the projective plane is not
linkless, as its $7$ minimal polyhedral maps (w.r.t. edge
contraction), determined by
Barnette \cite{BarnetteRP2polyhed}, have graphs equal to $6$ of the members
in Petersen family.

Examining the $21$ minimal triangulations of the torus, see
Lavrenchenko \cite{Lavrenchenko}, we note that $20$ of them have a
$K_6$ minor, and hence are not linkless, but the last one is
linkless, see Figure 4.1
(one checks that it contains no minor from Petersen's family).
Taking connected sums of this
triangulation, we obtain linkless graphs triangulating any
oriented surface of positive genus. By performing stellar
operations we obtain linkless graphs with arbitrarily many
vertices triangulating any oriented surface of positive genus.

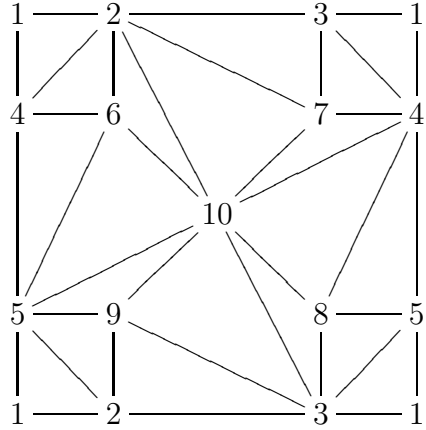
\begin{figure}\label{fig:T21}
\newcommand{\edge}[1]{\ar@{-}[#1]}
\newcommand{\lulab}[1]{\ar@{}[l]^<<{#1}}

\newcommand{\rulab}[1]{\ar@{}[r]^<<{#1}}
\newcommand{\rrulab}[1]{\ar@{}[rr]^<<{#1}}

\newcommand{\ldlab}[1]{\ar@{}[l]^<<{#1}}
\newcommand{\rdlab}[1]{\ar@{}[r]_<<{#1}}
\newcommand{\rrdlab}[1]{\ar@{}[rr]_<<{#1}}
\newcommand{\node}{*+[O][F-]{ }}
\centerline{ \xymatrix{
1  \edge{d} \edge{r} & 2 \edge{d} \edge{rr} \edge{ddr} \edge{rrd} && 3 \edge{d} \edge{r} \edge{dr} & 1 \edge{d}\\
4  \edge{dd} \edge{ur} \edge{r} & 6 \edge{dr} \edge{ddl} && 7 \edge{r} \edge{dl} & 4 \edge{dll} \edge{dd} \edge{ddl} \\
&& 10 \edge{dll} \edge{dl} \edge{ddr} \edge{dr} && \\
5 \edge{d} \edge{r} \edge{dr} & 9 \edge{drr} \edge{d} && 8 \edge{r} \edge{d} & 5 \edge{d} \edge{dl} \\
1 \edge{r} & 2 \edge{rr} && 3 \edge{r} & 1
} } \caption {Linkless graph of a torus}
\end{figure}


\section{Open problems}
\begin{enumerate}
\item Can the $S_i$'s in Theorem \ref{structureLemma} be taken to be homology spheres?

\item Can the intersections in Theorem \ref{structureLemma} be guaranteed to be CM?

In view of Proposition \ref{prop:Yonatan}, if the intersections
$S_i\cap(\cup_{j<i}S_j)$ in Theorem \ref{structureLemma} can be
taken to be CM, and the $S_i$'s can be taken to be homology
spheres, then Conjecture \ref{g-conj 2CM} would be reduced to the
conjecture that homology spheres have the weak-Lefschetz property; see Conjecture \ref{conj-g-hierarchy}(2).

\item Must a graph with a generic $(r-2)$-stress contain a subdivision of $K_r$ for $2\leq r\leq 6$?

The answer is positive for $r=2,3,4$ as in this case $G$ has a $K_r$
minor iff $G$ contains a subdivision of $K_r$ (\cite{Diestel},
Proposition 1.7.2). Mader proved that every graph on $n$ vertices
with more than $3n-6$ edges contains a subdivision of $K_{5}$
\cite{Mader2}. A positive answer in the case $r=5$ would strengthen
this result.

\item Let $G$ be a graph and let $k$ be a positive integer. Show that
$\mu(G)\leq k$ implies that $G$ is generically $k$-stress free.

\item Assume that $G$ has an edge and each edge belongs to at least $5$
triangles. Show that either $G$ has a $K_{7}$ minor, or $G$ is a clique
sum over $K_{l}$ for some $l\leq 6$.

If true, it implies that $\{6,7\}\in\Delta(G)$ forces a $K_7$ minor in $G$.

\item Is the graph of a triangulated non orientable $2$-manifold always not linkless?

\item Prove Charney-Davis conjecture \cite{Charney-Davis} for clique $3$-spheres using rigidity (shifting) arguments in order to give a simpler proof than in \cite{Davis-Okun}. We repeat their conjecture: Let $K$ be a $(d-1)$ dimensional clique (homology) sphere (that
is, all its missing faces are $1$ dimensional), where $d$ is even.
Show that $(-1)^{\frac{d}{2}}\sum_{i=0}^{d}h_{i}(K)
\geq 0$. Equivalently, $\sum_{0\leq k\leq \frac{d}{2}}
(-1)^{\frac{d}{2}-k}g_{k}(K) \geq 0$. In case $d=4$, the conjecture reads $f_{1}(K)\geq 5f_{0}(K)-16$ (to be
compared with the LBT for spheres: $f_{1}(K)\geq 4f_{0}(K)-10$).

\end{enumerate}


\chapter{Lefschetz Properties and Basic Constructions on Simplicial Spheres}\label{chapter:HL_WL}
\section{Basics of Lefschetz properties}\label{sec:Lefschetz}
Our motivating problem is the following well known McMullen's
$g$-conjecture for spheres. Recall that by \emph{homology sphere}
(or \emph{Gorenstein$^*$ complex}) we mean a pure simplicial complex
$L$ such that for every face $F\in L$ (including the empty set),
$\rm{lk}(F,L)$ has the same homology (say with integer coefficients)
as of a $\rm{dim}(\rm{lk}(F,L))$-sphere.
\begin{conj}(McMullen \cite{McMullen-g-conj})\label{conj-g}
Let $L$ be a homology sphere, then its $g$-vector is an $M$-sequence.
\end{conj}

An algebraic approach to this problem is to associate with $L$ a
standard ring whose Hilbert function is $g(L)$, the $g$-vector of
$L$. This was worked out successfully by Stanley \cite{St} in his
celebrated proof of Conjecture \ref{conj-g} for the case where $L$
is the boundary complex of a simplicial polytope. The hard-Lefschetz
theorem for toric varieties associated with rational polytopes,
translates in this case to the following property of face rings,
called \emph{hard-Lefschetz}.

Let $K$ be a simplicial complex on the vertex set $[n]$. Let
$A=\mathbb{R}[x_1,..,x_n]$ be the polynomial ring, each variable has
degree one. Recall that the face ring of $K$ is
$\mathbb{R}[K]=A/I_{K}$ where $I_K$ is the ideal in $A$ generated by
the monomials whose support is not an element of $K$. Let
$\Theta=(\theta_1,..,\theta_d)$ be an l.s.o.p. of $\mathbb{R}[K]$ -
it exists, e.g. \cite{StanleyGreenBook}, Lemma 5.2, and generic
1-forms $y_1,...,y_d$ from the basis $Y$ of $A_1$ (recall from
subsection \ref{subsec:dual symm}) will do. Denote $H(K)=
H(K,\Theta)=\mathbb{R}[K]/(\Theta)=H(K)_0\oplus H(K)_1\oplus...$
where the grading is induced by the degree grading in $A$, and
$(\Theta)$ is the ideal in $\mathbb{R}[K]$ generated by the images
of the elements of $\Theta$ under the projection $A\rightarrow
\mathbb{R}[K]$. $K$ is called \emph{Cohen-Macaulay} (CM for short)
over $\mathbb{R}$ if for an (equivalently, every) l.s.o.p. $\Theta$,
$\mathbb{R}[K]$ is a free $\mathbb{R}[\Theta]$-module. If $K$ is CM
then $\ddim_{\mathbb{R}}H(K)_i=h_i(K)$. (The converse is also true:
$h$ is an $M$-vector iff $h=h(K)$ for some CM complex $K$
\cite{StanleyGreenBook}, Theorem 3.3.) For $K$ a CM simplicial
complex with symmetric $h$-vector, if there exists an l.s.o.p.
$\Theta$ and an element $\omega\in A_1$ such that the multiplication
maps $\omega^{d-2i}: H(K,\Theta)_i\longrightarrow
H(K,\Theta)_{d-i}$, $m\mapsto \omega^{d-2i}m$, are isomorphisms for
every $0\leq i\leq \lfloor d/2\rfloor$, we say that $K$ has the
\emph{hard-Lefschetz property}, or that $K$ is HL.

As was shown by Stanley \cite{St}, for $K$
 the boundary complex of a simplicial $d$-polytope $P$,
 the l.s.o.p $\Theta$ induced by the embedding of $P_0$ in $\mathbb{R}^d$
 and $\omega=\sum_{1\leq i\leq n}x_i$ demonstrate that $K$ is HL;
 hence so do generic $y_1,...,y_{d+1}\in Y$. In terms of $GIN$ this is
equivalent to requiring that non of the monomials
$y_{d+1}^{d-2k-1}y_{d+2}^{k+1}$ are in $GIN(K)$, where $k=0,1,...$.
Indeed, these monomials are not in $GIN(K)$ iff the maps
$y_{d+1}^{d-2i}: H(K)_i\longrightarrow H(K)_{d-i}$ are onto, and
when $h(K)$ is symmetric this happens iff these maps are
isomorphisms.

Let us translate the hard-Lefschetz property from terms of $GIN$
into terms of symmetric shifting, as in \cite{skira}. Let
$\Delta(d,n)$ be the pure $(d-1)$-dimensional simplicial complex
with set of vertices $[n]$ and facets $\{S: S\subseteq [n], |S|=d,\
k\notin S\Rightarrow [k+1,d-k+2]\subseteq S\}$. Equivalently,
$\Delta(d,n)$ is the maximal pure $(d-1)$-dimensional simplicial
complex with vertex set $[n]$ which does not contain any of the sets
$T_d,...,T_{\lceil d/2\rceil}$, where
\begin{equation}\label{eq:T}
T_{d-k}=\{k+2,k+3,...,d-k,d-k+2,d-k+3,...,d+2\},\  0\leq k\leq \lfloor d/2\rfloor.
\end{equation}
Note that $\Delta(d,n)\subseteq\Delta(d,n+1)$, and define $\Delta(d)=\cup_n\Delta(d,n)$.
Kalai refers to the relation
\begin{equation}\label{eq:shUBT}
\Delta(K)\subseteq \Delta(d)
\end{equation}
as the \emph{shifting theoretic upper bound theorem}. Using the map
from $GIN(K)$ to $\Delta^s(K)$, we have just seen that for CM
$(d-1)$-dimensional complexes with symmetric $h$-vector,
$\Delta^s(K)\subseteq \Delta(d)$ is equivalent to $K$ being HL.

To justify the terminology in (\ref{eq:shUBT}), note that the
boundary complex of the cyclic $d$-polytope on $n$ vertices, denoted
by $C(d,n)$, satisfies $\Delta^s(C(d,n))=\Delta(d,n)$. This follows
from the fact that $C(d,n)$ is HL. Recently Murai \cite{MuraiCyclic}
proved that also $\Delta^e(C(d,n))=\Delta(d,n)$, as was conjectured
by Kalai \cite{skira}. It follows that if $K$ has $n$ vertices and
(\ref{eq:shUBT}) holds, then the $f$-vectors satisfy $f(K)\leq
f(C(d,n))$ componentwise.

For $K$ as above, weaker than the hard-Lefschetz property is to
require only that multiplications $y_{d+1}:
H(K)_{i-1}\longrightarrow H(K)_{i}$ are injective for $1\leq i\leq
\lceil d/2\rceil$ and surjective for $\lceil d/2\rceil<i\leq d$,
called here \emph{unimodal weak-Lefschetz} property (sometimes it is
called weak-Lefschetz in the literature). Even weaker is just to
require that multiplications $y_{d+1}: H(K)_{i-1}\longrightarrow
H(K)_{i}$ are injective for $1\leq i\leq \lfloor d/2\rfloor$, which
we refer to as the \emph{weak-Lefschetz property}, and say that $K$
is WL. (Injectivity for $i\leq \lceil d/2\rceil$ in the case of
Gorenstein$^*$ complexes implies also surjective maps for $\lceil
d/2\rceil<i\leq d$; see the proof of Theorem \ref{thmSwartz} below.)
This is equivalent to the following, in the case of symmetric
shifting \cite{BK-homology}:
\begin{eqnarray}\label{eq:shWL}
(1)\ S\in \Delta(K), |S|=k \ \Rightarrow \ [d-k]\cup S\in \Delta(K), \nonumber
\\
(2)\ S\in \Delta(K), |S|=k<\lfloor d/2\rfloor \ \Rightarrow \ \{d-k+1\}\cup S\in \Delta(K).
\end{eqnarray}
Condition (1) holds when $K$ is CM, and condition (2) holds iff $K$
is WL. As was noticed in \cite{BK-homology}, (\ref{eq:shWL}) is
implied by requiring that $\Delta(K)$ is pure and every
$S\in\Delta(K)$ of size less than $\lfloor d/2\rfloor$ is contained
in at least $2$ facets of $\Delta(K)$.

Note that if $L$ is a homology sphere, it is in particular CM with a
symmetric $h$-vector. If in addition it has the weak-Lefschetz
property, then in the standard ring
$S(L)=\mathbb{R}[K]/(\Theta,y_{d+1},A_{1+\lfloor
d/2\rfloor})=H(L,\Theta)/(y_{d+1},H_{1+\lfloor
d/2\rfloor})=S_0\oplus S_1\oplus...$ the following holds:
$g_i(L)=\ddim_\mathbb{R}S_i$ for all $0\leq i\leq \lfloor
d/2\rfloor$, and Conjecture \ref{conj-g} holds for $L$.

We summarize the discussion above in the following hierarchy of conjectures, where assertion $(i)$ implies assertion $(i+1)$:
\begin{conj}\label{conj-g-hierarchy}
Let $L$ be a homology $(d-1)$-sphere. Then:

(1) If $S\in\Delta(L)$, $|S|=k\leq \lfloor d/2\rfloor$ and $S\cap[d-k+1]=\emptyset$ then $S\cup[k+2,d-k+1]\in\Delta(L)$.

This is equivalent to $\Delta(K)\subseteq \Delta(d)$, and in the
symmetric case this is equivalent to $L$ being HL.

(2) If $S\in\Delta(L)$, $|S|=k< \lfloor d/2\rfloor$ and
$S\cap[d-k+1]=\emptyset$ then $S\cup[\lceil
d/2\rceil+2,d-k+1]\in\Delta(L)$. In the symmetric case this is
equivalent to $L$ being WL.

(3) $g(L)$ is an $M$-vector.
\end{conj}

\section{Hard Lefschetz and join}\label{sec:HL&Join}
Let $S$ be a Cohen-Macaulay $(d-1)$-simplicial complex over a field $k$.
 If there exists a degree one
element $\omega$ such that multiplication
\begin{equation}\label{eq:i-Lefschetz}
\omega^{d-2i}:H(S)_i\rightarrow H(S)_{d-i}
\end{equation}
is an isomorphism (for some l.s.o.p.) we say that $S$ is
\emph{$i$-Lefschetz} and that $\omega$ is an \emph{$i$-Lefschetz
element} of $H(S)$. If (\ref{eq:i-Lefschetz}) holds for every $0\leq
i\leq d/2$ then $S$ is \emph{HL} and $\omega$ is an
\emph{HL-element} of $H(S)$.

Let us recall a few ring theoretic terms, see e.g.
\cite{StanleyGreenBook} for details. Let $A=k[x_1,...,x_n]$. An
$A$-module $M$ of dimension $d$ is \emph{Cohen-Macaulay} (CM) if for
generic (i.e., for some algebraically independent) $y_1,...,y_d\in
A_1$ $M$ is a free module over the subring $k[y_1,..,y_d]$.   If in
addition $\ddim_k\soc M/(y_1,...,y_d)M=1$ where $\soc M:=\{u\in M:
A_+ u=0\}$, $A_+=A_1\oplus A_2\oplus...$ and $(y_1,...,y_d)$ is the
obvious ideal in $A$, then $M$ is \emph{Gorenstein}. Note that
Gorenstein$^*$ complexes have Gorenstein face rings (as
$A$-modules).
\begin{lem}\label{lemma:Poincare pairing}
Let $M$ be a $d$-dimensional Gorenstein module over the polynomial
ring $R=\mathbb{R}[x_1,...,x_n]$, with an l.s.o.p. $\Theta$. Denote
$H=M/(\Theta)M$. Then for every
$0\leq i\leq d/2$ the pairing $H_i\times H_{d-i}\rightarrow
\mathbb{R}$, $(x,y)\mapsto \alpha(xy)$ is non-degenerated under any
fixed isomorphism $\alpha: H_d\cong \mathbb{R}$.
\end{lem}
$Proof$: $M$ is Gorenstein, so, by definition,
$\rm{dim}_{\mathbb{R}}\rm{soc}H=1$. As $M$ is $d$-dimensional,
$\rm{dim}_{\mathbb{R}}H_d\geq 1$, but $H_d\subseteq \rm{soc}H$, thus
$\rm{soc}H=H_d$. As $R_+$ is generated by $\{x_1,...,x_n\}$, we get
that for every $0\leq i<d$ and $0\neq u\in H_i$ there exists $x_j$
such that $x_ju\neq 0$, and inductively there exists a monomial $m$
of degree $d-i$ such that $mu\neq 0$, thus the pairing is
non-degenerated. $\square$

\begin{lem}\label{lemma:w-inner product}
Let $K$ be a $(d-1)$-dimensional Gorenstein$^*$ complex with an
l.s.o.p. $\Theta$ and an HL element $\omega$ over the reals. Let
$H=\mathbb{R}[K]/(\Theta)$ and fix an isomorphism $\alpha: H_d\cong
\mathbb{R}$. Then for every $0\leq i\leq d/2$ there is an induced
non degenerated bilinear form on $H_i$ given by $<x,y>=\alpha(\omega^{d-2i}xy)$.
\end{lem}
$Proof$: Clearly $<,>$ is bilinear and symmetric. For $0\neq x\in
H_i$, by assumption $0\neq \omega^{d-2i}x\in H_{d-i}$, and by Lemma
\ref{lemma:Poincare pairing} there exists $y\in H_i$ such that
$\alpha(\omega^{d-2i}xy)\neq 0$, hence $<,>$ is not degenerated.
$\square$

\begin{lem}\label{lemma:w-invariant subspaces}
Under the assumptions of Lemma \ref{lemma:w-inner product}, $H$
decomposes into a direct sum of $\mathbb{R}[\omega]$-invariant
spaces, each is of the form $$V_m=\mathbb{R}m\oplus \mathbb{R}\omega
m\oplus...\oplus \mathbb{R}\omega^{d-2i} m$$ for $m\in
\mathbb{R}[K]/(\Theta)$ of degree $i$ for some $0\leq i \leq d/2$.
\end{lem}
$Proof$: $V_1$ ($1\in H_0$) is an $\mathbb{R}[\omega]$-invariant
space which contain $H_0$. Assume that for $1\leq i\leq d/2$ we have
already constructed a direct sum of $\mathbb{R}[\omega]$-invariant
spaces, $\tilde{V}_{i-1}$, which contains
$\tilde{H}_{i-1}:=H_0\oplus...\oplus H_{i-1}$, in which each $V_m$
contains some nonzero element of $\tilde{H}_{i-1}$. We now extend
the construction to have these properties w.r.t. $\tilde{H}_i$. By
assumption $\omega H_{i-1} \subseteq \tilde{V}_{i-1}\cap H_i$. Let
$m_1,...,m_t\in H_i$ form a basis to the subspace of $H_i$
orthogonal to $\omega H_{i-1}$ w.r.t. the inner product from Lemma
\ref{lemma:w-inner product}. Let
\begin{equation}\label{eq:V_i-sum}
\tilde{V}_{i}=\tilde{V}_{i-1}+ V_{m_1}+...+ V_{m_t}.
\end{equation}
We first show that each $V_{m_j}$ is $\mathbb{R}[\omega]$-invariant,
i.e. that $\omega^{d-2i+1}m_j=0$ for $1\leq j\leq t$. By Lemma
\ref{lemma:Poincare pairing} it is enough to show that for every
$x\in H_{i-1}$ $\alpha(x\omega^{d-2i+1}m_j)=0$. As
$\alpha(x\omega^{d-2i+1}m_j)=<\omega x,m_j>$ indeed it equals zero.

Next we show that the sum in (\ref{eq:V_i-sum}) is direct. As
$m_1,...,m_t\in H_i$ are linearly independent and
$\omega^{d-2i}:H_i\rightarrow H_{d-i}$ is injective, then the sum
$W_i:=V_{m_1}+...+ V_{m_t}$ is direct. To show that
$\tilde{V}_{i-1}\oplus W_i$, we check that for every $i\leq l \leq
d-i$ $H_l\cap \tilde{V}_{i-1}\cap W_i=0$ (for $l<i$ and for $l>d-i$
$W_i\cap H_l=0$). Indeed, an element in the intersection is of the
form $\omega^{l-i}\omega x=\omega^{l-i} y$ where $x\in H_{i-1}$ and
$y\in H_i$ is orthogonal to $\omega x$. Injectivity of
$\omega^{l-i}:H_i\rightarrow H_{l}$ implies $y=\omega x$, which
equals zero by orthogonality.

As the $h$-vector of $K$ is symmetric, $\tilde{V}_{\lfloor d/2
\rfloor}=H$, giving the desired decomposition. $\square$
\\
\textbf{Remark}: Even if $K$ is not HL we still get a decomposition
into a direct sum of irreducible $\mathbb{R}[\omega]$-invariant
spaces $V_m=\mathbb{R}m\oplus \mathbb{R}\omega m\oplus...\oplus
\mathbb{R}\omega^l m$, but no longer $l=d-2\ddeg(m)$.
\begin{thm}\label{thm:*} (With Eric Babson)
Let $K$ and $L$ be Gorenstein$^*$ complexes on disjoint sets of
vertices, of dimensions $d_K-1,d_L-1$, with l.s.o.p's
$\Theta_K,\Theta_L$ and HL elements $\omega_K,\omega_L$
respectively; over the reals. Then:

(0) $K*L$ has a symmetric $h$-vector and dimension $d_K+d_L-1$.

(1) $\Theta_K\biguplus\Theta_L$ is an l.s.o.p for $K*L$ (over
$\mathbb{R}$).

(2) $\omega_K+\omega_L$ is an HL element of
$\mathbb{R}[K*L]/(\Theta_K\biguplus\Theta_L)$.
\end{thm}
$Proof:$ The $h$-polynomials satisfy $h(t,K*L)=h(t,K)h(t,L)$, hence
the symmetry of $h(t,K*L)$ follows from that of $h(t,K)$ and
$h(t,L)$:
$$t^{d_K+d_L}h(\frac{1}{t},K*L)=t^{d_K}h(\frac{1}{t},K)
t^{d_L}h(\frac{1}{t},L)=h(t,K)h(t,L)=h(t,K*L).$$

For a set $I$ let $A_I:=\mathbb{R}[x_i: i\in I]$ be a polynomial
ring. The isomorphism $A_{K_0}\bigotimes_{\mathbb{R}}A_{L_0}\cong
A_{K_0\biguplus L_0}$, $a_K\otimes a_L\mapsto a_Ka_L$ induces a
structure of an $A=A_{K_0\biguplus L_0}$ module on
$\mathbb{R}[K]\bigotimes_{\mathbb{R}}\mathbb{R}[L]$, isomorphic to
$\mathbb{R}[K*L]$, by $m_K\otimes m_L\mapsto m_Km_L$ and
$(a_K\otimes a_L)(m_K\otimes m_L)=a_Km_K\otimes a_Lm_L$. (E.g.
$a_K\in A_{K_0}\subseteq A$ acts like $a_K\otimes 1$ on
$\mathbb{R}[K]\bigotimes_{\mathbb{R}}\mathbb{R}[L]$. )

The above isomorphism induces an isomorphism of $A$-modules
\begin{equation}\label{eq:JoinIsom}
\mathbb{R}[K*L]/(\Theta_K\biguplus\Theta_L) \cong
\mathbb{R}[K]/(\Theta_K) \bigotimes_{\mathbb{R}}
\mathbb{R}[L]/(\Theta_L),
\end{equation}
proving (1). Actually, $\mathbb{R}[K*L]$
is both a finitely generated and free
$\mathbb{R}[\Theta_K\biguplus\Theta_L]$-module, by Cohen-Macaulayness.

By Lemma \ref{lemma:w-invariant subspaces},
$\mathbb{R}[K]/(\Theta_K)$ decomposes into a direct sum of
$\mathbb{R}[\omega_K]$-invariant spaces, each is of the form
$V_m=\mathbb{R}m\bigoplus \mathbb{R}\omega_K m\bigoplus...\bigoplus
\mathbb{R}\omega_K^{d_K-2i} m$ for $m\in \mathbb{R}[K]/(\Theta_K)$
of degree $i$ for some $0\leq i \leq d_K/2$; and similarly for
$\mathbb{R}[L]/(\Theta_L)$.

The $\mathbb{R}[\omega_K]$-module $V_m$ is isomorphic to the
$\mathbb{R}[\omega]$-module $\mathbb{R}[\partial
\sigma^{d_K-2i}]/(\theta)$ by $\omega_K\mapsto \omega$ and $m\mapsto
1$, where $\sigma^j$ is the $j$-simplex, $\theta$ is an l.s.o.p.
induced by the positions of the vertices in an embedding of
$\sigma^{d_K-2i}$ as a full dimensional geometric simplex in
$\mathbb{R}^{d_K-2i}$ with the origin in its interior, and
$\omega=\sum_{v\in \sigma_0}x_v$ is an HL element for
$\mathbb{R}[\partial \sigma^{d_K-2i}]/(\theta)$. Thus, to prove (2)
it is enough to prove it for the join of boundaries of two simplices
with l.s.o.p.'s as above and the HL elements having weight $1$ on
each vertex of the ground set.

Note that the join $\partial \sigma^{k}*\partial \sigma^{l}$ is
combinatorially isomorphic to the boundary of the polytope
$P:=\rm{conv}(\sigma^{k}\cup_{\{0\}} \sigma^{l})$ where $\sigma^{k}$
and $\sigma^{l}$ are embedded in orthogonal spaces and intersect
only in the origin which is in the relative interior of both.
McMullen's proof of the $g$-theorem for simplicial polytopes
\cite{McMullen-g-proof1,McMullen-g-proof2} states that $\sum_{v\in
P_0}x_v = \omega_{\partial \sigma^{k}}+\omega_{\partial \sigma^{l}}$
is indeed an HL element of $\mathbb{R}[\partial \sigma^{k}*\partial
\sigma^{l}]/(\Theta_{\partial P})$ where $\Theta_{\partial P}$ is
the l.s.o.p. induced by the positions of the vertices in the
polytope $P$. By the definition of $P$, $\Theta_{\partial
P}=\Theta_{\partial \sigma^{k}}\uplus \Theta_{\partial \sigma^{l}}$.
Thus (2) is proved. $\square$
\\
\textbf{Remark}: As a nonzero multiple of an HL element is again HL,
then in Theorem \ref{thm:*}(2) any element $a\omega_K+b\omega_L$
where $a,b\in \mathbb{R}$, $ab\neq 0$, will do.

\begin{cor}\label{cor:*}
Let $K$ and $L$ be HL simplicial/homology/piecewise linear spheres of
dimensions $k,l$ respectively. Then their join $K*L$ is an HL
$(k+l+1)$-simplicial/homology/piecewise linear sphere.
\end{cor}
$Proof:$ As simplicial/homology/piecewise linear spheres are
Gorenstein$^*$, the corollary follows at once from Theorem
\ref{thm:*} and the fact that join of simplicial/homology/piecewise linear spheres is
again a simplicial/homology/piecewise linear sphere of appropriate dimension. $\square$


\section{Weak Lefschetz and gluing}
The proof of the following proposition is similar to the proof that
pure shellable complexes are Cohen-Macaulay due to Stanley
\cite{Stanley-ShellableCM}; see also \cite{Bruns-Herzog}, Theorem
5.1.13.
\begin{prop}(with Yhonatan Iron)\label{prop:Yonatan}
Let $K$, $L$ and $K\cap L$ be simplicial complexes of the same
dimension $d-1$. Assume that $K$ and $L$ are weak-Lefschetz. If $K\cap L$ is CM then $K\cup L$
is weak-Lefschetz.
\end{prop}
$Proof:$
For any two complexes $K,L$ with $(K\cup L)_0=[n]$ the inclusions $K\cap L\subseteq K,L\subseteq K\cup L$ induce a short exact sequence of $A=\mathbb{R}[x_1,...,x_n]$ modules
\begin{equation} \label{eq:YonatanShortExact}
\begin{CD}
0@> >>\mathbb{R}[K\cup L]@> >>\mathbb{R}[K]\oplus \mathbb{R}[L]@> >>\mathbb{R}[K\cap L]@> >>0
\end{CD}
\end{equation}
(maps are given by projections). The above four complexes have the
same dimension, hence they have a common l.s.o.p. $\Theta$ (as
intersection of finitely many nonempty Zariski open sets is
nonempty). As the functor $\otimes_AA/\Theta$ is right exact, we
obtain the following commutative diagram of exact sequences for each
$1\leq i\leq \lfloor d/2\rfloor$:
\begin{equation} \label{eq:YonatanTorExact}
\begin{CD}
\rm{Tor}(\mathbb{R}[K\cap L],\frac{A}{\Theta})_{i-1}@>\delta_{i-1} >>\frac{\mathbb{R}[K\cup L]}{(\Theta)} _{i-1}@> >>\frac{\mathbb{R}[K]}{(\Theta)} _{i-1}\oplus \frac{\mathbb{R}[L]}{(\Theta)} _{i-1}@> >>\frac{\mathbb{R}[K\cap L]}{(\Theta)} _{i-1}@> >>0\\
@VV V @VV\omega V @VV(\omega,\omega)V @VV\omega V @VV V \\
\rm{Tor}(\mathbb{R}[K\cap L],\frac{A}{\Theta})_{i}@>\delta_i >>\frac{\mathbb{R}[K\cup L]}{(\Theta)}_{i}@> >>\frac{\mathbb{R}[K]}{(\Theta)}_{i}\oplus \frac{\mathbb{R}[L]}{(\Theta)}_{i}@> >>\frac{\mathbb{R}[K\cap L]}{(\Theta)}_{i}@> >>0\\
\end{CD}
\end{equation}
where the horizontal arrows preserve grading and the vertical arrows
are multiplication by a generic $\omega\in A_1$, i.e. $\omega$ is a
WL element for both $K$ and $L$. (For the middle terms we used
distributivity of $\otimes$ and $\oplus$.)

In order to show that $\omega:(\frac{\mathbb{R}[K\cup
L]}{(\Theta)})_{i-1}\rightarrow (\frac{\mathbb{R}[K\cup
L]}{(\Theta)})_{i}$ is injective, it is enough to show that
$\delta_{i-1}=0$, which of course holds if
$\rm{Tor}(\mathbb{R}[K\cap L],A/\Theta)_{i-1}=0$. Note that for an
$A$-module $M$
$\rm{Tor}(M,A/\Theta)=\rm{Ker}(M\otimes_A(\Theta)\rightarrow M)$. As
$M=\mathbb{R}[K\cap L]$ is CM, it is a free
$\mathbb{R}[\Theta]$-module, hence $\rm{Tor}(\mathbb{R}[K\cap
L],A/\Theta)_{i-1}=0$ for every $i$. $\square$
\\
\textbf{Remarks:} (1) Note that the above proof provides an even
more general condition on $K\cap L$ which already guarantees that
$K\cup L$ is WL.

(2) Compare Proposition \ref{prop:Yonatan} to \cite{Swartz}, remark
after proof of Theorem 3.9: there the gluing corresponds to an ear
decomposition of homology spheres and balls.


\section{Lefschetz properties and connected sum}\label{sec:WL,HL&Connected Sum}
Let $K$ and $L$ be pure simplicial complexes which intersect in a
common facet $<\sigma>=K\cap L$. Their \emph{connected sum over
$\sigma$} is $K\#_{\sigma}L=(K\cup L)\setminus \{\sigma\}$.

\begin{thm}\label{thm:connected-sum}
Let $K$ and $L$ be Gorenstein$^*$ complexes over $\mathbb{F}$ which
intersect in a common facet $<\sigma>=K\cap L$, of dimension $d-1$.
Let $A=\mathbb{F}[x_v:v\in (K\cup L)_0]$. Then:

(0) $K\#_{\sigma}L$ is Gorenstein$^*$ of dimension $d$; in
particular its $h$-vector is symmetric.

(1) Let $\Theta$ be a common l.s.o.p for $K$, $L$, $<\sigma>$ and
$K\#_{\sigma}L$ over $A$ (it exists) and assume that $\omega$ is an
HL element for both $K$ and $L$ w.r.t. $\Theta$. Then $\omega$ is an
$i$-Lefschetz element of $\mathbb{F}[K\#_{\sigma}L]/(\Theta)$ for
$0<i\leq d/2$.

(2) $K\#_{\sigma}L$ is HL.
\end{thm}

$Proof$: A straightforward Mayer-Vietoris and Euler characteristic
argument shows that $K\#_{\sigma}L$ is Gorenstein$^*$, and hence has
a symmetric $h$-vector. It is also easy to compute directly that
$h(K\#_{\sigma}L)=h(K)+h(L)-(1,0,0,...,0,1)$, a sum of symmetric
vectors, and hence is symmetric; also $h_0=h_d=1$.

For a simplicial complex $L$ let $\mathbb{F}(L):=\bigoplus_{a:
\rm{supp}(a)\in L}\mathbb{F}x^a$ be an $A_{L_0}=\mathbb{F}[x_v: v\in
L_0]$ module defined by $x_v(x^a)=\{^{x_vx^a \ \rm{if}\ v\cup
\rm{supp}(a)\in L}_{0 \ \rm{otherwise}}$.
Note that $\mathbb{F}(L)\cong \mathbb{F}[L]$ as $A_{L_0}$-modules.

Then the following is an exact sequence of $A$-modules:
\begin{equation}\label{eq:conn-sum-exact}
\begin{CD}
0 \rightarrow \mathbb{F}(<\sigma>) @>(\iota,-\iota)>>
(\mathbb{F}(K)\oplus \mathbb{F}(L)) @>\iota_K+\iota_L>>
\mathbb{F}(K\cup_{\sigma}L) \rightarrow 0
\end{CD}
\end{equation}
where the $\iota$'s denote the obvious inclusions. As a finite
intersection of Zariski nonempty open sets is nonempty, $\Theta$ as
in (1) exists (see Lemma \ref{lemOmegaZariski}). When we mod out
$\Theta$ from (\ref{eq:conn-sum-exact}), which is the same as tensor
(\ref{eq:conn-sum-exact}) with $\otimes_A A/\Theta$, we obtain an
exact sequence:
\begin{equation}\label{eq:conn-sum-modTheta-exact}
(\mathbb{F}(<\sigma>)/(\Theta))\rightarrow
(\mathbb{F}(K)/(\Theta)\oplus \mathbb{F}(L)/(\Theta)) \rightarrow
(\mathbb{F}(K\cup_{\sigma}L)/(\Theta))\rightarrow 0
\end{equation}
where in the middle term we used distributivity of $\otimes$ and
$\oplus$.  Note that $\mathbb{F}(<\sigma>)/(\Theta)\cong \mathbb{F}$
is concentrated in degree $0$ and that
$(\mathbb{F}(K\#_{\sigma}L)/(\Theta))_{<d} \cong
(\mathbb{F}(K\cup_{\sigma}L)/(\Theta))_{<d}$. Thus, for $0<i\leq
d/2$ we obtain the following commutative diagram:
\begin{equation}\label{eq:conn-sum-omega-isom}
\begin{CD}
(\frac{\mathbb{F}(K\#_{\sigma}L)}{(\Theta)})_i@>\cong >>(\frac{\mathbb{F}(K\cup_{\sigma}L)}{(\Theta)})_i@>\cong >>(\frac{\mathbb{F}(K)}{(\Theta)})_i\bigoplus (\frac{\mathbb{F}(L)}{(\Theta)})_{i}\\
@VV\omega^{d-2i} V @VV\omega^{d-2i} V @VV\omega^{d-2i}\oplus \omega^{d-2i}V \\
(\frac{\mathbb{F}(K\#_{\sigma}L)}{(\Theta)})_{d-i}@>\cong >>(\frac{\mathbb{F}(K\cup_{\sigma}L)}{(\Theta)})_{d-i}@>\cong >>(\frac{\mathbb{F}(K)}{(\Theta)})_{d-i}\bigoplus (\frac{\mathbb{F}(L)}{(\Theta)})_{d-i} \\
\end{CD}
\end{equation}
where the right vertical arrow is an isomorphism by assumption.
Hence, the left vertical arrow is an isomorphism as well, meaning
that $\omega$ is an $i$-Lefschetz element of
$\mathbb{F}[K\#_{\sigma}L]/(\Theta)$ for $0<i\leq d/2$.

For $i=0$, as $K\# L$ is Cohen-Macaulay with l.s.o.p. $\Theta$ and
$h_d=1$, then there exists a  $0$-Lefschetz element $\tilde{\omega}$
(i.e. $\tilde{\omega}^{d}\neq 0$. This is equivalent to
$[2,d+1]\in\Delta^s(K\# L)$, which reflects the fact that $K\# L$
has non-vanishing top homology.). By Lemma \ref{lemOmegaZariski} the
sets of $0$-Lefschetz elements and of $(0<)$-Lefschetz elements are
Zariski open. The fact that they are nonempty implies that so is
their intersection, i.e. $K\# L$ is HL. $\square$
\\
\textbf{Remark}: (2) follows also from the symmetric case of
Corollary \ref{corKcupL}. The proof given here in our 'special case'
is simpler.

\begin{cor}\label{cor:connected-sum}
Let $K$ and $L$ be HL simplicial spheres of the same dimension $d$.
Then their connected sum $K\# L$ is also an HL $d$-sphere. $\square$
\end{cor}
\begin{cor}\label{cor:WLconnected-sum}
Let $K$ and $L$ be WL spheres of the same dimension, which intersect in a
common facet $<\sigma>=K\cap L$. Then $K\#_{\sigma}L$
is WL.
\end{cor}
$Proof$: Imitate the proof of Theorem \ref{thm:connected-sum}.
$\square$

\section{Swartz lifting theorem and beyond}
In this section we show that for proving Conjecture \ref{conj-g} it
suffices to show that $y_{d+1}: H(K,\Theta)_{\lfloor
d/2\rfloor}\rightarrow H(K,\Theta)_{\lceil d/2\rceil}$ is an
isomorphism for $K$ a homology $(d-1)$-sphere with $d$ odd and
generic l.s.o.p. $\Theta$ and $y_{d+1}$ in $A_1$. We end this
section by stating a stronger conjecture about the structure of the
set of pairs $(\Theta,\omega)$ of an l.s.o.p. and a $\lfloor
d/2\rfloor$-Lefschetz element (stronger than being nonempty), which
hopefully would be easier to prove.

Consider the multiplication maps $\omega_i:
H(K,\Theta)_{i}\longrightarrow H(K,\Theta)_{i+1}$, $m\mapsto
\omega_i m$ where $\omega_i \in A_1$. Let $\ddim(K)=d-1$. Denote by
$\Omega_{UWL}(K,i)$ the set of all $(\Theta,\omega_i)\in
A_1^{\ddim(K)+2}$ such that $\Theta$ is an l.s.o.p. of
$\mathbb{R}[K]$, $\mathbb{R}[K]$ is a free
$\mathbb{R}[\Theta]$-module, and $\omega_i: H(K)_{i}\longrightarrow
H(K)_{i+1}$ is injective for $i< d/2$ and surjective for $i\geq
d/2$. Denote by $\Omega_{HL}(K,i)$ the set of all
$(\Theta,\omega)\in (A_{K_0})_1^{d+1}$ such that $\Theta$ is an
l.s.o.p. of $\mathbb{R}[K]$, $\mathbb{R}[K]$ is a free
$\mathbb{R}[\Theta]$-module, and $\omega^{d-2i}:
H(K)_{i}\longrightarrow H(K)_{d-i}$ is injective ($0\leq i\leq
\lfloor d/2\rfloor$). For $d$ odd $\Omega_{UWL}(K,\lfloor
d/2\rfloor)=\Omega_{HL}(K,\lfloor d/2\rfloor)$, which we simply
denote by $\Omega(K,\lfloor d/2\rfloor)$.

The following was proved by Swartz \cite{Swartz}, Proposition 3.6
for $\Omega_{HL}(K,i)$; similar arguments can be used to prove the
same conclusion for $\Omega_{UWL}(K,i)$.
\begin{lem}(Swartz)\label{lemOmegaZariski}
For every simplicial complex $K$ and for every $i$,
$\Omega_{UWL}(K,i)$ is a Zariski open set. For $0\leq i\leq \lfloor
\frac{\ddim(K)+1}{2}\rfloor$, $\Omega_{HL}(K,i)$ is a Zariski open
set. (They may be empty, e.g. if $K$ is not pure.)
\end{lem}
\begin{thm}(Swartz)\label{thmSwartz}
Let $d\geq 1$. If for every homology $2d$-sphere $L$, $\Omega(L,d)$
is nonempty, then for every $t>2d$ and for every homology $t$-sphere
$K$, $\Omega_{UWL}(K,m)$ is nonempty for every $m\leq d$.
\end{thm}
$Proof$: By \cite{Swartz-SpheresToManifolds}, Theorem 4.26 and
induction on $t$, $\Omega_{UWL}(K,(t+1)-(d+1))$ is nonempty, i.e.
multiplication $\omega :H(K)_{t-d}\rightarrow H(K)_{t-d+1}$ is
surjective for a generic $\omega\in A_1$. As the ring $H(K)$ is
standard, $\Omega_{UWL}(K,(t+1)-(m+1))$ is nonempty for every $m\leq
d$. Hence, for the canonical module $\Omega(K)$, multiplication by a
generic degree $1$ element
$\omega:(\Omega(K)/\Theta\Omega(K))_{m}\rightarrow
(\Omega(K)/\Theta\Omega(K))_{m+1}$ is injective in the first $d$
degrees. As $K$ is a homology sphere, $\Omega(K)\cong
\mathbb{R}[K]$, hence $\Omega_{UWL}(K,m)$ is nonempty for every
$m\leq d$. $\square$

For more information about canonical modules we refer to
\cite{StanleyGreenBook}.

Combined with Lemma \ref{lemOmegaZariski}, and the fact that a
finite intersection of Zariski nonempty open sets is nonempty, if
the conditions of Theorem \ref{thmSwartz} are met for every $d\geq
1$ then every homology sphere is unimodal WL, and hence Conjecture
\ref{conj-g} follows.

We wish to show further, that if 'all' even dimensional spheres
satisfy the condition in Theorem \ref{thmSwartz} then 'all' spheres
are HL. By 'all' we mean a family of Gorenstein$^*$ simplicial
complexes which contains all boundaries of simplices and which is
closed under joins and links (e.g. homology /simplicial /PL
spheres). The following lemma provides a step in this direction.
\begin{lem}\label{lem:WL->HL}
Let $S$ be a Gorenstein$^*$ simplicial complex with an l.s.o.p.
$\Theta_S$ over $\mathbb{R}$. If $H(S,\Theta_S)$ is $(\lfloor
\frac{\ddim S +1}{2}\rfloor)$-Lefschetz but not HL then there exists
a simplex $\sigma$ such that $S*\partial\sigma$ is of even dimension
$2j$, and for every l.s.o.p. $\Theta_{\partial \sigma}$ of $\partial
\sigma$, $\mathbb{R}[S*\partial \sigma]/(\Theta_S \cup
\Theta_{\partial \sigma})$ has no $j$-Lefschetz element; in
particular $S*\partial\sigma$ is not unimodal WL. (We would like to
obtain this conclusion for \emph{every} l.s.o.p. of $S*\partial
\sigma$!)
\end{lem}
$Proof:$ Denote the dimension of $S$ by $d-1$ and recall that
$A_{S_0}=\mathbb{R}[x_v:v\in S_0]$. By Lemma \ref{lemOmegaZariski}
$\Omega_{HL}(S,i)$ is a Zariski open set for every $0\leq i\leq
\lfloor d/2\rfloor$. The assumption that $S$ is not HL (but is
$(\lfloor \frac{d}{2}\rfloor)$-Lefschetz) implies that there exists
$0\leq i_0\leq \lfloor d/2\rfloor -1$ such that
$\Omega_{HL}(S,i_0)=\emptyset$ (as a finite intersection of Zariski
nonempty open sets is nonempty). Hence, for the fixed l.s.o.p.
$\Theta_S$ and every $\omega_S \in (A_{S_0})_1$, there exists $0\neq
m=m(\omega_S)\in H_{i_0}(S)$ such that $\omega_S^{d-2i_0}m=0$.

Let $T=S*\partial \sigma$ where $\sigma$ is the $(d-2i_0
-1)$-simplex. Note that $\ddim(\sigma)\geq 1$ (as $S$ is $(\lfloor
\frac{\ddim S +1}{2}\rfloor)$-Lefschetz), hence $\partial \sigma\neq
\emptyset$. Then $T$ is of even dimension $2d-2i_0-2$. We have seen
(Theorem \ref{thm:*}) that for any l.s.o.p. $\Theta_{\partial
\sigma}$ of $\partial \sigma$, $\Theta_T:=\Theta_S \cup
\Theta_{\partial \sigma}$ is an l.s.o.p. of $T$. Every $\omega_T\in
(A_{T_0})_1$ has a unique expansion
$\omega_T=\omega_S+\omega_{\partial \sigma}$ where $\omega_S\in
(A_{S_0})_1$ and $\omega_{\partial \sigma} \in (A_{\partial
\sigma_0})_1$. Recall the isomorphism (\ref{eq:JoinIsom}) of
$A_{T_0}$-modules $\mathbb{R}[T]/(\Theta_T)\cong
\mathbb{R}[S]/(\Theta_S)\otimes_{\mathbb{R}} \mathbb{R}[\partial
\sigma]/(\Theta_{\partial \sigma})$. Let $m(\omega_T)\in
(\frac{\mathbb{R}[T]}{(\Theta_T)})_{d-i_0-1}$ be
$$m(\omega_T):=\sum_{0\leq j\leq d-2i_0-1}(-1)^j \omega_S^{d-2i_0-1-j}m\otimes \omega_{\partial \sigma}^j 1.$$
Note that the sum $\omega_T m(\omega_T)$ is telescopic, thus
$\omega_T m(\omega_T)=\omega_S^{d-2i_0}m\otimes 1+ (-1)^{d-2i_0-1}
m\otimes \omega_{\partial \sigma}^{d-2i_0}1= 0+0=0$. For a generic
$\omega_T$, the projection of $\omega_{\partial \sigma}$ on
$\mathbb{R}[\partial \sigma]/(\Theta_{\partial \sigma})$ is nonzero,
hence so is the projection of $\omega_{\partial \sigma}^{d-2i_0-1}$,
and we get that $m(\omega_T)\neq 0$. Thus, Zariski topology tells us
that for \emph{every} $\omega_T\in (A_{T_0})_1$, there exists $0\neq
m(\omega_T)\in (\frac{\mathbb{R}[T]}{(\Theta_T)})_{d-i_0-1}$ such
that $\omega_T m(\omega_T)=0$. $\square$


We conjecture that the following stronger property holds for $\Omega(L,d)$:
\begin{conj}\label{conj:OmegaHypperplane}
Let $L$ be a homology $2d$-sphere ($d\geq 1$) on $n$ vertices. Then
$\Omega(L,d)$ intersects every hyperplane in the vector space
$A_1^{2d+2}\cong \mathbb{R}^{(2d+2)n}$.
\end{conj}
For $L$ the boundary of a simplex, the complement of $\Omega(L,d)$
is the set of all $(z_1,...,z_{2d+2})\in A_1^{2d+2}$ such that
$\det(z_1,...,z_{2d+2})=0$. As $\det\in \mathbb{R}[z_{i,j}]_{1\leq
i,j\leq 2d+2}$ is an irreducible polynomial, in particular it has no
linear factor, hence $\Omega(L,d)$ intersects every hyperplane. By
an unpublished argument of Swartz, it follows that if $L'$ is
obtained from $L$ by a bistellar move, and $\Omega(L,d)$ intersects
every hyperplane then $\Omega(L',d)$ is nonempty. We need to show
that $\Omega(L',d)$ intersects every hyperplane, in order to
conclude that the $g$-conjecture holds for PL-spheres.
'Unfortunately', $\Omega(L,d)$ may not be connected, as its
complement is a codimension one algebraic variety.


\section{Lefschetz properties and Stellar subdivisions}\label{sec:WL&Stellar}
Roughly speaking, we will show that Stellar subdivisions preserve
the HL property.
\begin{prop}\label{propAlgContruction}
Let $K$ be a simplicial complex. Let $K'$ be obtained from $K$ by
identifying two distinct vertices $u$ and $v$ in $K$, i.e. $K'=\{T:
u\notin T \in K\}\cup\{(T\setminus \{u\})\cup \{v\}: u\in T \in
K\}$. Let $d\geq 2$. Assume that $\{d+2,d+3,...,2d+1\}\notin
\Delta(K')$ and that $\{d+1,d+2,...,2d-1\}\notin \Delta(\lk(u,K)\cap
\lk(v,K))$. Then $\{d+2,d+3,...,2d+1\}\notin \Delta(K)$. (Shifting
is over $\mathbb{R}$.)
\end{prop}
\textbf{Remark}: The case $d=2$ and $\ddim(K)=1$ follows from
Lemmata \ref{Whiteley} (symmetric case) and \ref{extWHiteley}
(exterior case).
\\
\\
\emph{Proof for symmetric shifting}: (with Eric Babson) Let $\psi: K_0\longrightarrow \mathbb{R}^{2d}$ be a generic
embedding, i.e. all minors of the representing matrix w.r.t. a fixed
basis are nonzero. It induces the following map:
\begin{eqnarray}\label{eq:HighRigidity}
\psi_K^{2d}: \oplus_{T\in K_{d-1}}\mathbb{R}T \longrightarrow
\oplus_{F\in \binom{K_0}{d-1}}\mathbb{R}^{2d}/\sspan(\psi(F)),\nonumber \\
1T\mapsto \sum_{F\in \binom{K_0}{d-1}}\delta_{F\subseteq T} \overline{\psi(T\setminus F)} F
\end{eqnarray}
where $\delta_{F\subseteq T}$ equals $1$ if $F\subseteq T$ and $0$
otherwise.

Recall that $\{d+2,d+3,...,2d+1\}\notin \Delta^s(K)$ iff $y_{2d+1}^d
\notin GIN(K)$, where $Y=\{y_i\}_i$ is a generic basis for $A_1$,
$A=\mathbb{R}[x_v: v\in K_0]$.
 By Lee \cite{Lee} Theorems 10,12,15 and Tay, White and
Whiteley \cite{TayWhiteWhiteley-skel1} Proposition 5.2, $y_{2d+1}^d
\notin GIN(K)$ iff $\Ker \phi_K^{2d}=0$ for some $\phi:
K_0\longrightarrow \mathbb{R}^{2d}$ (equivalently, every $\phi$ in
some Zariski non-empty open set of embeddings).

Consider the following degenerating map: for $0<t\leq 1$ let
$\psi_t: K_0\longrightarrow \mathbb{R}^{2d}$ be defined by
$\psi_t(i)=\psi(i)$ for every $i\neq u$ and
$\psi_t(u)=\psi(v)+t(\psi(u)-\psi(v))$. Thus $\psi_1=\psi$,
and $\lim_{t\mapsto 0}(\psi_t(u)-\psi_t(v))=\psi(u)-\psi(v)$. Let
$\psi_0=\lim_{t\mapsto 0}\psi_t$.

Let $\psi_{K,t}^{2d}: \oplus_{T\in K_{d-1}}\mathbb{R}T
\longrightarrow \oplus_{F\in
\binom{K_0}{d-1}}\mathbb{R}^{2d}/\sspan(\psi_t(F))$ be the map induced
by $\psi_t$; thus $\psi_{K,1}^{2d}=\psi_{K}^{2d}$. Denote
$\psi_0^{2d}=\lim_{t\mapsto 0}\psi_{K,t}^{2d}$. Assume for a moment
that $\psi_0^{2d}$ is injective. Then for a small enough perturbation of the
entries of a representing matrix of $\psi_0^{2d}$, the columns of the resulted
matrix would be independent, i.e. the corresponding linear
transformation would be injective. In particular, there would exist
an $\epsilon>0$ such that for every $0<t<\epsilon$, $\Ker
\psi_{K,t}^{2d}=0$, and hence for every $\phi: K_0\longrightarrow
\mathbb{R}^{2d}$ in some Zariski non-empty open set of embeddings,
$\Ker \phi_K^{2d}=0$. Thus, the following Lemma \ref{lemCol(A)ind} completes the proof.
$\square$

\begin{lem}\label{lemCol(A)ind}
$\psi_0^{2d}$ is injective for a non-empty Zariski open set of
embeddings $\psi: K_0\longrightarrow \mathbb{R}^{2d}$.
\end{lem}
$Proof:$ For every $0<t\leq 1$ and every $F$ such that
$\{u,v\}\subseteq F \in \binom{K_0}{d-1}$,
$\sspan(\psi_t(F))=\sspan(\psi(F))$, and hence in the range of $\psi_0^{2d}$ we
mod out by $\sspan(\psi(F))$ for summands with such $F$. For summands
of $\{u,v\}\nsubseteq F \in \binom{K_0}{d-1}$,
we mod out by $\sspan(\psi_0(F))$. Note that for $T$ such that
$\{u,v\}\subseteq T\in K_{d-1}$,
$$\psi_0^{2d}(T)|_{T\setminus
v}=(\psi(u)-\psi(v)) + \sspan(\psi(T\setminus
u))=-\psi_0^{2d}(T)|_{T\setminus u}.$$ For a linear transformation
$C$, denote by $[C]$ its representing matrix w.r.t. given bases. In
$[\psi_0^{2d}]$ bases are indexed by sets as in
(\ref{eq:HighRigidity}). First add rows $F'\uplus\{u\}$ to rows
$F'\uplus\{v\}$, then delete the rows $F$ containing $u$, to obtain
a matrix $[B]$, of a linear transformation $B$. In particular, we
delete all rows $F$ such that $\{u,v\}\subseteq F$.

Note that $K'_0=K_0\setminus \{u\}$, thus, for the obvious bases,
$[B]$ is obtained from $[(\psi|_{K'_0})_{K'}^{2d}]$ by doubling the
columns indexed by $T'\uplus\{v\}\in K'_{d-1}$ where both
$T'\uplus\{v\},T'\uplus\{u\}\in K_{d-1}$, and by adding a zero
column for every $T'\uplus\{u,v\}\in K_{d-1}$. For short, denote
$\psi_{K'}^{2d}=(\psi|_{K'_0})_{K'}^{2d}$. More precisely, the
linear maps $B$ and $\psi_{K'}^{2d}$ are related as follows: they
have the same range. The domain of $B$ is $\dom(B)=\dom(\psi_0^{2d})=D_1\oplus
D_2\oplus D_3$ where
\newline $D_1=\oplus\{\mathbb{R}T: T\in K_{d-1}, \{u,v\}\nsubseteq T, (u\in T)\Rightarrow (T\setminus u)\cup v \notin K\}$,
\newline $D_2=\oplus\{\mathbb{R}T: T\in K_{d-1}, u\in T, v\notin T, (T\setminus u)\cup v \in K\}$,
\newline $D_3=\oplus\{\mathbb{R}T: T\in K_{d-1}, \{u,v\}\subseteq T\}$.
\newline
For a base element $1T$ of $D_1$, let $T'\in K'$ be obtained from
$T$ by replacing $u$ with $v$. Then $B(1T)=\psi_{K'}^{2d}(1T')$;
thus $\Ker B|_{D_1}\cong \Ker \psi_{K'}^{2d}$. For a base element $1T$ of
$D_2$, $B(1T)=\psi_{K'}^{2d}(1((T\setminus u)\cup v))$, and
$B|_{D_3}=0$.

Assume we have a linear dependency $\sum_{T\in K_{d-1}}\alpha_T
\psi_0^{2d}(T)=0$.
By assumption, $\{d+2,d+3,...,2d+1\}\notin \Delta^{s}(K')$, hence $\Ker
\psi_{K'}^{2d}=0$, thus $\alpha_T=0$ for every base element $T$
except possibly for $T\in D_3$ and for $T'\uplus \{u\}, T'\uplus
\{v\} \in K_{d-1}$, where $\alpha_{T'\uplus \{u\}}=-\alpha_{T'\uplus
\{v\}}$.

Let $\psi_0^{2d}|_{\res}$ be the restriction of $\psi_0^{2d}$ to the
subspace spanned by the base elements $T$ such that $v\in T$ and for
which it is (yet) not known that $\alpha_T=0$, followed by
projection into the subspace spanned by the $F\in \binom{K_0}{d-1}$
coordinates where $v\in F$ - just forget the other coordinates. As
$\psi_0^{2d}(T)|_F=0$ whenever $F \ni v \notin T$, if
$\psi_0^{2d}|_{\res}$ is injective, then $\alpha_T=0$ for all $T \in
K_{d-1}$. Thus, the Lemma \ref{lemCol(A|_res)ind} below completes
the proof. $\square$

\begin{lem}\label{lemCol(A|_res)ind}
$\psi_0^{2d}|_{\res}$ is injective for a non-empty Zariski open set of
embeddings $\psi: K_0\longrightarrow \mathbb{R}^{2d}$.
\end{lem}
$Proof:$ Let $G=(\{u\}*(\lk(u,K)\cap \lk(v,K)))_{\leq d-2}$. Note that
$v$ appears in the index set of every row and every column of
$[\psi_0^{2d}|_{\res}]$. Omitting $v$ from the indices of both of the bases
used to define $\psi_0^{2d}|_{\res}$, we notice that
$$\psi_0^{2d}|_{\res}\cong\overline{\psi_0^{2d}|_{\res}}: \oplus_{T\in G_{d-2}}\mathbb{R}T \longrightarrow
\oplus_{F\in
\binom{G_0}{d-2}}\mathbb{R}^{2d}/\sspan(\psi(F\uplus\{v\}))=$$
$$\oplus_{F\in \binom{G_0}{d-2}}(\mathbb{R}^{2d}/\sspan(\psi(v)))/\overline{span(\psi(F))},$$
$$1T\mapsto \sum_{F\in \binom{G_0}{d-2}}\delta_{F\subseteq T} \overline{\psi(T\setminus F)}F$$
where $\delta_{F\subseteq T}$ equals $1$ if $F\subseteq T$ and $0$
otherwise, and $\overline{\sspan(\psi(F))}$ is the image of
$\sspan(\psi(F))$ in the quotient space
$\mathbb{R}^{2d}/\sspan(\psi(v))$.

Consider the projection $\pi: \mathbb{R}^{2d} \longrightarrow
\mathbb{R}^{2d}/\sspan(\psi(v)) \cong \mathbb{R}^{2d-1}$. Let
$\bar{\psi}=\pi \circ \psi|_{G_0}: G_0 \longrightarrow
\mathbb{R}^{2d-1}$, and $\bar{\psi}_G^{2d-1}$ be the induced map
as defined in (\ref{eq:HighRigidity}).
Then $\pi$ induces $\pi_*\overline{\psi_0^{2d}|_{\res}} = \bar{\psi}_G^{2d-1}$.

By assumption, $\{d+1,...,2d-1\}\notin \Delta^{s}(\lk(u,K)\cap
\lk(v,K))$. As symmetric shifting commutes with constructing a cone
(Kalai \cite{skira} Theorem 2.2.8, and Babson, Novik and Thomas
\cite{Babson-Novik-Thomas-Cone} Theorem 3.7), $\{d+2,...,2d\} \notin
\Delta^{s}(G)$. Hence $y_{2d}^{d-1} \notin \GIN(G)$,
and by Lee \cite{Lee}, $\Ker \phi_G^{2d-1}=0$ for a generic $\phi$.
Thus, all liftings $\psi: K_0\longrightarrow \mathbb{R}^{2d}$ such
that $\bar{\psi}=\phi$ satisfy $\Ker \psi_0^{2d}|_{\res}\cong \Ker
\phi_G^{2d-1}=0$, and this set of liftings is a non-empty Zariski
open set. $\square$
\\
\textbf{Remark}: Clearly the set of all $\psi$ such that
$\psi_K^{2d}$ is injective is Zariski open. We exhibited conditions
under which it is non-empty.
\\
\\
\emph{Proof for exterior shifting}: The proof is similar to the
proof for the symmetric case. We indicate the differences.
$\psi:K_0\rightarrow \mathbb{R}^{d+1}$ defines the first $d+1$
generic $f_i$'s w.r.t. the $e_i$'s basis of $\mathbb{R}^{|K_0|}$ and
induces the following map:
\begin{equation}\label{eq:extHighRigidity}
\psi^{d+1}_{K,\ext}: \oplus_{T\in K_{d-1}}\mathbb{R}T \longrightarrow
\oplus_{1\leq i\leq d+1}\oplus_{F\in \binom{K_0}{d-1}}\mathbb{R}F,\ m\mapsto(f_1\lfloor m,...,f_{d+1}\lfloor m)
\end{equation}
By Proposition \ref{prop.1}, $\Ker \psi^{d+1}_{K,\ext}=\cap_{1\leq
i\leq d+1}\Ker_{d-1}f_i\lfloor =
\cap_{R<\{d+2,...,2d+1\}}\Ker_{d-1}f_R\lfloor$, hence, by
shiftedness, $\{d+2,...,2d+1\}\notin \Delta^e(K) \Leftrightarrow
\Ker \psi^{d+1}_{K,\ext}=0$.

Replacing $\psi(u)$ by $\psi(v)$ induces a map
$$\psi^{d+1}_{K,u}:\oplus_{T\in K_{d-1}}\mathbb{R}T \longrightarrow
\oplus_{1\leq i\leq d+1}\oplus_{F\in \binom{K_0}{d-1}}\mathbb{R}F.$$
By perturbation, if $\Ker \psi^{d+1}_{K,u}=0$ then $\Ker
\psi^{d+1}_{K,\ext}=0$ for generic $\psi$.

Let $[B_{\ext}]$ be obtained from the matrix $[\psi^{d+1}_{K,u}]$ by
adding the rows $F'\uplus u$ to the corresponding rows $F'\uplus v$
and deleting the rows $F$ with $\{u,v\}\subseteq F$. The domain of
$B_{\ext}$ is $D_1\oplus D_2\oplus D_3$ as for $B$ in the symmetric
case. For a base element $1T$ of $D_1$, let $T'\in K'$ be obtained
from $T$ by replacing $u$ with $v$. Then
$B_{\ext}(1T)=\psi_{K',\ext}^{d+1}(1T')$; thus $\Ker
B_{\ext}|_{D_1}\cong \Ker \psi_{K',\ext}^{d+1}$. For a base element
$1T$ of $D_2$, $B_{\ext}(1T)=\psi_{K',\ext}^{d+1}(1((T\setminus
u)\cup v))$, and as we may number $v=1,u=2$ then $B|_{D_3}=0$ (the
rows of $F'\uplus u$ and of $F'\uplus v$ have opposite sign in
$\psi^{d+1}_{K,u}$). Now we can repeat the arguments showing that
$\Ker \psi^{2d}_0=0$ by considering $B$ in the symmetric case, to
show that $\Ker \psi^{d+1}_{K,u}=0$ by considering $B_{\ext}$.
$\square$

\begin{cor}\label{cor:contraction}
Let $K$ be a $2d$-sphere for some $d\geq 1$, and let $a,b\in K$ be
two vertices which satisfy the \emph{Link Condition}, i.e that
$\rm{lk}(a,K)\cap\rm{lk}(b,K)=\rm{lk}(\{a,b\},K)$.
Let $K'$ be obtained from $K$ by contracting $a\mapsto b$.
Then:

(1) $K'$ is a $2d$-sphere, PL homeomorphic to $K$ (see Theorem \ref{thmDey}).

(2) If $K'$ is $d$-Lefschetz and $\lk(\{a,b\},K)$ is
$(d-1)$-Lefschetz, then $K$ is $d$-Lefschetz (by Proposition \ref{propAlgContruction}).
$\square$
\end{cor}

Let $K$ be a simplicial complex. Its \emph{Stellar subdivision at a
face $T\in K$} is the operation $K\mapsto K'$ where
$K'=\Stellar(T,K):=(K\setminus \st(T,K))\cup(\{v_T\}*\partial T*
\lk(T,K))$, where $v_T$ is a vertex not in $K$. Note that for $u\in
T\in K$, $u,v_T\in K'$ satisfy the Link Condition and their
identification results in $K$. Further,
$\lk(\{u,v_T\},K')=\lk(u,\partial T* \lk(T,K))=\partial(T\setminus
u)* \lk(T,K)$.

\begin{thm}\label{thm:Stellar(HL)=HL}
Let $S$ be a homology sphere and $F\in S$. If $S$ and $\lk(F,S)$ are HL then $\Stellar(F,S)$ is HL.
\end{thm}
$Proof:$ Let $T=\Stellar(F,S)$, denote its dimension by $d-1$, and
assume by contradiction that $T$ is not HL. As we have seen in the
proof of Lemma \ref{lem:WL->HL}, there exists $0\leq i_0\leq \lfloor
d/2\rfloor$ such that $\Omega_{HL}(T,i_0)=\emptyset$. First we show
that $i_0\neq \lfloor d/2\rfloor$: for $d$ even this is obvious. For
$d$ odd, note that for $u\in F$ the contraction $v_F\mapsto u$ in
$T$ results in $S$, which is $\lfloor d/2\rfloor$-Lefschetz.
Further, the $(d-3)$-sphere
$\lk(\{v_F,u\},T)=\lk(F,S)*\partial(F\setminus \{u\})$ is HL by
Theorem \ref{thm:*}, and in particular is $(\lfloor
d/2\rfloor-1)$-Lefschetz. Thus, by Corollary \ref{cor:contraction}
$T$ is $\lfloor d/2\rfloor$-Lefschetz, and hence $0\leq i_0\leq
\lfloor d/2\rfloor-1$.

Let $L=T*\partial \sigma$, where $\sigma$ is the
$(d-2i_0-1)$-simplex (then $L$ has even dimension $2d-2i_0-2$). By
Lemma \ref{lem:WL->HL}, for any two l.s.o.p.'s $\Theta_T$ and
$\Theta_{\partial \sigma}$ of $\mathbb{R}[T]$ and
$\mathbb{R}[\partial \sigma]$ respectively, $\mathbb{R}[L]/(\Theta_T
\cup \Theta_{\partial \sigma})$ has no $(d-i_0-1)$-Lefschetz
element.

On the other hand, we shall now prove the existence of such
l.s.o.p.'s and a $(d-i_0-1)$-Lefschetz element, to reach a
contradiction. This requires a close look on the proof of
Proposition \ref{propAlgContruction}.

Note that $L=\Stellar(F,S*\partial \sigma)$, and that for $u\in F$
the contraction $v_F\mapsto u$ in $L$ results in $S*\partial
\sigma$. Further, $\lk(\{v_F,u\},L)=\lk(F,S)*\partial(F\setminus
\{u\})* \partial \sigma$.

Applying Zariski topology considerations to subspaces of the space
of embeddings $\{f:L_0\rightarrow \mathbb{R}^{2d-2i_0}\}\cong
\mathbb{R}^{|L_0|\times (2d-2i_0)}$, we now show that there exists
an embedding $\psi: L_0\longrightarrow \mathbb{R}^{d}\oplus
\mathbb{R}^{d-2i_0-1}\oplus \mathbb{R}$ such that the following
three properties hold \emph{simultaneously}:

(1) $\psi|_{S_0} \subseteq \mathbb{R}^{d}\oplus 0 \oplus \mathbb{R}$
and induces an l.s.o.p. $\Theta_S$ of $\mathbb{R}[S]$ and an HL
element $\omega_S$ of $\mathbb{R}[S]/(\Theta_S)$; $\psi|_{\sigma_0}
\subseteq 0 \oplus \mathbb{R}^{d-2i_0-1}\oplus \mathbb{R}$ and
induces an l.s.o.p. $\Theta_{\partial \sigma}$ of
$\mathbb{R}[\partial \sigma]$ and an HL element $\omega_{\partial
\sigma}$ of $\mathbb{R}[S]/(\Theta_{\partial \sigma})$. By Theorem
\ref{thm:*}, $\omega_S+\omega_{\partial \sigma}$ is an HL element of
$\mathbb{R}[S*\partial \sigma]/(\Theta_S \cup \Theta_{\partial
\sigma})$.

In matrix language, the first $2d-2i_0-1$ columns of
$[\psi|_{S_0\cup \sigma_0}]$ form an l.s.o.p. of
$\mathbb{R}[S*\partial \sigma]$, and its last column is the
corresponding HL element.

(2) $0\neq \psi(v_F)\in \mathbb{R}^{d}\oplus 0 \oplus \mathbb{R}$
induces a map $\pi: \mathbb{R}^{2d-2i_0}\rightarrow
\mathbb{R}^{2d-2i_0}/\sspan \psi(v_F) \cong \mathbb{R}^{2d-2i_0-1}$
such that $\pi\circ\psi|_{S_0\cup \sigma_0}$ induces an element in
$\Omega(G,d-i_0-2)$ for $G=\{u\}* \lk(\{v_F,u\},L)$.

To see this, consider e.g. an embedding $\psi'$ with
$\psi'(v_F)=(1,0,...,0)$, $\psi'(u)=(0,1,0,...,0)$, $\psi'(s)$
vanishes on the first two coordinates for any $s\in S_0\setminus
\{u\}$ and in addition $[\psi']$ vanishes on all entries on which we
required in (1) that $[\psi]$ vanishes. By Theorem \ref{thm:*} there
exists such $\psi'$ so that its composition with the projection
$\pi':\mathbb{R}^{2d-2i_0}\rightarrow \mathbb{R}^{2d-2i_0}/\sspan
\{\psi(v_F),\psi(u)\}$ induces a pair $(\Theta,\omega)$ of an
l.s.o.p. and an HL element for
$\lk(\{v_F,u\},L)=\lk(F,S)*\partial(F\setminus \{u\})*
\partial \sigma$. By adding $x_u$ to this l.s.o.p. we obtain an l.s.o.p.
for $G$ where $\omega: H(G)_{d-i_0-2}\rightarrow H(G)_{d-i_0-1}$ is
injective. Now perturb $\psi'$ to obtain $\psi$ for which property
(2) hold.

The restriction of maps $\psi$ with property (2) to $\st(F,S)_0\cup
\{v_F\}$ is a nonempty Zariski open set in the space of embeddings
$\{f: \st(F,S)_0\cup \{v_F\}\rightarrow \mathbb{R}^{d}\oplus 0
\oplus \mathbb{R}\}$. The restriction of maps $\psi$ with property
(1) to $S_0$ is a nonempty Zariski open set in the space of
embeddings $\{f: S_0\rightarrow \mathbb{R}^{d}\oplus 0 \oplus
\mathbb{R}\}$. Hence, their projections on the linear subspace $\{f:
\st(F,S)_0\rightarrow \mathbb{R}^{d}\oplus 0 \oplus \mathbb{R}\}$
are nonempty Zariski open sets (in this subspace). The intersection
of these projections is again a nonempty Zariski open set, thus
there are maps $\psi$ for which both properties (1) and (2) hold.

(3) $\psi|_{S_0\cup \{v_F\}} \subseteq \mathbb{R}^{d}\oplus 0 \oplus
\mathbb{R}$ and the first $d$ columns of $[\psi]$ induce an l.s.o.p.
$\Theta_T$ of $\mathbb{R}[T]$.

The set of restrictions $\psi|_{T_0}$ of maps $\psi$ with property
(3) is nonempty Zariski open in the subspace $\{f: T_0\rightarrow
\mathbb{R}^{d}\oplus 0 \oplus 0\}$; hence, so is its projection on
the linear subspace $\{f: \st(F,S)_0\rightarrow \mathbb{R}^{d}\oplus
0 \oplus 0\}$. By similar considerations to the above, there are
maps $\psi$ for which all the properties (1), (2) and (3) hold.

The proof of Proposition \ref{propAlgContruction} together with
properties (1) and (2) tell us that for small enough $\epsilon$, the
map $\psi": L_0\longrightarrow \mathbb{R}^{2d-2i_0}$ defined by
$\psi"(v_F)=\psi(u)+\epsilon (\psi(v_F)-\psi(u))$ and
$\psi"(v)=\psi(v)$ for every other vertex $v\in L_0$, satisfy $\Ker
\psi"_{L}^{2d-2i_0}=0$ (see equation (\ref{eq:HighRigidity}) for the
definition of this map). As a nonempty Zariski open set is dense, by
looking on the subspace $\{f:T_0\rightarrow \mathbb{R}^{d}\oplus 0
\oplus \mathbb{R}\}\}$ , we can take $\psi(v_F)$ and $\epsilon$ such
that $\psi"$ satisfies property (3) as well.

Thus, the first $d$ columns of $[\psi"]$ induce an l.s.o.p.
$\Theta_T$ of $T$, the next $d-i_0-1$ columns induce an l.s.o.p.
$\Theta_{\partial \sigma}$ of $\partial \sigma$, and the last column
of $[\psi"]$ is a $(d-i_0-1)$-Lefschetz element of
$\mathbb{R}[L]/(\Theta_T \cup \Theta_{\partial \sigma})$. This
contradicts our earlier conclusion, which was based on assuming that
the assertion of this theorem is incorrect. $\square$

\begin{cor}\label{cor:Stellar}
Let $\mathcal{S}$ be a family of homology spheres which is closed
under taking links and such that all of its elements are HL. Let
$\mathbb{S}=\mathbb{S}(\mathcal{S})$ be the family obtained from
$\mathcal{S}\cup \{\partial \sigma^n: n\geq 1\}$ by taking the
closure under the operations: (0) taking links; (1) join; (2)
Stellar subdivisions. Then every element in $\mathbb{S}$ is HL.
\end{cor}
$Proof:$ We prove by double induction - on dimension, and on the
sequence of operations of type (0),(1) and (2) which define $S\in
\mathbb{S}$ - that $S$ and all its face links are HL. Let us call
$S$ with this property \emph{hereditary HL}.

Note that every $S\in \mathcal{S}$, every boundary of a simplex, and
every (homology) sphere of dimension $\leq 2$, is hereditary HL.
This provides the base of the induction.

Clearly if $S$ is hereditary HL, then so are all of its links, as
$\lk(Q,(\lk(F,S))=\lk(Q\uplus F,S)$. If $S$ and $S'$ are hereditary
HL then by Theorem \ref{thm:*} so is $S*S'$ (here we note that every
$T\in S*S'$ is of the form $T=F\uplus F'$ where $F\in S$ and $F'\in
S'$, and that $\lk(T,S*S')=\lk(F,S)*\lk(F',S')$). We are left to
show that if $F\in S$ and $S$ is hereditary HL, then so is
$T:=\Stellar(F,S)$. Assume $\ddim F\geq 1$, otherwise there is
nothing to prove. First we note that by the induction hypothesis for
every $v\in T_0$, $\lk(v,T)$ is hereditary HL:
\\ Case $v=v_F$: $\lk(v_F,T)=\lk(F,S)*\partial F$ is hereditary HL
by Theorem \ref{thm:*}, as argued above.
\\ Case $v\in F$: $\lk(v,T)=\Stellar(F\setminus \{v\},\lk(v,S))$ is hereditary HL
by the induction hypothesis on the dimension.
\\ Case $v\notin F$, $v\neq v_F$ and $F\in \lk(v,S)$:
$\lk(v,T)=\Stellar(F,\lk(v,S))$ is hereditary HL by the induction
hypothesis on the dimension.
\\ Otherwise: $\lk(v,T)=\lk(v,S)$ is hereditary HL.

We are left to show that $T$ is HL: $S$ is HL, and for $u\in F$
$\lk(\{v_F,u\},T)=\lk(F,S)*\partial(F\setminus \{u\})$ is HL by
Theorem \ref{thm:*}. Thus, by Theorem \ref{thm:Stellar(HL)=HL} $T$
is HL, and together with the above, $T$ is hereditary HL. $\square$
\\
 \\
\textbf{Remark}: The barycentric subdivision of a simplicial complex
$K$ can be obtained by a sequence of Stellar subdivisions: order the
faces of $K$ of dimension $>0$ by weakly decreasing size, and
perform Stellar subdivisions at those faces according to this order;
the barycentric subdivision of $K$ is obtained. Brenti and Welker
\cite{Brenti-Welker}, Corollary 3.5, showed that the $h$-polynomial
of the barycentric subdivision of a Cohen-Macaulay complex has only
simple and real roots, and hence is unimodal. In particular,
barycentric subdivision preserves non-negativity of the $g$-vector
for spheres with all links being HL. The above corollary shows that
the hereditary HL property itself is preserved.

\section{Open problems}
\begin{enumerate}
\item Show that for any $d\geq 1$ and any homology $2d$-sphere $L$,
$\Omega(L,d)$ intersects every hyperplane in the vector space
$\mathbb{R}^{|L_0|\times (2d+2)}$ of $2d+2$ degree one forms.

In Theorem \ref{thmSwartz} we have seen that to conclude the $g$-conjecture  $\Omega(L,d)\neq \emptyset$ is enough, but this stronger conjecture may be easier to prove in the PL case, by using bistellar moves. It tries to correct an unpublished argument of Swartz.

\item \emph{Shifting theoretic lower bound relation:} In Example \ref{ex:stacked} we computed the algebraic
shifting of a stacked $(d-1)$-sphere on $n$ vertices, denoted
$\Delta(S(d,n))$. Prove that if $K$ is a homology $(d-1)$-sphere on
$n$ vertices then $\Delta(S(d,n))\subseteq \Delta(K)$.

This conjecture immediately implies Barnette's lower bound theorem
for triangulated spheres. The symmetric case of this conjecture is
equivalent to the claim that the multiplication map $y_{d+1}^{d-2}:
H(K)_1\longrightarrow H(K)_{d-1}$, is an isomorphism. Rigidity
theory only tells us that $y_{d+1}: H(K)_1\longrightarrow H(K)_{2}$
is injective.

\item\label{prob:WL*HL=WL} Is the join of a unimodal WL complex with an HL complex always unimodal WL?

\item Is $\omega$ from Theorem \ref{thm:connected-sum} an HL element for $K\# L$?

\item Let $\mathbb{S}$ be a family of simplicial complexes which is closed under links and
joins and contains all boundaries of simplices (e.g. simplicial
spheres, homology spheres, PL-spheres). Prove that if all elements
of $\mathbb{S}$ are unimodal WL then all of them are HL.
\end{enumerate}


\chapter{Algebraic Shifting and the $g$-Conjecture: a Topological Approach}\label{chapter:KalaiSarkaria}
\section{Kalai-Sarkaria conjecture}
As we have seen in Section \ref{sec:Lefschetz}, if a simplicial
$(d-1)$-sphere $K$ satisfies $\Delta(K)\subseteq \Delta(d)$ then
$g(K)$ is an $M$-sequence. A stronger conjecture was stated,
independently, by Kalai and Sarkaria \cite{skira}, Conjecture 27:
\begin{conj}(Kalai, Sarkaria)\label{conj:Kalai-Sarkaria}
If $K$ is a simplicial complex with $n$ vertices and $||K||$ can be
embedded is the $(d-1)$-sphere, then $\Delta(K)\subseteq
\Delta(d,n)$. Equivalently, $T_{d-k}\notin \Delta(K)$ for every
$0\leq k\leq \lfloor d/2\rfloor$ (see equation (\ref{eq:T})).
\end{conj}
Note that by Swartz lifting theorem, Theorem \ref{thmSwartz}, to
conclude Conjecture \ref{conj-g} for simplicial spheres it is enough
to show that for $d$ odd $T_{\lceil d/2\rceil}:=\{\lceil
d/2\rceil+2,...,d+2\}\notin \Delta(K)$. The later conjecture
trivially holds for $d=1$, as $3$ points cannot be embedded into $2$
points. It holds for $d=3$ by Theorem \ref{mainThm}, as
$\{4,5\}\in\Delta(K)$ implies that the graph of $K$ has a
$K_5$-minor, hence $K$ does not embed in $S^2$. It is open for
$d=5,7,9,..$.

Sarkaria suggested to relate the Van-Kampen obstruction to
embeddability of $||K||$ to that of $||\Delta(K)||$. We recall this
obstruction in the next section, and later relate it to a notion of
minors for simplicial complexes, and to a combinatorial problem
which would imply Conjecture \ref{conj:Kalai-Sarkaria}.

\section{Van Kampen's Obstruction}
\subsection{Deleted join and $\mathbb{Z}_2$ coefficients}\label{subsection:DeletedJoin}
The presentation here is based on work of Sarkaria
\cite{SarkariaMax,SarkariaUnpub} who attributes it to Wu \cite{Wu}
and all the way back to Van Kampen \cite{VanKampen}. It is a Smith
theoretic interpretation of Van Kampen's obstructions.

Let $K$ be a simplicial complex. The join $K*K$ is the simplicial
complex $\{S^1\uplus T^2: S,T\in K\}$ (the superscript indicates
two disjoint copies of $K$). The \emph{deleted join} $K_*$ is the
subcomplex $\{S^1\uplus T^2: S,T\in K, S\cap T=\emptyset\}$. The
restriction of the involution $\tau:K*K\longrightarrow K*K$,
$\tau(S^1\cup T^2)=T^1\cup S^2$ to $K_*$ is into $K_*$. It induces
a $\mathbb{Z}_2$-action on the cochain complex
$C^*(K_*;\mathbb{Z}_2)$. For a simplicial cochain complex $C$ over
$\mathbb{Z}_2$ with a $\mathbb{Z}_2$-action $\tau$, let $C_S$ be
its subcomplex of \emph{symmetric cochains}, $\{c\in C:
\tau(c)=c\}$. Restriction induces an action of $\tau$ as the
identity map on $C_S$. Note that the following sequence is exact
in dimensions $\geq 0$:
$$0\longrightarrow C_S(K_*)\longrightarrow C(K_*)\stackrel{id+\tau}{\longrightarrow}
C_S(K_*)\longrightarrow 0$$ where $C_S(K_*)\longrightarrow C(K_*)$
is the trivial injection. (The only part of this statement that may
be untrue for a non-free simplicial cochain complex $C$ over
$\mathbb{Z}_2$ with a $\mathbb{Z}_2$-action $\tau$, is that
${id+\tau}$ is surjective.) Thus, there is an induced long exact
sequence in cohomology
$$H_S^0(K_*)\stackrel{\rm{Sm}}{\longrightarrow} H_S^1(K_*)\longrightarrow...\longrightarrow
H_S^q(K_*)\longrightarrow H^q(K_*) \longrightarrow
H_S^q(K_*)\stackrel{\rm{Sm}}{\longrightarrow} H_S^{q+1}(K_*)
\longrightarrow... .$$ Composing the connecting homomorphism
$\rm{Sm}$ $m$ times we obtain a map $\rm{Sm}^m:
H_S^0(K_*)\longrightarrow H_S^m(K_*)$. For the fundamental
$0$-cocycle $1_{K_*}$, i.e. the one which maps $\sum_{v\in
(K_*)_0}a_vv\mapsto \sum_{v\in (K_*)_0}a_v \in \mathbb{Z}_2$, let
$[1_{K_*}]$ denotes its image in $H_S^0(K_*)$.
$\rm{Sm}^m([1_{K_*}])$ is called the $m$-th \emph{Smith
characteristic class} of $K_*$, denoted also as $\rm{Sm}^m(K)$.

For any positive integer $d$ let $H(d)$ be the $(d-1)$-skeleton of
the $2d$-dimensional simplex. A well known result by Van Kampen and
Flores \cite{Flores,VanKampen} asserts that the Van Kampen
obstruction with $\mathbb{Z}$ coefficients (see the next subsection)
of $H(d)$ in dimension $(2d-1)$ does not vanish, and hence $H(d)$ is
not embeddable in the $2(d-1)$-sphere (note that the case $H(2)=K_5$
is part of the easier direction of Kuratowski's theorem). Here are
the analogous statements for $\mathbb{Z}_2$ coefficients.

\begin{thm}(Sarkaria \cite{SarkariaUnpub} Theorem 6.5, see also Wu \cite{Wu} pp.114-118.)\label{thmSmithH(d)}
For every $d\geq 1$, $\rm{Sm}^{2d-1}(1_{H(d)_*})\neq 0$.
\end{thm}

\begin{thm}(Sarkaria \cite{SarkariaUnpub} Theorem 6.4 and \cite{SarkariaMax} p.6)\label{thmSmith->nonEmb}
If a simplicial complex $K$ embeds in $\mathbb{R}^m$ (or in the
$m$-sphere) then $\rm{Sm}^{m+1}(1_{K_*})=0$.
\end{thm}
$Sketch\ of\ proof$: The definition of Smith class makes sense for
singular homology as well; the obvious map from the simplicial chain
complex to the singular one induces an isomorphism between the
corresponding Smith classes. The definition of deleted join makes
sense for subspaces of a Euclidean space as well (see e.g.
\cite{Matousek-BU}, 5.5); thus an embedding $||K||$ of $K$ into
$\mathbb{R}^m$ induces a continuous $\mathbb{Z}_2$-map from
$||K||_*$ into the join of $\mathbb{R}^m$ with itself minus the
diagonal, which is $\mathbb{Z}_2$-homotopic to the antipodal
$m$-sphere, $S^m$. The equivariant cohomology of $S^m$ over
$\mathbb{Z}_2$ is isomorphic to the ordinary cohomology of
$\mathbb{R}P^m$ over $\mathbb{Z}_2$, which vanishes in dimension
$m+1$. We get that $\rm{Sm}^{m+1}(S^m)$ maps to
$\rm{Sm}^{m+1}(1_{||K||_*})$ and hence the later equals to zero as
well. But $||K_*||$ and $||K||_*$ are $\mathbb{Z}_2$-homotopic,
hence $\rm{Sm}^{m+1}(1_{K_*})=0$. $\square$

\subsection{Deleted product and $\mathbb{Z}$ coefficients}\label{subsection:DeletedProduct}
More commonly in the literature, Van Kampen's obstruction is defined
via deleted products and with $\mathbb{Z}$ coefficients, where,
except for $2$-simplicial complexes, its vanishing is also
sufficient for embedding of the complex in a Euclidean space of
double its dimension.

The presentation of the background on the obstruction here is based
on the ones in \cite{NovikVanKampen}, \cite{Wu} and \cite{Ummel}.

Let $K$ be a finite simplicial complex. Its deleted product is
$K\times K \setminus \{(x,x): x\in K\}$, employed with a fixed-point
free $\mathbb{Z}_2$-action $\tau(x,y)=(y,x)$. It
$\mathbb{Z}_2$-deformation retracts into $K_{\times}=\cup\{S\times
T: S,T\in K, S\cap T=\emptyset\}$, with which we associate a cell
chain complex over $\mathbb{Z}$: $C_{\bullet}(K_{\times})=\bigoplus
\{\mathbb{Z}(S\times T): S\times T \in K_\times\}$ with a boundary
map $\partial(S\times T)=\partial S \times T + (-1)^{\rm{dim} S}
S\times
\partial T$, where $S\times T$ is a $\rm{dim}(S\times T)$-chain.
The dual cochain complex consists of the $j$-cochains
$C^j(K_{\times})=\rm{Hom}_{\mathbb{Z}}(C_j(K_{\times}),\mathbb{Z})$
for every $j$.

There is a $\mathbb{Z}_2$-action on $C_{\bullet}(K_{\times})$
defined by $\tau(S\times T)=(-1)^{\rm{dim}(S)\rm{dim}(T)}T\times S$.
As it commutes with the coboundary map, by restriction of the
coboundary map we obtain the subcomplexes of symmetric cochains
$C_s^{\bullet}(K_{\times})=\{c\in C^{\bullet}(K_{\times}):
\tau(c)=c\}$ and of antisymmetric cochains
$C_a^{\bullet}(K_{\times})=\{c\in C^{\bullet}(K_{\times}):
\tau(c)=-c\}$. Their cohomology rings are denoted by
$H_s^{\bullet}(K_{\times})$ and $H_a^{\bullet}(K_{\times})$
respectively. Let $H_{\rm{eq}}^{m}$ be $H_s^{m}$ for $m$ even and
$H_a^{m}$ for $m$ odd.

For every finite simplicial complex $K$ there is a unique
$\mathbb{Z}_2$-map, up to $\mathbb{Z}_2$-homotopy, into the infinite
dimensional sphere $i: K_{\times}\rightarrow S^{\infty}$, and hence
a uniquely defined map $i^*:
H_{\rm{eq}}^{\bullet}(S^{\infty})\rightarrow
H_{\rm{eq}}^{\bullet}(K_{\times})$. For $z$ a generator of
$H_{\rm{eq}}^{m}(S^{\infty})$ call
$o^m=o^m_{\mathbb{Z}}(K_{\times})=i^*(z)$ the Van Kampen
obstruction; it is uniquely defined up to a sign. It turns out to
have the following explicit description: fix a total order $<$ on
the vertices of $K$. It evaluates elementary symmetric chains of
even dimension $2m$ by
\begin{equation}\label{eqVKs}
o^{2m}((1+\tau)(S\times T))= \{^{1\ \rm{if}\ \rm{the}\
\rm{unordered}\ \rm{pair}\ \{S,T\}\ \rm{is}\ \rm{of}\ \rm{the}\
\rm{form}\ s_0<t_0<..<s_m<t_m} _{0 \ \rm{for}\ \rm{other}\
\rm{pairs}\ \{S,T\}}
\end{equation}
and evaluates elementary antisymmetric chains of odd dimension
$2m+1$ by
\begin{equation}\label{eqVKa}
o^{2m+1}((1-\tau)(S\times T))= \{^{1\ \rm{if}\ \{S,T\}\ \rm{is}\
\rm{of}\ \rm{the}\ \rm{form}\ t_0<s_0<t_1<..<t_m<s_m<t_{m+1}} _{0 \
\rm{for}\ \rm{other}\ \rm{pairs}\ \{S,T\}}
\end{equation}
where the $s_l$'s are elements of $S$ and the $t_l$'s are elements
of $T$. Its importance to embeddability is given in the following
classical result:
\begin{thm}\label{thmVK}\cite{VanKampen,Shapiro,Wu}
If a simplicial complex $K$ embeds in $\mathbb{R}^m$ then
$H_{\rm{eq}}^{\bullet}(K_{\times})\ni
o^m_{\mathbb{Z}}(K_{\times})=0$. If $K$ is $m$-dimensional and
$m\neq 2$ then $o^{2m}_{\mathbb{Z}}(K_{\times})=0$ implies that $K$
embeds in $\mathbb{R}^{2m}$.
\end{thm}

\section{Relation to minors of simplicial complexes}
\subsection{Definition of minors and statement of results}
The concept of graph minors has proved be to very fruitful.
A famous result by Kuratowski asserts that a graph can be embedded
into a $2$-sphere if and only if it contains neither of the graphs
$K_5$ and $K_{3,3}$ as minors. We wish to generalize the notion of
graph minors to all (finite) simplicial complexes in a way that
would produce analogous statements for embeddability of higher
dimensional complexes in higher dimensional spheres. We hope that
these higher minors will be of interest in future research, and
indicate some results and problems to support this hope.
\newline
Let $K$ and $K'$ be simplicial complexes. $K\mapsto K'$ is called a
\emph{deletion} if $K'$ is a subcomplex of $K$. $K\mapsto K'$ is
called an \emph{admissible contraction} if $K'$ is obtained from $K$
by identifying two distinct vertices of $K$, $v$ and $u$, such that
$v$ and $u$ are not contained in any missing face of $K$ of
dimension $\leq \rm{dim}(K)$. (A set $T$ is called a missing face of
$K$ if it is not an element of $K$ while all its proper subsets
are.) Specifically, $K'=\{T: u\notin T \in K\}\cup\{(T\setminus
\{u\})\cup \{v\}: u\in T \in K\}$. An equivalent formulation of the
condition for admissible contractions is that the following holds:
\begin{equation}\label{eqSkelLinkCond}
(\rm{lk}(v,K) \cap \rm{lk}(u,K))_{\rm{dim}(K)-2} =
\rm{lk}(\{v,u\},K).
\end{equation}
For $K$ a graph,
(\ref{eqSkelLinkCond}) just means that $\{v,u\}$ is an edge in $K$.

We say that a simplicial complex $H$ is a \emph{minor} of $K$, and
denote it by $H<K$, if $H$ can be obtained from $K$ by a sequence of
admissible contractions and deletions (the relation $<$ is a partial
order). Note that for graphs this is the usual notion of a minor.
\newline
\textbf{Remarks}: (1) In equation (\ref{eqSkelLinkCond}), the
restriction to the skeleton of dimension at most $\rm{dim}(K)-2$ can
be relaxed by restriction to the skeleton of dimension at most
$\rm{min}\{\rm{dim}(\rm{lk}(u,K)),\rm{dim}(\rm{lk}(v,K))\}-1$,
making the condition for admissible contraction \emph{local}, and
weaker. All the results and proofs in this section hold verbatim for
this notion of a minor as well.

(2) In the definition of a minor, without loss of generality we may
replace the local condition from the remark above by the following
stronger local condition, called the \emph{Link Condition} for
$\{u,v\}$:
\begin{equation}\label{eqLinkCond}
\rm{lk}(u,K)\cap \rm{lk}(v,K)=\rm{lk}(\{u,v\},K).
\end{equation}
To see this, let $K\mapsto K'$ be an admissible contraction which is
obtained by identifying the vertices $u$ and $v$ where
$\rm{dim}(\rm{lk}(u,K))\leq \rm{dim}(\rm{lk}(v,K))$. Delete from $K$
all the faces $F\uplus\{u\}$ such that $F\uplus\{u,v\}$ is a missing
face of dimension $\rm{dim}(\rm{lk}(u,K))+2$, to obtain a simplicial
complex $L$. Note that $\{u,v\}$ satisfies the Link Condition in
$L$, and the identification of $u$ with $v$ in $L$ results in $K'$.
I thank an anonymous referee for this remark.

We first relate this minor notion to Van Kampen's obstruction with $\mathbb{Z}_2$ coefficients.
\begin{thm}\label{thm o(L)}
Let $H$ and $K$ be simplicial complexes. If $H<K$ and $\rm{Sm}^m(H)\neq 0$
then $\rm{Sm}^m(K)\neq 0$.
\end{thm}

\begin{cor}\label{cor Emb}
For every $d\geq 1$, if $H(d)<K$ then $K$ is not embeddable in the
$2(d-1)$-sphere. $\square$
\end{cor}
\textbf{Remark}:\label{rem Emb} Corollary \ref{cor Emb} would also
follow from the following conjecture:
\begin{conj}\label{conjMinorEmb}
If $H<K$ and $K$ is embeddable in the $m$-sphere then $H$ is
embeddable in the $m$-sphere.
\end{conj}
The analogue of Theorem \ref{thm o(L)}
with $\mathbb{Z}$ coefficients holds:
\begin{thm}\label{thm o_Z}
Let $H$ and $K$ be simplicial complexes. If $H<K$ and
$o^m_{\mathbb{Z}}(H_{\times})\neq 0$ then
$o^m_{\mathbb{Z}}(K_{\times})\neq 0$.
\end{thm}
From Theorems \ref{thm o_Z} and \ref{thmVK} it follows that
Conjecture \ref{conjMinorEmb} is true when $2\rm{dim}(H) = m \neq 4$
(and, trivially, when $2\rm{dim}(H) < m$).

In view of Theorems \ref{thm o(L)}, \ref{thmSmithH(d)} and
\ref{thmSmith->nonEmb}, to conclude Conjecture
\ref{conj:Kalai-Sarkaria} in the symmetric case $T_{\lceil
d/2\rceil}$ where $d$ odd, and hence also the $g$-conjecture for
simplicial spheres Conjecture \ref{conj-g}, it suffices to prove the
following combinatorial conjecture:
\begin{conj}\label{conj567->Minor}
Let $K$ be a simplicial complex.
For every $d\geq 1$, if $H(d)\subseteq \Delta(K)$ then $H(d)<K$.
\end{conj}
This conjecture holds for $d=1,2$ and is otherwise open. Assume that
$d_0$ is the minimal $d$ for which Conjecture \ref{conj567->Minor}
fails, and let $K$ be a minimal counterexample w.r.t. the number of
vertices. W.l.o.g. $\rm{dim}(K)=d_0-1$.  Then $H(d_0)\subseteq
\Delta(K)$ but $H(d_0)\nless K$. We may assume further that

(1) (Maximality) For every missing face $T$ of $K$ of dimension $\leq d_0-1$, $H(d_0)< K\cup \{T\}$.

(2) (Links) For every two distinct vertices $v,u\in K_0$ such that $\{v,u\}$ is not contained in any
missing face of $K$ of dimension $\leq d_0-1$, $H(d_0-1)< \rm{lk}_K(u)\cap \rm{lk}_K(v)$.

(1) follows from the fact that if $K\subseteq L$ then
$\Delta(K)\subseteq \Delta(L)$. (2) follows from the minimality and
from Proposition \ref{propAlgContruction}. Indeed, if
$H(d_0-1)\nless \rm{lk}_K(u)\cap \rm{lk}_K(v)$ then
$H(d_0-1)\nsubseteq \Delta(\rm{lk}_K(u)\cap \rm{lk}_K(v))$ and as
$H(d_0)\subseteq \Delta(K)$ we obtain that the contraction of $v,u$
results in $K'$ for which $H(d_0)\subseteq \Delta(K')$ and
$H(d_0)\nless K'$, contradicting the minimality of $K$.

This led us to suspect that counterexamples for $d=3$ may be
provided by the following complexes. Let $M_L$ be the vertex
transitive neighborly $4$-sphere on $15$ vertices
$\rm{manifold}\_(4,15,5,1)$ found by Frank Lutz \cite{Lutz}. Note
that every edge in $M_L$ is contained in a missing triangle. Let $K$
be the $2$-skeleton of $M_L$ union with a missing triangle. It is
easy to find triangles such that every edge in $K$ is contained in a
missing triangle. As $M_L$ is neighborly, by counting and the fact
that $\Delta(K)$ is shifted, we conclude that
$\{5,6,7\}\in\Delta(K)$, hence $H(3)\subseteq\Delta(K)$. Is
$H(d)<K$? Note that deletions must be performed before any
contraction is possible.

\subsection{Proof of Theorem \ref{thm o(L)}}
The idea is to define an injective chain map $\phi:
C_*(H;\mathbb{Z}_2)\longrightarrow C_*(K;\mathbb{Z}_2)$ which
induces $\phi(\rm{Sm}^m(1_{K_*}))=\rm{Sm}^m(1_{H_*})$ for every
$m\geq 0$.

\begin{lem}\label{lemAdContr->InjChainMap}
Let $K\mapsto K'$ be an admissible contraction. Then it induces an
injective chain map $\phi: C_*(K';\mathbb{Z}_2)\longrightarrow
C_*(K;\mathbb{Z}_2)$.
\end{lem}
$Proof$: Fix a labeling of the vertices of $K$, $v_0,v_1,..,v_n$,
such that $K'$ is obtained from $K$ by identifying $v_0\mapsto v_1$
where $\rm{dim}(\rm{lk}(v_0,K))\leq \rm{dim}(\rm{lk}(v_1,K))$.

Let $F\in K'$. If $F\in K$, define $\phi(F)=F$. If $F\notin K$,
define $\phi(F)=\sum\{(F\setminus v)\cup v_0: v\in F, (F\setminus
v)\cup v_0\in K\}$. Note that if $F\notin K$ then $v_1\in F$ and
$(F\setminus v_1)\cup v_0 \in K$, so the sum above is nonzero.
Extend linearly to obtain a map $\phi:
C_*(K';\mathbb{Z}_2)\longrightarrow C_*(K;\mathbb{Z}_2)$.

First, let us check that $\phi$ is a chain map, i.e. that it
commutes with the boundary maps $\partial$. It is enough to verify
this for the basis elements $F$ where $F\in K'$. If $F\in K$ then
$\supp(\partial F)\subseteq K$, hence $\partial(\phi F)=\partial
F=\phi(\partial F)$. If $F\notin K$ then $\partial(\phi
F)=\partial(\sum\{(F\setminus v)\cup v_0: v\in F, (F\setminus
v)\cup v_0\in K\})$, and as we work over $\mathbb{Z}_2$, this
equals
\begin{equation}\label{eqPartialPhi}
\partial(\phi F) =
\sum\{F\setminus v : v\in F, (F\setminus v)\cup v_0\in K\} +
\end{equation}
\begin{equation}\nonumber
\sum\{(F\setminus \{u,v\})\cup v_0: u,v\in F, (F\setminus v)\cup
v_0\in K, (F\setminus u)\cup v_0\notin K\}.
\end{equation}
On the other hand $\phi(\partial F)=\phi(\sum\{F\setminus u: u\in
F, F\setminus u\in K\}) + \phi(\sum\{F\setminus u: u\in F,
F\setminus u\notin K\})$ and as we work over $\mathbb{Z}_2$, this
equals
\begin{equation}\label{eqPhiPartial}
\phi(\partial F) = \sum\{F\setminus u : u\in F, (F\setminus u)\in
K\} +
\end{equation}
\begin{equation}\nonumber
\sum\{(F\setminus \{u,v\})\cup v_0: u,v\in F, (F\setminus
\{u,v\})\cup v_0\in K, (F\setminus v)\in K,
 (F\setminus u)\notin K\}.
\end{equation}
It suffices to show that in equations (\ref{eqPartialPhi}) and
(\ref{eqPhiPartial}) the left summands on the RHSs are equal, as
well as the right summands on the RHSs. This follows from
observation \ref{obsAd} below. Thus $\phi$ is a chain map.

Second, let us check that $\phi$ is injective. Let $\pi_K$ be the
restriction map $C_*(K';\mathbb{Z}_2)\longrightarrow
\oplus\{\mathbb{Z}_2F: F\in K'\cap K\}$, $\pi_K(\sum\{\alpha_FF:
F\in K'\})=\sum\{\alpha_FF: F\in K'\cap K\}$. Similarly, let
$\pi_K^{\perp}$ be the restriction map
$C_*(K';\mathbb{Z}_2)\longrightarrow \oplus\{\mathbb{Z}_2F: F\in
K'\setminus K\}$. Note that for a chain $c\in C_*(K';\mathbb{Z}_2)$,
$c=\pi_K(c)+\pi_K^{\perp}(c)$ and $\supp(\phi(\pi_K(c)))\cap
\supp(\phi(\pi_K^{\perp}(c)))=\emptyset$. Assume that $c_1,c_2\in
C_*(K';\mathbb{Z}_2)$ such that $\phi(c_1)=\phi(c_2)$. Then
$\pi_K(c_1)=\phi(\pi_K(c_1))=\phi(\pi_K(c_2))=\pi_K(c_2)$, and
$\phi(\pi_K^{\perp}(c_1))=\phi(\pi_K^{\perp}(c_2))$. Note that if
$F_1,F_2 \notin K$ then $F_1,F_2\in K'$ and if $F_1\neq F_2$ then
$\supp(\phi(1F_1)) \ni (F_1\setminus v_1)\cup v_0 \notin
\supp(\phi(1F_2))$. Hence also
$\pi_K^{\perp}(c_1)=\pi_K^{\perp}(c_2)$. Thus $c_1=c_2$. $\square$

\begin{obs}\label{obsAd}
Let $K\mapsto K', v_0\mapsto v_1$ be an admissible contraction with
$\rm{dim}(\rm{lk}(v_0,K))\leq \rm{dim}(\rm{lk}(v_1,K))$. Let $K'\ni
F\notin K$ and $v\in F$. Then $(F\setminus v)\in K$ if and only if
$(F\setminus v)\cup v_0\in K$.
\end{obs}
$Proof$: Assume $F\setminus v\in K$. As $(F\setminus v_1)\cup v_0
\in K$ we only need to check the case $v\neq v_1$. We proceed by
induction on $\rm{dim}(F)$. As $\{v_0,v_1\}\in K$ whenever
$\rm{dim}(K)>0$ (and whenever $\rm{dim}(\rm{lk}(v_0,K))\geq 0$, if
we use the weaker local condition for admissible contractions), the
case $\rm{dim}(F)\leq 1$ is clear. (If $\rm{dim}(K)=0$ there is
nothing to prove. For the weaker local condition for admissible
contractions, if $\rm{lk}(v_0,K))=\emptyset$ then there is nothing
to prove.) By the induction hypothesis we may assume that all the
proper subsets of $(F\setminus v)\cup v_0$ are in $K$. Also
$v_0,v_1\in (F\setminus v)\cup v_0$. The admissibility of the
contraction implies that $(F\setminus v)\cup v_0\in K$. The other
direction is trivial. $\square$

\begin{lem}\label{lemInjChainMap->Smith}
Let $\phi: C_*(K';\mathbb{Z}_2)\longrightarrow C_*(K;\mathbb{Z}_2)$
be the injective chain map defined in the proof of Lemma
\ref{lemAdContr->InjChainMap} for an admissible contraction
$K\mapsto K'$. Then for every $m\geq 0$,
$\phi^*(\rm{Sm}^m([1_{K_*}]))=\rm{Sm}^m([1_{K'_*}])$ for the induced
map $\phi^*$.
\end{lem}
$Proof$: For two simplicial complexes $L$ and $L'$ and a field
$k$, the following map is an isomorphism of chain complexes:
$$\alpha=\alpha_{L,L',k}: C(L;k)\otimes_k C(L';k)\longrightarrow C(L*L';k), \ \ \alpha((1T)\otimes (1T'))=1(T\uplus T')$$
where $T\in L, T'\in L'$ and $\alpha$ is extended linearly. In
case $L=L'$ (in the definition of join we think of $L$ and $L'$ as
two disjoint copies of $L$) and $k$ is understood we denote
$\alpha_{L,L',k}=\alpha_{L}$.

Thus there is an induced chain map $\phi_*:
C_*(K'*K';\mathbb{Z}_2)\longrightarrow C_*(K*K;\mathbb{Z}_2)$,
$\phi_*=\alpha_{K}\circ \phi \otimes \phi \circ \alpha_{K'}^{-1}$
where $\phi \otimes \phi: C(K';\mathbb{Z}_2)\otimes_{\mathbb{Z}_2}
C(K';\mathbb{Z}_2)\longrightarrow
C(K;\mathbb{Z}_2)\otimes_{\mathbb{Z}_2} C(K;\mathbb{Z}_2)$ is
defined by $\phi \otimes \phi(c\otimes c')= \phi(c)\otimes \phi(c')$
(which this is a chain map).

Consider the subcomplex $C_*(K'_*;\mathbb{Z}_2)\subseteq
C_*(K'*K';\mathbb{Z}_2)$. We now verify that every $c\in
C_*(K'_*;\mathbb{Z}_2)$ satisfies $\phi_*(c)\in
C_*(K_*;\mathbb{Z}_2)$. It is enough to check this for chains of
the form $c=1(S^1\cup T^2)$ where $S,T\in K'$ and $S\cap T=
\emptyset$. For a collection of sets $A$ let $V(A)=\cup_{a\in
A}a$. Clearly if the condition
\begin{equation}\label{deljoinCond}
V(\supp (\phi(S)))\cap V(\supp (\phi(T)))= \emptyset
\end{equation}
is satisfied then we are done. If $v_1\notin S, v_1\notin T$, then
$\phi(S)=S, \phi(T)=T$ and (\ref{deljoinCond}) holds. If  $T\ni
v_1 \notin S$, then $\phi(S)=S$ and $V(\supp \phi(T))\subseteq
T\cup\{v_0\}$. As $v_0\notin S$ condition (\ref{deljoinCond})
holds. By symmetry, (\ref{deljoinCond}) holds when $S\ni v_1
\notin T$ as well.

With abuse of notation (which we will repeat) we denote the above
chain map by $\phi$, $\phi: C_*(K'_*;\mathbb{Z}_2)\longrightarrow
C_*(K_*;\mathbb{Z}_2)$. For a simplicial complex $L$, the involution
$\tau_L:L_*\longrightarrow L_*$, $\tau_L(S^1\cup T^2)=T^1\cup S^2$
induces a $\mathbb{Z}_2$-action on $C_*(L_*;\mathbb{Z}_2)$. It is
immediate to check that $\alpha_{L,L',k}$ and $\phi \otimes \phi$
commute with these $\mathbb{Z}_2$-actions, and hence so does their
composition, $\phi$. Thus, we have proved that $\phi:
C_*(K'_*;\mathbb{Z}_2)\longrightarrow C_*(K_*;\mathbb{Z}_2)$ is a
$\mathbb{Z}_2$-chain map.

Therefore, there is an induced map on the symmetric cohomology
rings $\phi: H_S^*(K_*)\longrightarrow H_S^*(K'_*)$ which commutes
with the connecting homomorphisms $\rm{Sm}:H_S^i(L)\longrightarrow
H_S^{i+1}(L)$ for $L=K_*,K'_*$.

Let us check that for the fundamental $0$-cocycles
$\phi([1_{K_*}])=[1_{K'_*}]$ holds. A representing cochain is
$1_{K_*}: \oplus_{v\in (K_*)_0}\mathbb{Z}_2 v \longrightarrow
\mathbb{Z}_2$, $1_{K_*}(1v)=1$. As $\phi|_{C_0(K'_*)}=id$ (w.r.t.
the obvious injection $(K'_*)_0\longrightarrow (K_*)_0$), for every
$u\in (K'_*)_0$ $(\phi
1_{K_*})(u)=1_{K_*}(\phi|_{C_0(K'_*)}(u))=1_{K_*}(u)=1$, thus
$\phi(1_{K_*})=1_{K'_*}$.

As $\phi$ commutes with the Smith connecting homomorphisms, for
every $m\geq 0$, $\phi(\rm{Sm}^m(1_{K_*}))=\rm{Sm}^m(1_{K'_*})$.
$\square$

\begin{thm}\label{thmMinor->Smith}
Let $H$ and $K$ be simplicial complexes. If $H<K$ then there
exists an injective chain map
$\phi:C_*(H;\mathbb{Z}_2)\longrightarrow C_*(K;\mathbb{Z}_2)$
which induces $\phi(\rm{Sm}^m(1_{K_*}))=\rm{Sm}^m(1_{H_*})$ for
every $m\geq 0$.
\end{thm}
$Proof$: Let the sequence $K=K^0\mapsto K^1\mapsto ...\mapsto K^t=H$
demonstrate the fact that $H<K$. If $K^i\mapsto K^{i+1}$ is an
admissible contraction, then by Lemmas \ref{lemAdContr->InjChainMap}
and \ref{lemInjChainMap->Smith} it induces an injective chain map
$\phi_i: C_*(K^{i+1};\mathbb{Z}_2)\longrightarrow
C_*(K^i;\mathbb{Z}_2)$ which in turn induces
$\phi_i(\rm{Sm}^m(1_{(K^{i})_*}))=\rm{Sm}^m(1_{(K^{i+1})_*})$ for
every $m\geq 0$. If $K^i\mapsto K^{i+1}$ is a deletion - take
$\phi_i$ to be the map induced by inclusion, to obtain the same
conclusions. Thus, the composition
$\phi=\phi_0\circ...\circ\phi_{t-1}:
C_*(H;\mathbb{Z}_2)\longrightarrow C_*(K;\mathbb{Z}_2)$ is as
desired. $\square$
\\
\\
\emph{Proof of Theorem \ref{thm o(L)}}: By Theorem
\ref{thmMinor->Smith}
$\phi(\rm{Sm}^{m}(1_{K_*}))=\rm{Sm}^{m}(1_{H_*})$. Thus if
$\rm{Sm}^{m}(1_{H_*})\neq 0$ then $\rm{Sm}^{m}(1_{K_*})\neq 0$.
 $\square$
\\
\textbf{Remark}: The conclusion of Theorem \ref{thm o(L)} would fail
if we allow arbitrary identifications of vertices. For example, let
$K'=K_5$ and let $K$ be obtained from $K'$ by splitting a vertex
$w\in K'$ into two new vertices $u,v$, and connecting $u$ to a
non-empty proper subset of $\rm{skel}_0(K')\setminus \{w\}$, denoted
by $A$, and connecting $v$ to $(\rm{skel}_0(K')\setminus
\{w\})\setminus A$. As $K$ embeds into the $2$-sphere,
$\rm{Sm}^3(K)=0$. By identifying $u$ with $v$ we obtain $K'$, but
$\rm{Sm}^3(K')\neq 0$. To obtain from this example an example where
the edge $\{u,v\}$ is present, let $L=\Cone(K)\cup \{u,v\}$, and let $L'$ be the complex
obtained form $L$ by identifying $u$ with $v$. Then $\rm{Sm}^4(L)=0$
while $\rm{Sm}^4(L')\neq 0$.

\begin{ex}\label{exFindMinor} Let $K$ be the simplicial complex
spanned by the following collection of $2$-simplices:
$(\binom{[7]}{3}\setminus \{127,137,237 \})\cup
\{128,138,238,178,278,378 \}$.
\end{ex}
$K$ is not a subdivision of $H(3)$, and its geometric realization
even does not contain a subspace homeomorphic to $H(3)$ (as there
are no $7$ points in $||K||$, each with a neighborhood whose
boundary contains a subspace which is homeomorphic to $K_6$).
Nevertheless, contraction of the edge $78$ is admissible and results
in $H(3)$. By Theorem \ref{thm o(L)} $K$ has a non-vanishing Van
Kampen's obstruction in dimension $5$, and hence is not embeddable
in the $4$-sphere.

\begin{ex}\label{exUli}
Let $T$ be a missing $d$-face in the cyclic $(2d+1)$-polytope on $n$
vertices, denoted by $C(2d+1,n)$, and let
$K=(C(2d+1,n))_d\cup \{T\}$. Then $\Delta(K)\nsubseteq
\Delta(2d+1,n)$ and $K$ is not embeddable in the $2d$-sphere.
\end{ex}
\emph{Proof - sketch}: As $\Delta(K)$ is shifted, by counting the number
of faces not greater or equal $\{d+3,...,2d+3\}$ in the product
partial order on $(d+1)$-tuples of $[n]$ we get $\Delta(K) \ni
\{d+3,...,2d+3\}\notin \Delta(2d+1,n)$.

By Gale evenness condition (see e.g.
\cite{Grunbaum},\cite{Ziegler}), $T$ is of the form
$\{t_1,...,t_{d+1}\}$ where $1<t_1, t_{d+1}<n$ and $t_i+1<t_{i+1}$
for every $i$. Let $s_i=t_i+1$ for $1\leq i\leq d$, $s_0=1$ and
$s_{d+1}=n$. In $K$ we can admissibly contract $j+1\mapsto j$ for
$j,j+1 \in [s_{i-1},t_i-1]$ for $1\leq i\leq d$ and similarly
$j\mapsto j+1$ for $j,j+1 \in [t_{d+1},s_{d+1}]$; to get a
simplicial complex $K'$ with the same description as $K$ but on one
vertex less, i.e. that its missing faces are the $(d+1)$-tuples of
pairwise non-adjacent vertices bigger than the smallest vertex and
smaller than the largest vertex except for one such set - $T$.
Successive application of these contractions results in $H(d+1)$ as
a minor of $K$, hence by Theorems \ref{thm o(L)}, \ref{thmSmithH(d)}
and \ref{thmSmith->nonEmb}, $K$ is not embeddable in the
$2d$-sphere. $\square$

Example \ref{exUli} is a special case of the following conjecture,
a work in progress of Uli Wagner and the author.

\begin{conj}\label{conjUli}
Let $K$ be a triangulated $2d$-sphere and let $T$ be a missing
$d$-face in $K$. Let $L= K_d\cup\{T\}$. Then $L$ does
not embed in $\mathbb{R}^{2d}$.
\end{conj}

\subsection{Proof of Theorem \ref{thm o_Z}}
Fix a total order on the
vertices of $K$, $v_0<v_1<..<v_n$ and consider an admissible
contraction $K\mapsto K'$ where $K'$ is obtained from $K$ by
identifying $v_0\mapsto v_1$ (shortly this will be shown to be
without loss of generality). Define a map $\phi$ as follows: for
$F\in K'$
\begin{equation}\label{phiTimes}
\phi(F)=\{^{F \ \ \rm{if}\ F\in K}_ {\sum\{\rm{sgn}(v,F)(F\setminus
v)\cup v_0:\ v\in F, (F\setminus v)\cup v_0\in K\}\ \ \rm{if}\
F\notin K}
\end{equation}
where $\rm{sgn}(v,F)=(-1)^{|\{t\in F: t<v\}|}$. Extend linearly to
obtain an injective $\mathbb{Z}$-chain map $\phi:
C_{\bullet}(K')\longrightarrow C_{\bullet}(K)$. (The check that this
map is indeed an injective $\mathbb{Z}$-chain map is similar to the
proof of Lemma \ref{lemAdContr->InjChainMap}.) In case we contract a
general $a\mapsto b$, for the signs to work out consider the map
$\tilde{\phi}=\pi^{-1}\phi\pi$ rather than $\phi$, where $\pi$ is
induced by a permutation on the vertices which maps $\pi(a)=v_0,\
\pi(b)=v_1$. Then $\tilde{\phi}$ is an injective $\mathbb{Z}$-chain
map.

As $\phi(S\times T):=\phi(S)\times \phi(T)$ commutes with the
$\mathbb{Z}_2$ action and with the boundary map on the chain complex
of the deleted product, $\phi$ induces a map
$H_{\rm{eq}}^{\bullet}(K_{\times})\rightarrow
H_{\rm{eq}}^{\bullet}(K'_{\times})$. It satisfies
$\phi^*(o^m_{\mathbb{Z}}(K_{\times}))=o^m_{\mathbb{Z}}(K'_{\times})$
for all $m\geq 1$. The checks are straightforward (for proving the
last statement, choose a total order with contraction which
identifies the minimal two elements $v_0\mapsto v_1$, and show
equality on the level of cochains). We omit the details.

If $K\mapsto K'$ is a deletion, consider the injection $\phi:
K'\rightarrow K$ to obtain again an induced map with
$\phi^*(o^m_{\mathbb{Z}}(K_{\times}))=o^m_{\mathbb{Z}}(K'_{\times})$.

Let the sequence $K=K^0\mapsto K^1\mapsto ...\mapsto K^t=H$
demonstrate the fact that $H<K$. By composing the corresponding maps
as above we obtain a map $\phi^*$ with
$\phi^*(o^m_{\mathbb{Z}}(K_{\times}))=o^m_{\mathbb{Z}}(H_{\times})$
and the result follows. $\square$

\section{Topology preserving edge contractions}\label{sec:TopPreservingEdgeContractions}
\subsection{PL manifolds}
The following theorem answers in the affirmative a question asked by
Dey et. al. \cite{Dey}, who already proved the dimension $\leq 3$
case.
\begin{thm}\label{thmDey}
Given an edge in a triangulation of a compact PL (piecewise
linear)-manifold without boundary, its contraction results in a
PL-homeomorphic space if and only if it satisfies the Link Condition
(\ref{eqLinkCond}).
\end{thm}
$Proof$: Let $M$ be a PL-triangulation
of a compact $d$-manifold without boundary. Let $ab$ be an edge of
$M$ and let $M'$ be obtained from $M$ by contracting $a\mapsto b$.
We will prove that if the Link Condition (\ref{eqLinkCond}) holds
for $ab$ then $M$ and $M'$ are PL-homeomorphic, and otherwise they
are not homeomorphic (not even 'locally homologic'). For $d=1$ the
assertion is clear. Assume $d>1$.

Denote $B(b)=\{b\}*\antist(b,\lk(a,M))$ and
$L=\antist(a,M)\cap B(b)$. Then $M'=\antist(a,M)\cup_{L}B(b)$.
As $M$ is a PL-manifold without boundary, $\lk(a,M)$ is a
$(d-1)$-PL-sphere (see e.g. \cite{Hudson}, Corollary 1.16). By
Newman's theorem (e.g. \cite{Hudson}, Theorem 1.26)
$\antist(b,\lk(a,M))$ is a $(d-1)$-PL-ball. Thus $B(b)$ is a
$d$-PL-ball. Observe that $\partial
(B(b))=\antist(b,\lk(a,M))\cup \{b\} *
\lk(b,\lk(a,M)) = \lk(a,M) = \partial (\clst(a,M))$.

The identity map on $\rm{lk}(a,M)$ is a PL-homeomorphism $h:
\partial (B(b))\rightarrow \partial (\clst(a,M))$, hence it
extends to a PL-homeomorphism $\tilde{h}: B(b)\rightarrow
\clst(a,M)$ (see e.g. \cite{Hudson}, Lemma 1.21).

Note that $L=\rm{lk}(a,M)\cup (\{b\}*(\rm{lk}(a,M)\cap
\rm{lk}(b,M)))$.

If $\rm{lk}(a)\cap \rm{lk}(b)=\rm{lk}(ab)$ (in $M$) then
$L=\rm{lk}(a,M)$, hence gluing together the maps $\tilde{h}$ and
the identity map on $\rm{ast}(a,M)$ results in a PL-homeomorphism
from $M'$ to $M$.

If $\rm{lk}(a)\cap \rm{lk}(b)\neq \rm{lk}(ab)$ (in $M$) then
$\rm{lk}(a,M)\subsetneqq L$. The case $L=B(b)$ implies that
$M'=\rm{ast}(a,M)$ and hence $M'$ has a nonempty boundary, showing
it is not homeomorphic to $M$. A small punctured neighborhood of a
point in the boundary of $M'$ has trivial homology while all small
punctured neighborhoods of points in $M$ has non vanishing
$(d-1)$-th homology. This is what we mean by 'not even locally
homologic': $M$ and $M'$ have homologically different sets of
small punctured neighborhoods.

We are left to deal with the case $\rm{lk}(a,M)\subsetneqq
L\subsetneqq B(b)$. As $L$ is closed there exists a point $t\in
L\cap \rm{int}(B(b))$ with a small punctured neighborhood $N(t,M')$
which is not contained in $L$. For a subspace $K$ of $M'$ denote by
$N(t,K)$ the neighborhood in $K$ $N(t,M')\cap K$. Thus
$N(t,M')=N(t,\rm{ast}(a,M))\cup_{N(t,L)}N(t,B(b))$. We get a
Mayer-Vietoris exact sequence in reduced homology:
\begin{equation}\label{eqMV}
H_{d-1}N(t,L)\rightarrow H_{d-1}N(t,\rm{ast}(a,M))\oplus
H_{d-1}N(t,B(b))\rightarrow H_{d-1}N(t,M')\rightarrow
\end{equation}
\begin{equation}\nonumber
H_{d-2}N(t,L)\rightarrow H_{d-2}N(t,\rm{ast}(a,M))\oplus
H_{d-2}N(t,B(b)).
\end{equation}
Note that $N(t,\rm{ast}(a,M))$ and $N(t,B(b))$ are homotopic to
their boundaries which are $(d-1)$-spheres. Note further that
$N(t,L)$ is homotopic to a proper subset $X$ of $\partial
(N(t,B(b)))$ such that the pair $(\partial (N(t,B(b))),X)$ is
triangulated. By Alexander duality $H_{d-1}N(t,L)=0$. Thus,
(\ref{eqMV}) simplifies to the exact sequence
$$0\rightarrow \mathbb{Z}\oplus \mathbb{Z}\rightarrow H_{d-1}N(t,M')\rightarrow
H_{d-2}N(t,L)\rightarrow 0.$$ Thus, $\rm{rank}(H_{d-1}N(t,M'))\geq
2$, hence $M$ and $M'$ are not locally homologic, and in
particular are not homeomorphic. $\square$
\\
\textbf{Remarks}: (1) Omitting the assumption in Theorem
\ref{thmDey} that the boundary is empty makes both implications
incorrect. Contracting an edge to a point shows that the Link
Condition is not sufficient. Contracting an edge on the boundary
of a cone over an empty triangle shows that the Link Condition is
not necessary.

(2) The necessity of the Link Condition holds also in the
topological category (and not only in the PL category), as the
proof of Theorem \ref{thmDey} shows. Indeed, for this part we only
used the fact that $B(b)$ is a pseudo manifold with boundary
$\rm{lk}(a,M)$ (not that it is a ball); taking the point $t$ to
belong to exactly two facets of $B(b)$. For sufficiency of the Link Condition in the
topological category, see Problem \ref{prob:Top} in Section \ref{sec:top open prob} below.

Walkup \cite{Walkup} mentioned, without details, the necessity of
the Link Condition for contractions in topological manifolds, as
well as the sufficiency of the Link Condition for the $3$
dimensional case (where the category of PL-manifolds coincides
with the topological one); see \cite{Walkup}, p.82-83.

\subsection{PL spheres}\label{subsecSpheres}
\begin{de}\label{defStrongContr}
Boundary complexes of simplices are \emph{strongly edge
decomposable} and, recursively, a triangulated PL-manifold $S$ is
\emph{strongly edge decomposable} if it has an edge which satisfies
the Link Condition (\ref{eqLinkCond}) such that both its link and
its contraction are strongly edge decomposable.
\end{de}

By Theorem \ref{thmDey} the complexes in Definition
\ref{defStrongContr} are all triangulated PL-spheres. Note that
every $2$-sphere is strongly edge decomposable.

Let $vu$ be an edge in a simplicial complex $K$ which satisfies the
Link Condition, whose contraction $u\mapsto v$ results in the
simplicial complex $K'$. Note that the $f$-polynomials satisfy
$$f(K,t)=f(K',t)+t(1+t)f(\rm{lk}(\{vu\},K),t),$$
hence the $h$-polynomials satisfy
\begin{equation}\label{eq-hContraction}
h(K,t)=h(K',t)+t h(\rm{lk}(\{vu\},K),t).
\end{equation}
We conclude the following:
\begin{cor}
The $g$-vector of strongly edge decomposable triangulated spheres is
non negative. $\square$
\end{cor}
Is it also an $M$-vector? Compare with Theorem
\ref{thm:Stellar(HL)=HL}. The strongly edge decomposable spheres
(strictly) include the family of triangulated spheres which can be
obtained from the boundary of a simplex by repeated Stellar
subdivisions (at any face); the later are polytopal, hence their
$g$-vector is an $M$-sequence. For the case of subdividing only at
edges (\ref{eq-hContraction}) was considered by Gal (\cite{Gal},
Proposition 2.4.3).

\section{Open problems}\label{sec:top open prob}
\begin{enumerate}

\item Prove that if $H<K$ and $K$ is embeddable in the $m$-sphere then $H$ is
embeddable in the $m$-sphere.

\item Let $K$ be a triangulated $2d$-sphere and let $T$ be a missing
$d$-face in $K$. Let $L= K_d\cup\{T\}$. Show that $L$ does not embed
in $\mathbb{R}^{2d}$.

\item\label{prob:Top}
Given an edge in a triangulation of a compact manifold without
boundary which satisfies the Link Condition, is it true that its
contraction results in a homeomorphic space? Or at least in a
space of the same homotopic or homological type?

A Mayer-Vietoris argument shows that such topological manifolds $M$
and $M'$ have the same Betti numbers; both $\clst(a,M)$ and $B(b)$
are cones and hence their reduced homology vanishes.

A candidate for a counterexample for Problem \ref{prob:Top} may be
the join $M=T*P$ where $T$ is the boundary of a triangle and $P$ a
triangulation of Poincar\'{e} homology $3$-sphere, where an edge
with one vertex in $T$ and the other in $P$ satisfies the Link
Condition. By the double-suspension theorem (Edwards
\cite{Edwards} and Cannon \cite{Cannon}) $M$ is a topological
$5$-sphere.

\item Show that the $g$-vector of strongly edge decomposable triangulated spheres is an $M$-vector.
\end{enumerate}


\chapter{Face Rings for Graded Posets}\label{chapter:ShiftingPosets}

\section{Some classic $f$-vector results for graded posets}
Let us review the characterization of $f$-vectors of finite
simplicial complexes, known as the
Sch\"{u}tzenberger-Kruskal-Katona theorem (see \cite{Bollobas} for a
proof and for references). For any two integers $k,n>0$ there exists a unique
expansion
\begin{equation} \label{nk}
n= {n_{k}\choose k} + {n_{k-1}\choose k-1}+...+ {n_{i}\choose i}
\end{equation}
such that $n_{k}>n_{k-1}>...>n_{i}\geq i \geq 1$ (details in
\cite{Bollobas}). Define the function $\partial_{k-1}$ by
 $$\partial_{k-1}(n)=
{n_{k}\choose k-1} + {n_{k-1}\choose k-2}+...+ {n_{i}\choose
i-1},\ \ \partial _{k-1}(0)=0.$$
\begin{thm}[Sch\"{u}tzenberger-Kruskal-Katona] \label{KKthm}
$f$ is the $f$-vector of some
simplicial complex iff $f$ ultimately vanishes and
\begin{equation} \label{partial_{k-1}(n)}
\forall k\geq 0\ \   0\leq \partial_{k}(f_{k})\leq f_{k-1}.
\end{equation}
\end{thm}
For a ranked meet semi-lattice $P$, finite at every rank, let
$f_{i}$ be the number of elements with rank $i+1$ in $P$, and set
${\rm rank}(\hat{0})=0$ where $\hat{0}$ is the minimum of $P$. The
$f$-vector of $P$ is $(f_{-1},f_0,f_1,...)$.

$P$ has the \emph{diamond property} if for every $x,y\in P$ such that
$x<y$ and ${\rm rank} (y)-{\rm rank} (x)=2$ there exist at least
two elements in the open interval $(x,y)$. The closed interval is
denoted by $[x,y]=\{z\in P: x\leq z\leq y \}$.

We identify a simplicial complex with the poset of its faces
ordered by inclusion. The following generalization of Theorem
\ref{KKthm} is due to Wegner \cite{Wegner}.

\begin{thm}[Wegner]\label{Wthm}
Let $P$ be a finite ranked meet semi-lattice with the diamond
property. Then its $f$-vector ultimately vanishes and satisfies
(\ref{partial_{k-1}(n)}).
\end{thm}
For $\hat{x}\in P$ define $P(\hat{x})=\{x\in P: \hat{x}\leq x\}$
and let $y'\prec  y$ denote $y$ covers $y'$.

\begin{lem}\label{cond* lemma}
For a ranked meet semi-lattice $P$, the diamond property is
equivalent to satisfying the following condition:

(\textbf{*}) For every $\hat{x}\in P$, $x$ which covers $\hat{x}$
and $y$ such that $y \in P(\hat{x})$ and $y\not =\hat{x}$, there
exists $y'\in P(\hat{x})$ such that $y'\prec y$ and $x\nleq y'$.
\end{lem}
A multicomplex (on a finite ground set) can be considered as an
order ideal of monomials $I$ (i.e. if $m|n \in I$ then also $m\in
I$) on a finite set of variables. Its $f$-vector is defined by
$f_{i}=|\{m\in I: {\rm deg} (m)=i+1\}|$ (again $f_{-1}=1$). Define
the function  $\partial^{k-1}$ by $$\partial^{k-1}(n)=
{n_{k}-1\choose k-1} + {n_{k-1}-1\choose k-2}+...+ {n_{i}-1\choose
i-1},\ \
\partial^{k-1}(0)=0,$$ w.r.t the expansion (\ref{nk}).
\begin{thm}[Macaulay](simpler proofs in \cite{CL,Stanley-Hilbert})\label{Mthm}
 $f$ is
the $f$-vector of some multicomplex iff $f_{-1}=1$ and
\begin{equation} \label{partial^{k-1}(n)}
\forall k\geq 0\ \    0\leq \partial^{k}(f_{k})\leq f_{k-1}.
\end{equation}
\end{thm}

\begin{figure}\label{Fig1}
\newcommand{\edge}[1]{\ar@{-}[#1]}
\newcommand{\lulab}[1]{\ar@{}[l]^<<{#1}}
\newcommand{\rulab}[1]{\ar@{}[r]^<<{#1}}
\newcommand{\ldlab}[1]{\ar@{}[l]^<<{#1}}
\newcommand{\rdlab}[1]{\ar@{}[r]_<<{#1}}
\newcommand{\node}{*+[O][F-]{ }}
\centerline{ \xymatrix{
 & & & & \circ \rulab{y} \edge{d} \edge{ddlll} & & & &\circ \rulab{y} \edge{d} \edge{dlll} &\\
& \circ \lulab{x_{i+1}} \edge{d}\ & & & \bullet \rulab{\exists y'} \edge{ddlll} & \circ \lulab{x_{r}} \edge{d} & & & \bullet \rulab{\exists y'} \edge{dlll} &\\
& \circ \lulab{x_{i}} \edge{d} & & & & \circ \lulab{x_{r-1}} \edge{dd}\\
& \circ \lulab{x_{i-1}} \edge{d}\\
& \circ \lulab{\hat{x}} & & & & \circ \lulab{\hat{x}} } } \caption
{The parallelogram property for $i<r$ (left) and for $i=r$
(right).}
\end{figure}
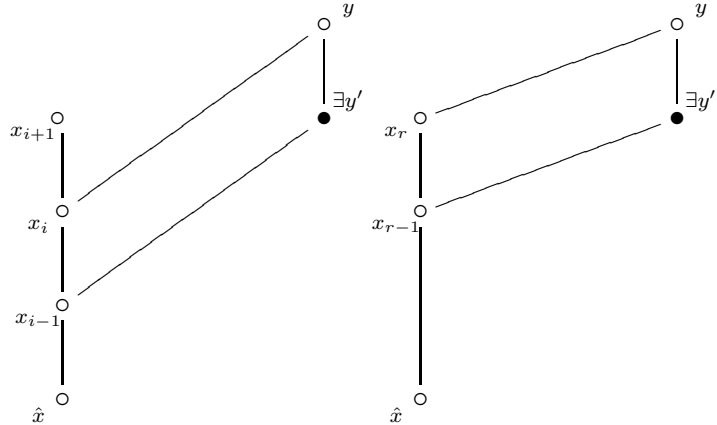

\begin{de}(Parallelogram property)\label{def**}
A ranked poset $P$ is said to have the \emph{parallelogram property}
if the following condition holds:

(\textbf{*}\textbf{*}) For every $\hat{x}\in P$ and $y\in
P(\hat{x})$ such that $y\not =\hat{x}$, if the chain
$\{\hat{x}=x_{0}\prec x_{1}\prec ...\prec x_{r} \}$ equals the
closed interval $[\hat{x},x_{r}]$ ($r>0$) and is maximal w.r.t.
inclusion such that $r<{\rm rank}(y)$ (the rank of $y$ in the
poset $P(\hat{x})$), and if $x_{i}<y$ and $x_{i+1}\nleqslant y$ for
some $0< i\leq r$, then there exists $y'\in P(\hat{x})$ such that
$y'\prec y$, $x_{i-1}<y'$ and $x_{i}\nleq y'$. For $i=r$ interpret
$x_{r+1}\nleqslant y$ as: $[\hat{x},y]$ is not a chain.
\end{de}
See Figure \ref{Fig1} for an illustration of the parallelogram
property. Note that condition (\textbf{*}) of Lemma \ref{cond*
lemma} implies condition (\textbf{*}\textbf{*}) of Definition
\ref{def**} (with $1$ being the only possible value of $r$).
Posets of multicomplexes, polyhedral complexes, and rooted trees, satisfy the parallelogram property.
\\
We identify a multicomplex with the poset of its monomials ordered
by division. We now generalize Theorem \ref{Mthm}; the proof is combinatorial.
\begin{thm}[\cite{Nevo-GeneralizedMacaulay}, Theorem 1.6] \label{THM}
Let $P$ be a ranked meet semi-lattice, finite at every rank, with
the parallelogram property. Then its $f$-vector satisfies
(\ref{partial^{k-1}(n)}) and $f_{-1}(P)=1$.
\end{thm}

For generalizations of Macaulay's theorem in a different direction
('compression'), see e.g. \cite{CL,WW}.

In Section \ref{sec:BasicsShifting} we presented the symmetric and
exterior face rings of a simplicial complex, to which we applied
the shifting operator. The graded components of the face ring have
dimensions corresponding to the $f$-vector. Shifting changes the
bases of these components, hence preserves the $f$-vector, and
results in a shifted complex, for which e.g. the inequalities $\partial_{k}(f_{k})\leq f_{k-1}$ are easier to prove.

Analogous algebraic object and operator for more general graded
posets are desirable in order to prove $f$-vector theorems for them.
An algebraic
object that will correspond to Macaulay inequalities, may help to
settle the following well known conjecture:
\begin{conj}\label{M-polytopes}
The toric $g$-vector of a (non-simplicial) polytope is an
$M$-sequence, i.e. it satisfies Macaulay inequalities
(\ref{partial^{k-1}(n)}).
\end{conj}
Recently Karu \cite{Karu} proved that the toric $g$-vector of a
polytope is nonnegative, using a (complicated) graded module. Can we
shift this structure, in order to prove Conjecture
\ref{M-polytopes}? A problem is that Karu's structure is a module
over the polynomial ring with $d$ (=dimension) variables, and not
with $n$ (=number of vertices) variables.

In the next two sections we make initial steps in this program.

\section{Algebraic shifting for geometric meet semi-lattices}\label{sec:ShiftGeomLattice}
We will associate an analogue of the exterior face ring to
geometric ranked meet semi-lattices, which coincides with the
usual construction for the case of simplicial complexes. Applying
algebraic shifting we
construct a canonically defined shifted simplicial complex, having
the same $f$-vector as its geometric meet semi-lattice.

Let $(L,<,r)$ be a ranked atomic meet semi-lattice with $L$ the set of its
elements, $<$ the partial order relation and $r:L\rightarrow
\mathbb{N}$ its rank function. We denote it in short by $L$. $L$ is called $geometric$ if
\begin{equation}\label{rank-ineq}
r(x\wedge y)+r(x\vee y)\leq r(x)+r(y)
\end{equation}
for every $x,y\in L$ such that $x\vee y$ exists. For example, the
intersections of a finite collection of hyperplanes in a vector
space form a geometric meet semi-lattice w.r.t. the reverse
inclusion order and the codimension rank. Face posets of simplicial
complexes are important examples of geometric meet semi-lattices,
where (\ref{rank-ineq}) holds with equality.

Adding a maximum to a ranked meet semi-lattice makes it a lattice,
denoted by $\hat{L}$, but the maximum may not have a rank. Denote by
$\hat{0},\hat{1}$ the minimum and maximum of $\hat{L}$,
respectively, and by $L_i$ the set of rank $i$ elements in $L$.
$r(\hat{0})=0$.

We now define the algebra $\bigwedge L$ over a field $k$ with
characteristic $2$. Let $V$ be a vector space over $k$ with basis
$\{e_u: u\in L_1\}$. Let $I_L=I_1+I_2+I_3$ be the ideal in the
exterior algebra $\bigwedge V$ defined as follows. Choose a total
ordering of $L_1$, and denote by $e_S$ the wedge product
$e_{s_1}\mathbf{\wedge}...\mathbf{\wedge} e_{s_{|S|}}$ where
$S=\{s_1<...<s_{|S|}\}$. Define:
\begin{equation}\label{I1-eq}
I_1=(e_S: S\subseteq L_1, \vee S=\hat{1}\in \hat{L}),
\end{equation}
\begin{equation}\label{I2-eq}
I_2=(e_S: S\subseteq L_1, \vee S\in L, r(\vee S)\neq |S|),
\end{equation}
\begin{equation}\label{I3-eq}
I_3=(e_S-e_T: T,S\subseteq L_1, \vee T=\vee S\in L, r(\vee
S)=|S|=|T|, S\neq T).
\end{equation}

(As ${\rm char} (k)=2$, $e_S-e_T$ is independent of the ordering
of the elements in $S$ and in $T$.) Let $\bigwedge L = \bigwedge
V/I_L$. As $I_L$ is generated by homogeneous elements, $\bigwedge
L$ inherits a grading from $\bigwedge V$. Let $f(\bigwedge
L)=(f_{-1},f_0,..)$ be its graded dimensions vector, i.e.
$f_{i-1}$ is the dimension of the degree $i$ component of
$\bigwedge L$.
\newline \textbf{Remark}:
If $L$ is the poset of a simplicial complex, then $I_L=I_1$ and $\bigwedge L$
is the classic exterior face ring of $L$, as in \cite{55}.

The following proposition will be used for showing that $\bigwedge
L$ and $L$ have the same $f$-vector. Its easy proof by induction
on the rank is omitted.
\begin{prop}\label{rank}
Let $L$ be a geometric ranked meet semi-lattice. Let $l\in L$ and
let $S$ be a minimal set of atoms such that $\vee S=l$, i.e. if
$T\subsetneq S$ then $\vee T<l$. Then $r(l)=|S|$. $\square$
\end{prop}
\textbf{Remark}: The converse of Proposition \ref{rank} is also
true: Let $L$ be a ranked atomic meet semi-lattice such that every
$l\in L$ and every minimal set of atoms $S$ such that $\vee S=l$
satisfy $r(l)=|S|$. Then $L$ is geometric.

\begin{prop}\label{basis}
$f(\bigwedge L)=f(L)$.
\end{prop}
$Proof$:
Denote by $\tilde{w}$ the
projection of $w\in \bigwedge V$ on $\bigwedge L$.
We will show that
picking $S(l)$ such that $S(l)\subseteq L_1, \vee S(l)=l, |S(l)|=r(l)$
for each $l\in L$ gives a basis over $k$ of $\bigwedge L$, $E=\{\tilde{e}_{S(l)}: l\in L\}$.

As $\{\tilde{e}_{S}: S\subseteq L_1\}$ is a basis of $\bigwedge
V$, it is clear from the definition of $I_L$ that $E$ spans
$\bigwedge L$. To show that $E$ is independent, we will prove
first that the generators of $I_{L}$ as an ideal, that are
specified in (\ref{I2-eq}), (\ref{I1-eq}) and (\ref{I3-eq}),
actually span it as a vector space over $k$.

As $x\vee \hat{1}=\hat{1}$ for all $x\in L$, the generators of
$I_1$ that are specified in (\ref{I1-eq}) span it as a $k$-vector
space. Next, we show that the generators of $I_2$ and $I_1$ that
are specified in (\ref{I2-eq}) and in (\ref{I1-eq}) respectively,
span $I_1+I_2$ as a $k$-vector space: if $e_S$ is such a generator
of $I_2$ and $U\subseteq L_1$ then either $e_U\mathbf{\wedge}
e_S\in I_1$ (if $U\cap S\neq \emptyset$ or if $\vee(U\cup
S)=\hat{1}$) or else, by Proposition \ref{rank}, $r(\vee(U\cup
S))<|U\cup S|$ and hence $e_U\mathbf{\wedge} e_S$ is also such a
generator of $I_2$.

Let $e_S-e_T$ be a generator of $I_3$ as specified in
(\ref{I3-eq}) and let $U\subseteq L_1$. If $U\cap T\neq \emptyset$
then $e_T\mathbf{\wedge} e_U=0$ and $e_S\mathbf{\wedge} e_U$ is
either zero (if $U\cap S\neq \emptyset$) or else a generator of
$I_1+I_2$, by Proposition \ref{rank}; and similarly when $U\cap
S\neq \emptyset$. If $U\cap T= \emptyset = U\cap S$ then $\vee
(S\cup U)=\vee (T\cup U)$ and $|S\cup U|=|T\cup U|$. Hence, if
$e_S\mathbf{\wedge} e_U-e_T\mathbf{\wedge} e_U$ is not the obvious
difference of two generators of $I_1$ or of $I_2$ as specified in
(\ref{I1-eq}) and (\ref{I2-eq}), then it is a generator of $I_3$
as specified in (\ref{I3-eq}). We conclude that these generators
of $I_{L}$ as an ideal span it as a vector space over $k$.

Assume that $\sum_{l\in L}a_l \tilde{e}_{S(l)}=0$, i.e.
$\sum_{l\in L}a_l e_{S(l)}\in I_L$ where $a_l\in k$ for all $l\in
L$. By the discussion above, $\sum_{l\in L}a_l e_{S(l)}$ is in the
span (over $k$) of the generators of $I_3$ that are specified in
(\ref{I3-eq}). But for every $l\in L$ and every such generator $g$
of $I_3$, if $g=\sum\{b_S e_S: \vee S\in L, r(\vee S)=|S|\}$
($b_S\in k$ for all $S$) then $\sum\{b_S: \vee S=l\}=0$. Hence
$a_l=0$ for every $l\in L$. Thus $E$ is a basis of $\bigwedge L$,
hence $f(\bigwedge L)=f(L)$. $\square$

Now let us shift. Note that exterior algebraic shifting, which was
defined for the exterior face ring, can be applied to any graded
exterior algebra finitely generated by degree $1$ elements. It
results in a simplicial complex with an $f$-vector that is equal to
the vector of graded dimensions of the algebra. This shows that any
such graded algebra satisfies Kruskal-Katona inequalities! We apply
this construction to $\bigwedge L$:

Let $B=\{b_u: u\in L_1\}$ be a basis of $V$. Then $\{\tilde{b}_S:
S\subseteq L_1\}$ spans $\bigwedge L$. Choosing a basis from this
set in the greedy way w.r.t. the lexicographic order $<_L$ on equal
sized sets, defines a collection of sets:
$$\Delta_B(L)=\{S: \tilde{b}_S\notin {\rm span}_k\{\tilde{b}_T: |T|=|S|, T<_L S\}\}.$$
$\Delta_B(L)$ is a simplicial complex, and by Proposition
\ref{basis} $f(\Delta_B(L))=f(L)$. For a generic $B$,
$\Delta_B(L)$ is shifted.
Moreover, the construction is canonical, i.e. is independent both
of the chosen ordering of $L_1$ and of the generically chosen
basis $B$. It is also independent of the characteristic $2$ field
that we picked. We denote $\Delta(L)=\Delta_B(L)$ for a generic
$B$. For proofs of the above statements we refer to Bj\"{o}rner
and Kalai \cite{BK} (they proved for the case where $L$ is a
simplicial complex, but the proofs remain valid for any graded
exterior algebra finitely generated by degree $1$ elements).

We summarize the above discussion in the following theorem:
\begin{thm}\label{GL}
Let $L$ be a geometric meet semi-lattice, and let $k$ be a field
of characteristic $2$. There exists a canonically defined shifted
simplicial complex $\Delta(L)$ associated with $L$, with
$f(\Delta(L))=f(L)$. $\square$
\end{thm}
\textbf{Remarks}: (1) The fact that $L$ satisfies Kruskal-Katona
inequalities follows also without using our algebraic construction,
from the fact that it satisfies the diamond property and applying
Theorem \ref{Wthm}. The diamond property is easily seen to hold for
all ranked atomic meet semi-lattices.

(2) A different operation, which does depend on the ordering of
$L_1$ and results in a simplicial complex with the same $f$-vector,
was described by Bj\"{o}rner \cite{Bj-MatroidApplications}, Chapter
7, Problem 7.25: totally order $L_1$. For each $x\in L$ choose the
lexicographically least subset $S_x\subseteq L_1$ such that $\vee
S_x=x$ ($S_{\hat{0}}=\emptyset$). Define $\Delta_<(L)=\{S_x: x\in
L\}$. Then $\Delta_<(L)$ is a simplicial complex with the same
$f$-vector as $L$. An advantage in our operation is that it is
canonical (and results in a shifted simplicial complex). To see that
these two operations are indeed different, let $L$ be the face poset
of a simplicial complex. Then for any total ordering of $L_1$,
$\Delta_<(L)=L$. But if the simplicial complex is not shifted (e.g.
a $4$-cycle), then $\Delta(L)\neq L$.

\section{Algebraic shifting for generalized multicomplexes}
We will associate an analogue of the symmetric (Stanley-Reisner)
face ring with a common generalization of multicomplexes and
geometric meet semi-lattices. Applying an algebraic shifting
operation, we construct a multicomplex having the same $f$-vector
as the original poset.

Let $\mathbb{P}$ be the following family of posets: to construct
$P\in \mathbb{P}$ start with a geometric meet semi-lattice $L$.
Associate with each $l\in L$ the (square free) monomial
$m(l)=\prod_{a<l,a\in L_1}x_a$, and equip it with rank
$r(m(l))=r(l)$. Denote this collection of monomials by $M_0$. Now
repeat the following procedure finitely or countably many times to
construct $(M_0\subseteq M_1\subseteq...)$: Choose $m\in M_i$ and
$a\in L$ such that $x_a|m$, $\frac{x_a}{x_b}m\in M_i$ for all
$b\in L_1$ such that $x_b|m$, and $x_am\notin M_i$. $M_{i+1}$ is
obtained from $M_i$ by adding $x_am$, setting its rank to be
$r(x_am)=r(m)+1$ and let it cover all the elements
$\frac{x_a}{x_b}m$ where $b\in L_1$ such that $x_b|m$. Define
$P=\cup M_i$.

Note that the posets in $\mathbb{P}$ are ranked (not necessarily
atomic) meet semi-lattices with the parallelogram property, and
that $\mathbb{P}$ includes all multicomplexes (start with $L$, a
simplicial complex) and geometric meet semi-lattices ($P=M_0$).

For $P\in \mathbb{P}$ define the following analogue of the
Stanley-Reisner ring: Assume for a moment that $P$ is finite. Fix
a field $k$, and denote $P_1=\{1,..,n\}$. Let $A=k[x_1,..,x_n]$ be
a polynomial ring. For $j$ such that $1\leq j\leq n$ let $r_j$ be
the minimal integer number such that $x_j^{r_j+1}$ does not divide
any of the monomials $p\in P$. Note that each $i\in P$ of rank $1$
belongs to a unique maximal interval which is a chain; whose top
element is $x_i^{r_i}$. By abuse of notation, we identify the
elements in such intervals with their corresponding monomials in
$A$.

We add a maximum $\hat{1}$ to $P$ to obtain $\hat{P}$ and define the following ideals in $A$:

$I_0=(\prod_{i=1}^n x_i^{a_i}: \exists j\ 1\leq j\leq n,\
a_j>r_j),$

$I_1=(\prod_{i=1}^n x_i^{a_i}: \forall j\ a_j\leq r_j, \ \vee_{i=1}^n x_i^{a_i}=\hat{1}\in \hat{P}),$

$I_2=(\prod_{i=1}^n x_i^{a_i}: \vee_{i=1}^n x_i^{a_i}\in P, r(\vee_{i=1}^n x_i^{a_i})\neq\sum_i a_i),$

$I_3=(\prod_{i=1}^n x_i^{a_i}-\prod_{i=1}^n x_i^{b_i}: \vee_{i=1}^n x_i^{a_i}=\vee_{i=1}^n x_i^{b_i}\in P, r(\vee_{i=1}^n x_i^{a_i})=\sum_i a_i=\sum_i b_i),$

$I_P=I_0+I_1+I_2+I_3.$

Define $k[P]:=A/I_P$. As $I_P$ is homogeneous, $k[P]$ inherits a
grading from $A$. Let $f(k[P])=(f_{-1},f_0,..)$ where $f_i=\dim_k
\{m\in k[P]: r(m)=i+1 \}$ ($f_{-1}=1$).

The proof of the following proposition is similar to the proof of Proposition \ref{basis}, and is omitted.
\begin{prop}\label{basis2}
$f(k[P])=f(P)$. $\square$
\end{prop}
Denote by $\tilde{w}$ the projection of $w\in A$ on $k[P]$. Let
$B=\{y_1,..,y_n\}$ be a basis of $A_1$. Then
$$\Delta_B(P):=\{\prod_{i=1}^n y_i^{a_i}:
\prod_{i=1}^n\tilde{y_i}^{a_i} \notin {\rm
span}_k\{\prod_{i=1}^n\tilde{y_i}^{b_i}: \sum_{i=1}^n
a_i=\sum_{i=1}^n b_i, \prod_{i=1}^n y_i^{b_i}<_L \prod_{i=1}^n
y_i^{a_i} \}\}$$ is an order ideal of monomials with an $f$-vector
$f(P)$. (The lexicographic order on monomials of equal degree is
defined by $\prod_{i=1}^n y_i^{b_i}<_L \prod_{i=1}^n y_i^{a_i}$
iff there exists $j$ such that for all $1\leq t<j\ a_t=b_t$ and
$b_j>a_j$.) To prove this, we reproduce the argument of Stanley
for proving Macaulay's theorem (\cite{St}, Theorem 2.1): as the
projections of the elements in $\Delta_B(P)$ form a $k$-basis of
$k[P]$, then by Proposition \ref{basis2} $f(\Delta_B(P))=f(P)$. If
$m\notin \Delta_B(P)$ then $m=\sum\{a_n n: {\rm deg} (n)={\rm deg}
(m), n<_L m\}$, hence for any monomial $m'$ $m'm=\sum\{a_n m'n:
{\rm deg} (n)={\rm deg} (m), n<_L m\}$. But ${\rm deg} (m'm)={\rm
deg} (m'n)$ and $m'n<_L m'm$ for these $n$'s, hence $m'm\notin
\Delta_B(P)$, thus $\Delta_B(P)$ is an order ideal of monomials.
\newline \textbf{Remark}: For $B$ a generic basis the construction is canonical in
 the same sense as defined for the exterior case.

Combining Proposition \ref{basis2} with Theorem \ref{Mthm} we obtain
\begin{cor}\label{MP}
Every $P\in \mathbb{P}$ satisfies Macaulay inequalities
(\ref{partial^{k-1}(n)}). $\square$
\end{cor}
\textbf{Remark}:
If $P$ is infinite, let $P_{\leq r}:=\{p\in P: r(p)\leq r\}$ and construct $\Delta(P_{\leq r})$ for each $r$.
Then $\Delta(P_{\leq r})\subseteq \Delta(P_{\leq r+1})$ for every $r$, and $\Delta(P):=\cup_r\Delta(P_{\leq r})$
is an order ideal of monomials with $f$-vector $f(P)$. Hence, Corollary \ref{MP} holds in this case too.

\section{Open problems}
\begin{enumerate}
\item The Kruskal-Katona inequalities hold for any meet semi
lattice with the diamond property (Theorem \ref{Wthm}). Can an algebraic
proof be given?

\item Similarly, can Theorem \ref{THM} be proved algebraically?

\item Prove Conjecture \ref{M-polytopes} by applying shifting to a suitable algebraic object.
\end{enumerate}


\bibliographystyle{plain}
\bibliography{biblio}
\end{document}